\magnification1200
\baselineskip15pt
\newread\AUX\immediate\openin\AUX=\jobname.aux
\def\ref#1{\expandafter\edef\csname#1\endcsname}
\ifeof\AUX\immediate\write16{\jobname.aux gibt es nicht!}\else
\input \jobname.aux
\fi\immediate\closein\AUX
\def\today{\number\day.~\ifcase\month\or
  Januar\or Februar\or M{\"a}rz\or April\or Mai\or Juni\or
  Juli\or August\or September\or Oktober\or November\or Dezember\fi
  \space\number\year}
\font\sevenex=cmex7
\scriptfont3=\sevenex
\font\fiveex=cmex10 scaled 500
\scriptscriptfont3=\fiveex

\def\hidew#1{\setbox0=\hbox to 0pt{\hss$#1$\hss}\box0}
\def\nospace{\thickmuskip0mu\thinmuskip0mu\medmuskip0mu}
\def\a{\alpha}
\def\e{{\bf e}}
\def\rho{\varrho}

\def\phi{\varphi}
\def\epsilon{\varepsilon}
\def\theta{\vartheta}
\def\uauf{\lower1.7pt\hbox to 3pt{%
\vbox{\offinterlineskip
\hbox{\vbox to 8.5pt{\leaders\vrule width0.2pt\vfill}%
\kern-.3pt\hbox{\lams\char"76}\kern-0.3pt%
$\raise1pt\hbox{\lams\char"76}$}}\hfil}}
\def\cite#1{\expandafter\ifx\csname#1\endcsname\relax
{\bf?}\immediate\write16{#1 ist nicht definiert!}\else\csname#1\endcsname\fi}
\def\expandwrite#1#2{\edef\next{\write#1{#2}}\next}
\def\neverexpand{\noexpand\noexpand\noexpand}
\def\strip#1\ {}
\def\ncite#1{\expandafter\ifx\csname#1\endcsname\relax
{\bf?}\immediate\write16{#1 ist nicht definiert!}\else
\expandafter\expandafter\expandafter\strip\csname#1\endcsname\fi}
\newwrite\AUX
\immediate\openout\AUX=\jobname.aux
\newcount\Abschnitt\Abschnitt0
\def\beginsection#1. #2 \par{\advance\Abschnitt1\GNo0%
\vskip0pt plus.10\vsize\penalty-250
\vskip0pt plus-.10\vsize\bigskip\vskip\parskip
\edef\TEST{\number\Abschnitt}
\expandafter\ifx\csname#1\endcsname\TEST\relax\else
\immediate\write16{#1 hat sich geaendert!}\fi
\expandwrite\AUX{\neverexpand\ref{#1}{\TEST}}
\leftline{\bf\number\Abschnitt. \ignorespaces#2}%
\nobreak\smallskip\noindent\SATZ1}
\def\Proof:{\par\noindent{\it Proof:}}
\def\Remark:{\ifdim\lastskip<\medskipamount\removelastskip\medskip\fi
\noindent{\bf Remark:}}
\def\Remarks:{\ifdim\lastskip<\medskipamount\removelastskip\medskip\fi
\noindent{\bf Remarks:}}
\def\Definition:{\ifdim\lastskip<\medskipamount\removelastskip\medskip\fi
\noindent{\bf Definition:}}
\def\Example:{\ifdim\lastskip<\medskipamount\removelastskip\medskip\fi
\noindent{\bf Example:}}
\def\Examples:{\ifdim\lastskip<\medskipamount\removelastskip\medskip\fi
\noindent{\bf Examples:}}
\newcount\SATZ\SATZ1
\def\proclaim #1. #2\par{\ifdim\lastskip<\medskipamount\removelastskip
\medskip\fi
\noindent{\bf#1.\ }{\it#2}\Par
\ifdim\lastskip<\medskipamount\removelastskip\goodbreak\medskip\fi}
\def\Aussage#1{%
\expandafter\def\csname#1\endcsname##1.{\ifx?##1?\relax\else
\edef\TEST{#1\penalty10000\ \number\Abschnitt.\number\SATZ}
\expandafter\ifx\csname##1\endcsname\TEST\relax\else
\immediate\write16{##1 hat sich geaendert!}\fi
\expandwrite\AUX{\neverexpand\ref{##1}{\TEST}}\fi
\proclaim {\number\Abschnitt.\number\SATZ. #1\global\advance\SATZ1}.}}
\Aussage{Theorem}
\Aussage{Proposition}
\Aussage{Corollary}
\Aussage{Lemma}
\font\la=lasy10
\def\strich{\hbox{$\vcenter{\hbox
to 1pt{\leaders\hrule height -0,2pt depth 0,6pt\hfil}}$}}
\def\dashedrightarrow{\hbox{%
\hbox to 0,5cm{\leaders\hbox to 2pt{\hfil\strich\hfil}\hfil}%
\kern-2pt\hbox{\la\char\string"29}}}

\def\Bindestrich{\penalty10000-\hskip0pt}
\let\_=\Bindestrich
\def\.{{\sfcode`.=1000.}}

\def\Par{\par}
\def\:={\mathrel{\raise0,9pt\hbox{.}\kern-2,77779pt
\raise3pt\hbox{.}\kern-2,5pt=}}
\def\=:{\mathrel{=\kern-2,5pt\raise0,9pt\hbox{.}\kern-2,77779pt
\raise3pt\hbox{.}}} 
\def\into{\hookrightarrow}
\def\pfeil{\rightarrow}

\def\pf#1{\buildrel#1\over\rightarrow}

\def\Ugleich{\hbox{$\cup$\kern.5pt\vrule depth -0.5pt}}
\def\|#1|{\mathop{\rm#1}\nolimits}
\def\<{\langle}
\def\>{\rangle}
\let\Times=\times
\def\times{\mathop{\Times}}
\let\Otimes=\otimes
\def\otimes{\mathop{\Otimes}}
\catcode`\@=11
\def\hex#1{\ifcase#1 0\or1\or2\or3\or4\or5\or6\or7\or8\or9\or A\or B\or
C\or D\or E\or F\else\message{Warnung: Setze hex#1=0}0\fi}
\def\fontdef#1:#2,#3,#4.{%
\alloc@8\fam\chardef\sixt@@n\FAM
\ifx!#2!\else\expandafter\font\csname text#1\endcsname=#2
\textfont\the\FAM=\csname text#1\endcsname\fi
\ifx!#3!\else\expandafter\font\csname script#1\endcsname=#3
\scriptfont\the\FAM=\csname script#1\endcsname\fi
\ifx!#4!\else\expandafter\font\csname scriptscript#1\endcsname=#4
\scriptscriptfont\the\FAM=\csname scriptscript#1\endcsname\fi
\expandafter\edef\csname #1\endcsname{\fam\the\FAM\csname text#1\endcsname}
\expandafter\edef\csname hex#1fam\endcsname{\hex\FAM}}
\catcode`\@=12 

\fontdef Ss:cmss10,,.
\fontdef Fr:eufm10,eufm7,eufm5.
\def\fa{{\Fr a}}
\def\fb{{\Fr b}}


\def\ft{{\Fr t}}
\def\fu{{\Fr u}}

\def\fX{{\Fr X}}
\def\fY{{\Fr Y}}

\fontdef bbb:msbm10,msbm7,msbm5.
\fontdef mbf:cmmib10,cmmib7,.
\fontdef msa:msam10,msam7,msam5.
\def\CC{{\bbb C}}\def\DD{{\bbb D}}

\def\NN{{\bbb N}}\def\PP{{\bbb P}}
\def\QQ{{\bbb Q}}\def\RR{{\bbb R}}

\def\ZZ{{\bbb Z}}
\def\cA{{\cal A}}\def\cC{{\cal C}}

\def\cP{{\cal P}}
\def\cR{{\cal R}}

\mathchardef\leer=\string"0\hexbbbfam3F
\mathchardef\subsetneq=\string"3\hexbbbfam24
\mathchardef\semidir=\string"2\hexbbbfam6E
\mathchardef\dirsemi=\string"2\hexbbbfam6F
\mathchardef\haken=\string"2\hexmsafam78
\mathchardef\auf=\string"3\hexmsafam10
\let\OL=\overline
\def\overline#1{{\hskip1pt\OL{\hskip-1pt#1\hskip-1pt}\hskip1pt}}
\def\aq{{\overline{a}}}

\def\cq{{\overline{c}}}


\def\Gq{{\overline{G}}}

\def\pq{{\overline{p}}}\def\Pq{{\overline{P}}}

\def\Rq{{\overline{R}}}

\def\uq{{\overline{u}}}
\def\vq{{\overline{v}}}

\def\Xq{{\overline{X}}}
\def\Yq{{\overline{Y}}}
\def\zq{{\overline{z}}}
%
\abovedisplayskip 9.0pt plus 3.0pt minus 3.0pt
\belowdisplayskip 9.0pt plus 3.0pt minus 3.0pt
\newdimen\Grenze\Grenze2\parindent\advance\Grenze1em
\newdimen\Breite
\newbox\DpBox
\def\NewDisplay#1$${\Breite\hsize\advance\Breite-\hangindent
\setbox\DpBox=\hbox{\hskip2\parindent$\displaystyle{#1}$}%
\ifnum\predisplaysize<\Grenze\abovedisplayskip\abovedisplayshortskip
\belowdisplayskip\belowdisplayshortskip\fi
\global\futurelet\nexttok\WEITER}
\def\WEITER{\ifx\nexttok\qed\expandafter\leftQEDdisplay
\else\leftdisplay\fi}
\def\leftdisplay{\hskip-\hangindent\leftline{\box\DpBox}$$}
\def\leftQEDdisplay{\hskip-\hangindent
\line{\copy\DpBox\hfill\lower\dp\DpBox\copy\QEDbox}%
\belowdisplayskip0pt$$\bigskip\let\nexttok=}
\everydisplay{\NewDisplay}
\newbox\QEDbox
\newbox\nichts\setbox\nichts=\vbox{}\wd\nichts=2mm\ht\nichts=2mm
\setbox\QEDbox=\hbox{\vrule\vbox{\hrule\copy\nichts\hrule}\vrule}
\def\qed{\leavevmode\unskip\hfil\null\nobreak\hfill\copy\QEDbox\medbreak}
\newdimen\HIindent
\newbox\HIbox
\def\setHI#1{\setbox\HIbox=\hbox{#1}\HIindent=\wd\HIbox}
\def\HI#1{\par\hangindent\HIindent\hangafter=0\noindent\leavevmode
\llap{\hbox to\HIindent{#1\hfil}}\ignorespaces}

\newdimen\maxSpalbr
\newdimen\altSpalbr
\newcount\Zaehler

\def\beginrefs{%
\expandafter\ifx\csname Spaltenbreite\endcsname\relax
\def\Spaltenbreite{1cm}\immediate\write16{Spaltenbreite undefiniert!}\fi
\expandafter\altSpalbr\Spaltenbreite
\maxSpalbr0pt
\gdef\alt{}
\def\\##1\relax{%
\gdef\neu{##1}\ifx\alt\neu\global\advance\Zaehler1\else
\xdef\alt{\neu}\global\Zaehler=1\fi\xdef\SigText{##1\the\Zaehler}}
\def\L|Abk:##1|Sig:##2|Au:##3|Tit:##4|Zs:##5|Bd:##6|S:##7|J:##8||{%
\def\SigText{##2}\global\setbox0=\hbox{##2\relax}
\edef\TEST{[\SigText]}
\expandafter\ifx\csname##1\endcsname\TEST\relax\else
\immediate\write16{##1 hat sich geaendert!}\fi
\expandwrite\AUX{\neverexpand\ref{##1}{\TEST}}
\setHI{[\SigText]\ }
\ifnum\HIindent>\maxSpalbr\maxSpalbr\HIindent\fi
\ifnum\HIindent<\altSpalbr\HIindent\altSpalbr\fi
\HI{[\SigText]}
\ifx-##3\relax\else{##3}: \fi
\ifx-##4\relax\else{##4}{\sfcode`.=3000.} \fi
\ifx-##5\relax\else{\it ##5\/} \fi
\ifx-##6\relax\else{\bf ##6} \fi
\ifx-##8\relax\else({##8})\fi
\ifx-##7\relax\else, {##7}\fi\Par}
\def\B|Abk:##1|Sig:##2|Au:##3|Tit:##4|Reihe:##5|Verlag:##6|Ort:##7|J:##8||{%
\def\SigText{##2}\global\setbox0=\hbox{##2\relax}
\edef\TEST{[\SigText]}
\expandafter\ifx\csname##1\endcsname\TEST\relax\else
\immediate\write16{##1 hat sich geaendert!}\fi
\expandwrite\AUX{\neverexpand\ref{##1}{\TEST}}
\setHI{[\SigText]\ }
\ifnum\HIindent>\maxSpalbr\maxSpalbr\HIindent\fi
\ifnum\HIindent<\altSpalbr\HIindent\altSpalbr\fi
\HI{[\SigText]}
\ifx-##3\relax\else{##3}: \fi
\ifx-##4\relax\else{##4}{\sfcode`.=3000.} \fi
\ifx-##5\relax\else{(##5)} \fi
\ifx-##7\relax\else{##7:} \fi
\ifx-##6\relax\else{##6}\fi
\ifx-##8\relax\else{ ##8}\fi\Par}
\def\Pr|Abk:##1|Sig:##2|Au:##3|Artikel:##4|Titel:##5|Hgr:##6|Reihe:{%
\def\SigText{##2}\global\setbox0=\hbox{##2\relax}
\edef\TEST{[\SigText]}
\expandafter\ifx\csname##1\endcsname\TEST\relax\else
\immediate\write16{##1 hat sich geaendert!}\fi
\expandwrite\AUX{\neverexpand\ref{##1}{\TEST}}
\setHI{[\SigText]\ }
\ifnum\HIindent>\maxSpalbr\maxSpalbr\HIindent\fi
\ifnum\HIindent<\altSpalbr\HIindent\altSpalbr\fi
\HI{[\SigText]}
\ifx-##3\relax\else{##3}: \fi
\ifx-##4\relax\else{##4}{\sfcode`.=3000.} \fi
\ifx-##5\relax\else{In: \it ##5}. \fi
\ifx-##6\relax\else{(##6)} \fi\PrII}
\def\PrII##1|Bd:##2|Verlag:##3|Ort:##4|S:##5|J:##6||{%
\ifx-##1\relax\else{##1} \fi
\ifx-##2\relax\else{\bf ##2}, \fi
\ifx-##4\relax\else{##4:} \fi
\ifx-##3\relax\else{##3} \fi
\ifx-##6\relax\else{##6}\fi
\ifx-##5\relax\else{, ##5}\fi\Par}
\bgroup
\baselineskip12pt
\parskip2.5pt plus 1pt
\hyphenation{Hei-del-berg}
\sfcode`.=1000
\beginsection References. References

}
\def\endrefs{%
\expandwrite\AUX{\neverexpand\ref{Spaltenbreite}{\the\maxSpalbr}}
\ifnum\maxSpalbr=\altSpalbr\relax\else
\immediate\write16{Spaltenbreite hat sich geaendert!}\fi
\egroup}


\newcount\GNo\GNo=0
\def\eqno#1{\global\advance\GNo1%
\edef\FTEST{$(\number\Abschnitt.\number\GNo)$}%
\ifx?#1?\relax\else
\expandafter\ifx\csname#1\endcsname\FTEST\relax\else
\immediate\write16{#1 hat sich geaendert!}\fi
\expandwrite\AUX{\neverexpand\ref{#1}{\FTEST}}\fi
\llap{\hbox to 40pt{\FTEST\hfill}}}

\catcode`@=11
\def\eqalignno#1{\null\vcenter{\openup\jot\m@th\ialign{\eqno{##}\hfil
&\strut\hfil$\displaystyle{##}$&$\displaystyle{{}##}$\hfil\crcr#1\crcr}}\,}
\catcode`@=12

\def\noq{{\overline{n}}}
\def\nuq{{\underline{n}}}
\def\moq{{\overline{m}}}
\def\muq{{\underline{m}}}
\def\Mod{\|mod|}
\def\r{r}
\def\rr{\cdot2r}
\def\odd{{\rm odd}}
\def\even{{\rm even}}
\def\LE{\sqsubseteq}
\mathchardef\lee=\string"3\hexmsafam35
\mathchardef\Less=\string"3\hexmsafam40
\def\LEE{\mathrel{\vcenter{\hbox{$\prec$}\vskip-10pt\hbox{$=$}}}}

\font\BF=cmbx10 scaled \magstep2
\font\CSC=cmcsc10 

\baselineskip15pt

{\baselineskip1.5\baselineskip\rightskip0pt plus 5truecm
\leavevmode\vskip0truecm\noindent
\BF Semisymmetric polynomials and the invariant theory of
matrix vector pairs

}
\vskip1truecm
\leftline{{\CSC Friedrich Knop}%
\footnote*{\rm Partially supported by a grant of the NSF}}
\leftline{Department of Mathematics, Rutgers University, New Brunswick NJ
08903, USA}
\leftline{knop@math.rutgers.edu}
\vskip1truecm

\beginsection Introduction. Introduction

\noindent In this paper we investigate a new family of multivariable
polynomials. These polynomials, denoted $R_\lambda(z_1,\ldots,z_n;r)$,
depend on a parameter $r$ and are indexed by a partition $\lambda$ of
length $n$. Up to a scalar, $R_\lambda$ is characterized by the
following elementary properties:

\item{$\bullet$}$R_\lambda$ is symmetric in the odd variables
$z_1,z_3,z_5,\ldots$ as well as in the even variables
$z_2,z_4,z_6,\ldots$. Polynomials having this kind of symmetry are
called {\it semisymmetric}.

\item{$\bullet$}For the partition
$\lambda=(\lambda_1,\ldots,\lambda_n)$ define the odd degree as
$|\lambda|_\odd:=\sum_{i\ \odd}\lambda_i$. Then the degree of
$R_\lambda(z)$ is $|\lambda|_\odd$.

\item{$\bullet$}Consider the vector
$\rho:=((n-1)\r,(n-2)\r,\ldots,\r,0)$. Then $R_\lambda(z)$ vanishes at
all points of the form $z=\rho+\mu$ where $\mu$ is any partition with
$\mu\ne\lambda$ and $|\mu|_\odd\le|\lambda|_\odd$.

\noindent The simplest nontrivial example comes from the partition
$(1)=(1,0,\ldots,0)$ in which case
$R_\lambda(z;r)=\sum_{i=1}^n(-1)^{i-1}z_i-\lfloor n^2/4\rfloor$. It is
clearly semisymmetric, has degree $|\lambda|_\odd =1$ and vanishes at
$z=\rho+\mu$ where $\mu=(0)$ or $\mu=(1^2)$.

The polynomials $R_\lambda(z)$ are analogous to the polynomials
$P_\lambda(z;r)$ which were previously introduced in\cite{KnSa}. In
fact, the definition of $P_\lambda$ is the same except that
$P_\lambda$ is symmetric in {\it all\/} variables $z_1,\ldots,z_n$ and
the odd degree $|\lambda|_\odd$ is replaced by the (full) degree
$|\lambda|=\sum_i\lambda_i$. The $P_\lambda$ are called {\it
shifted Jack polynomials\/} since their highest degree components are
the Jack polynomials. This is in contrast to their semisymmetric
counterparts: even their highest degree components form a genuinely new
class of multivariable homogeneous polynomials.

All the polynomials mentioned above have a representation theoretic
origin. Let $G$ be a connected reductive group acting on a finite
dimensional vector space $V$. We are interested in the case when this
action is {\it multiplicity free,\/} i.e., every simple $G$\_module
occurs at most once in the algebra of polynomial functions of
$V$. Then the algebra $\PP\DD^G$ of $G$\_invariant differential
operators on $V$ is commutative. Moreover, one can define a (Harish
Chandra) isomorphism which identifies $\PP\DD^G$ with the space
$\CC[\fa^*]^{W_V}$ of $W_V$\_invariant polynomials on a finite
dimensional vector space $\fa^*$ where $W_V$ is a finite reflection
group.

The point is now that there exist very particular invariant
differential operators $D_\lambda$ on $V$ which form a basis of
$\PP\DD^G$. The idea of their construction goes back to Capelli. The
operators $D_\lambda$ correspond, via the Harish Chandra isomorphism,
to polynomials $p_\lambda\in\CC[\fa^*]^{W_V}$ and it is these
polynomials which one would like to understand. Already the top
homogeneous component $\pq_\lambda$ of $p_\lambda$ is very important
since it has the following representation theoretic meaning. Consider
the symbol $E_\lambda$ of the differential operator $D_\lambda$. By
construction, this is a $G$\_invariant function on the cotangent
bundle $T_V^*=V\oplus V^*$. It can be considered as a generalization
of a zonal spherical function. On the other hand, one can define a
(Chevalley) isomorphism of $\CC[V\oplus V^*]^G$ with
$\CC[\fa^*]^{W_V}$ and under this isomorphism $E_\lambda$ corresponds
to $\pq_\lambda$.

The investigation of the polynomials $p_\lambda$ and $\pq_\lambda$ is
greatly facilitated by the fact that multiplicity free actions on
vector spaces are classified. This is due to the efforts of Kac
\cite{Kac}, Benson-Ratcliff \cite{BenRat}, and Leahy \cite{Leahy}. The
most important numerical invariant of a multiplicity free action is
the dimension of $\fa^*$ which is called its {\it rank\/}. It follows
from the classification that there are only seven series in which the
rank is unbounded. These series are listed in the following
table. More precisely, an (indecomposable) multiplicity free action
which is not in the table has rank less or equal~$7$.

\medskip
\leftline{
\vbox{\halign{$#$\ \hfil&$(#)$\qquad\hfil&$#$\qquad\hfil&$#$\qquad\hfil&\hfil$#$\hfil\cr
G&\omit\hfill&V&\|rank|&r\cr
\noalign{\smallskip\hrule\medskip\noindent{\it Classical cases:}\hfill}
GL_p(\CC)&p\ge2&S^2(\CC^p)&p&{1\over2}\cr
GL_p(\CC)\times GL_q(\CC)&p,q\ge1&\CC^p\otimes\CC^q&\|min|(p,q)&1\cr
GL_p(\CC)&p\ge2&\Lambda^2(\CC^p)&\lfloor{p\over2}\rfloor&2\cr
\noalign{\noindent{\it Semiclassical cases:}\hfill}
GL_p(\CC)\times GL_q(\CC)&p,q\ge1&(\CC^p\otimes\CC^q)\oplus\CC^q&\|min|(2p+1,2q)&{1\over2}\cr
GL_p(\CC)&p\ge2&\Lambda^2(\CC^p)\oplus\CC^p&p&1\cr
\noalign{\noindent{\it Quasiclassical cases:}\hfill\smallskip}
GL_p(\CC)\times GL_q(\CC)&p,q\ge1&(\CC^p\otimes\CC^q)\oplus(\CC^q)^*&\|min|(2p+1,2q)&\cr
GL_p(\CC)&p\ge2&\Lambda^2(\CC^p)\oplus(\CC^p)^*&p&\cr
}}}

\medskip\noindent
As indicated in the table, the seven series fall into three classes:
classical, semiclassical, and quasiclassical. The reason for that is
that all the cases in each class can be treated uniformly: there are
polynomials depending on a free parameter\footnote*{In the odd rank
quasiclassical case there are two parameters.} $r$ such that
the polynomials $p_\lambda$ of each particular case are obtained by
specializing the parameter as indicated in the table.

In the classical class, the space $V$ consists of matrices: symmetric,
rectangular, or skewsymmetric. This case has been treated in
\cite{KnSa} and the polynomials $p_\lambda$ are basically the shifted
Jack polynomials $P_\lambda(z;r)$.

The purpose of the present paper is to study the semiclassical
case. Here an element of $V$ is a pair $(A,v)$ where $A$ is a
rectangular or skewsymmetric matrix and $v$ is a vector. The
polynomials $p_\lambda$ are our $R_\lambda(z;r)$ described in the
beginning.

The quasiclassical class is left mostly to future research. Here an
element of $V$ is a pair $(A,\alpha)$ where $A$ is a matrix as above
and $\alpha$ is a covector (linear form). Preliminary investigations
indicate that this case is much more involved than the other two
cases. Nevertheless, the small cases, more precisely the cases with
$\|rank|\le4$, are covered also in the present paper since for those
the combinatorics of the quasiclassical and the semiclassical class
turn out to be isomorphic. In particular, we can also say something
about the action of $GL_1(\CC)\times GL_q(\CC)$ and
$\CC^q\oplus(\CC^q)^*$, a case already considered by Vilenkin--\v
Sapiro~\cite{VS}.

The zonal spherical functions $E_\lambda$ have numerous different
descriptions. This means that the results of this paper are also
relevant for the action of $GL_{p-1}(\CC)$ on $GL_p(\CC)$ {\it by
conjugation} or for the action of $Sp_{2p}(\CC)$ on
$X:=SL_{2p+1}(\CC)/Sp_{2p}(\CC)$. This is remarkable since the space
$X$ is only spherical and not symmetric.

Now we describe the our results about the polynomials $R_\lambda$. The
most important result is the construction of $n$ commuting {difference
operators\/} of which the $R_\lambda$ are simultaneous
eigenfunctions. These differential operators are defined by
explicit formulas \cite{E90}. An analogous result has already been the
main statement of \cite{KnSa}. The result here is similar but much
more involved.

Except for some elementary results, like existence and uniqueness of
the $R_\lambda$, most proofs hinge on the difference operators. The first
immediate consequence is that the top homogeneous component
$\Rq_\lambda$ of $R_\lambda$ is a simultaneous eigenfunction of $n$
commuting {\it differential operators\/} of order
$1,1,2,2,3,3,\ldots$. These are semisymmetric analogues of the
Sekiguchi-Debiard operators which characterize Jack polynomials.
Observe that Heckman and Opdam define analogues of the
Sekiguchi-Debiard operators for any finite root system but our
semisymmetric case is {\it not\/} covered by their construction.

Another rather immediate consequence of the difference operators is
the {\it Extra Vanishing Theorem}. Remember, that $R_\lambda$ is {\it
defined\/} to vanish at all points of the form $z=\rho+\mu$ where
$\mu\ne\lambda$ and $|\mu|_\odd\le|\lambda|_\odd$. It turns out, that
$R_\lambda$ actually vanishes at many more points. In
section~\cite{DifferenceOperators} we define an order relation
$\lambda\LE\mu$ on the set of partitions such that
$R_\lambda(\rho+\mu)=0$ whenever $\lambda\not\LE\mu$. This order
relation should be regarded as a semisymmetric analogue of the
familiar containment relation for partitions.

A property which can be considered as dual to extra vanishing is
called {\it triangularity\/}. By {\it definition}, the polynomial
$R_\lambda$ can be expressed as a linear combination of monomials
$z^\mu$ whose degree is less or equal $|\lambda|_\odd$. As it turns
out much fewer monomials are needed. This phenomenon is called
triangularity since it can be rephrased as saying that the base change
matrix from monomials to $R_\lambda$'s is triangular. In
section~\cite{Triangularity} we actually prove two versions of
triangularity. For the first, we define a map
$\lambda\mapsto[\lambda]$ from the set of partitions into the set
$\NN^n$ of compositions such that if $z^\mu$ appears in $R_\lambda$
then $\mu\le[\lambda]$. Here ``$\le$'' is the usual (inhomogeneous)
dominance order on $\NN^n$.

The monomials $z^\mu$ are not semisymmetric. Therefore, one can
attempt to formulate triangularity strictly within the set of
semisymmetric polynomials. To do this, define the {\it elementary
semisymmetric polynomials\/} as $\e_1:=\Rq_{(1)}$,
$\e_2:=\Rq_{(1^2)}$, etc. These can be computed explicitly:
$\e_{2i-1}=e_i^\odd-e_i^\even$ and $\e_{2i}=e_i^\even$ where
$e_i^{\odd/\even}$ is the usual elementary symmetric function of
degree $i$ in the variables $\{z_{2j-1}\mid j\,\}$ or $\{z_{2j}\mid
j\,\}$, respectively. Now consider all monomials in the $\e_i$:
$$\eqno{}
\e_\mu:=\e_1^{\mu_1-\mu_2}\e_2^{\mu_2-\mu_3}
\ldots \e_{n-1}^{\mu_{n-1}-\mu_n}\e_n^{\mu_n}
$$
where $\mu$ is a partition. In section~\cite{Triangularity} we define
a new order relation $\mu\preceq\lambda$ on the set of partitions
which is a semisymmetric analogue of the classical dominance
order. Using the explicit form of the difference operators we are able
to prove that each $R_\lambda$ is a linear combinations of $\e_\mu$
with $\mu\preceq\lambda$. This is the second triangularity result,
alluded to above. It should be noted that the second form easily
implies the first one but not conversely. This has to be seen in
contrast to the classical case of (shifted) Jack polynomials where
both forms of triangularity are actually equivalent.

In section~\cite{Binomial}, we prove what could be considered as the
second main result of this paper: the {\it duality formula} \cite{E1}
$$\eqno{E91}
{R_\lambda(-\underline{\a}-z)\over R_\lambda(-\underline{\a}-\rho)}=
\sum_\mu(-1)^{|\mu|_\odd}
{R_\mu(\rho+\lambda)\over R_\mu(\rho+\mu)}
{R_\mu(z)\over R_\mu(-\underline{\a}-\rho)}.
$$
where $\a$ is an arbitrary parameter and
$\underline{\a}=(\a,\ldots,\a)$. In other words, the formula expresses
the transformation $z_i\mapsto-\a-z_i$ of the space of semisymmetric
polynomials in terms of its basis $R_\lambda$. The classical
analogue has been established Okounkov~\cite{Ok} whose proof we follow
closely. This holds also for some of its consequences described below.
The key to the proof of the duality formula are again the difference
operators.

A first consequence is an explicit interpolation formula
(\cite{Interpol}{\it iii)}). It allows to explicitly calculate the
expansion of an arbitrary semisymmetric polynomial in terms of
$R_\lambda$'s. More precisely,
$$\eqno{}
f(z)=\sum_\mu(-1)^{|\mu|_\odd}{\hat f(\rho+\mu)\over R_\mu(\rho+\mu)}\,R_\mu(z)
$$
where
$$\eqno{}
\hat f(\rho+\mu)=\sum_\nu(-1)^{|\nu|_\odd}
{R_\nu(\rho+\mu)\over R_\nu(\rho+\nu)}f(\rho+\nu).
$$
The key to this formula is the observation that the transformation
$z\mapsto -\underline\alpha-z$ is involutory.

Another consequence of the duality theorem is the fact that the
expression
$$\eqno{}
R_\lambda(-\underline\alpha-\rho-\nu)\over R_\lambda(-\underline\alpha-\rho)
$$
is symmetric in $\lambda$ and $\nu$ (just put $z=\rho+\nu$ in
\cite{E91}). This observation allows sometimes to interchange in
formulas the argument $z$ with the index $\lambda$. For example, let
$D$ be one of the fundamental difference operators. Then the
eigenvalue equation $D(R_\lambda)=c(\lambda)R_\lambda$ turns into a
{\it Pieri type formula}, i.e., a formula which expresses the
multiplication operator $p\mapsto fp\,$ on the space of semisymmetric
polynomials (with fixed $f$) in terms of the basis $R_\lambda$. This
way, we obtain in section~\cite{PieriFormula} the expansion of
$f(z)R_\lambda(z)$ in terms of $R_\mu$'s where $f(z)$ is one of
$e_i^\odd(z)$, $e_i^\even(z)$, or $R_{1^i}(z)$. As a byproduct of
these investigation we prove in section~\cite{EvaluationFormula} a
formula for the value of $R_\lambda(z)$ in $z=-\underline\a-\rho$.

Most of these results have consequences for the homogeneous
polynomials $\Rq_\lambda$. The evaluation formula specializes to a
formula for the value of $\Rq_\lambda(z)$ in $z=(1,\ldots,1)$. The
duality formula implies the semisymmetric binomial theorem which
expresses the effect of the transformation $z_i\mapsto \a+z_i$ in
terms of the homogeneous basis $\Rq_\lambda$ (see
section~\cite{EvaluationFormula} for these two statements). Its
classical analogue is due to Okounkov-Olshanski~\cite{OO2}. Finally,
we obtain expansions of $f(z)\Rq_\lambda$ in terms of $\Rq_\mu$'s
where $f(z)$ is one of $e_i^\odd(z)$, $e_i^\even(z)$, or $\e_i(z)$

Finally, scattered all over the paper, we derive several explicit
formulas. More precisely, we determine $R_\lambda(z)$ when
$\lambda=(a\, 1^{m-1})$, $a,m\ge1$ is ``hook'' (\cite{Elementary} for
$a=1$, \cite{Hook1} for $m$ odd, \cite{Hook2} for $m$
even). Furthermore, we calculate $\Rq_\lambda(z)$ where $\lambda=(a\,
b)$ is a two row diagram (\cite{Tworow}) or any $\Rq_\lambda$ if
$n=3$. For $n=3$ these are expressible in terms of Jacobi polynomials,
for $n=4$ we get one of Horn's hypergeometric function.

Even though most of the theory is parallel to that of (shifted) Jack
polynomials there are some differences. One of them is that the
$R_\lambda$ don't specialize for $r=0$ to anything easy. It seems that
matters become rather more involved. Also, neither $R_\lambda$ nor its
top homogeneous term $\Rq_\lambda$ seems to have any obvious
positivity properties (see \cite{KnSa} and \cite{KnSa2} for the
classical case). Another remarkable difference occurs when the number
of variables is even. Then the specialization of $\Rq_\lambda(z)$ at
the point $z=(1,\ldots,1)$ may be zero. This has the consequence that
the shifted polynomials $R_\lambda(z)$ cannot be defined via the
binomial formula \cite{E31} since not all of them occur in this
formula. A major open problem is orthogonality: Jack polynomials are
most commonly defined by an orthogonalization process with respect to
some explicit scalar product. Such a scalar product is still missing
for the $\Rq_\lambda$'s.

\medskip \noindent {\bf Acknowledgment:} I would like to thank Yasmine
Sanderson and the referee for valuable comments concerning the
exposition of this paper.

\beginsection Shifted. Shifted semisymmetric functions

Let $k$ be a field of characteristic zero. Consider the polynomial
ring $\cP:=k[z_1,\ldots,z_n]$. On it, the symmetric group $S_n$ acts
by permutation of the variables. The {\it semisymmetric group\/} is
the subgroup $W$ of $S_n$ which doesn't mix even and with odd entries:
$\pi\in W$ if $\pi(i)\equiv i \Mod 2$ for all $i$. Throughout this
paper, we are adopting the following notation: we put $\nuq:=\lfloor
n/2\rfloor$ and $\noq:=n-\nuq=\lceil n/2\rceil$. Then we have $W\cong
S_\noq\times S_\nuq$. For $z\in k^n$ we let
$z_\odd:=(z_1,z_3,\ldots,z_{2\noq-1})\in k^\noq$ and
$z_\even:=(z_2,z_4,\ldots,z_{2\nuq})\in k^\nuq$.

We are going to study the {\it ring of semisymmetric polynomials\/}
$\cP^W$. Clearly, as an algebra, $\cP^W$ is a polynomial ring
generated by $e_i(z_\odd)$, $i=1,\ldots,\noq$ and
$e_i(z_\even)$, $i=1,\ldots,\nuq$ where $e_i$ is
the $i$\_th elementary symmetric polynomial.

Let $\Lambda$ be the set of {\it partitions\/} of length $n$, i.e.,
$n$\_tuples of integers $\lambda=(\lambda_i)$ with
$\lambda_1\ge\ldots\ge\lambda_n\ge0$. We are going to consider
$\Lambda$ as a subset of $k^n$. The degree of $\lambda$ is
$|\lambda|=\sum_i\lambda_i$. We also define the odd degree
$|\lambda|_\odd:=|\lambda_\odd|$ and the even degree
$|\lambda|_\even:=|\lambda_\even|$. The odd degree will be for
semisymmetric polynomials what the degree is for symmetric
polynomials.

Finally, we choose once and for all a parameter $r\in k$
with\footnote*{In fact, only the slightly weaker condition $r\ne-{p\over
2q}$ where $p$ and $q$ are integers with $1\le p$ and $1\le q<
{n\over2}$ is needed.}  $r\not\in\QQ_{<0}$ and put
$$\eqno{}
\rho:=((n-1)r,(n-2)r,\ldots,2r,r,0).
$$

\Theorem Existence. For any $d\in\NN$ let $\Lambda(d)$ be the set of
$\lambda\in\Lambda$ with $|\lambda|_\odd\le d$. Let $\Lambda(d)\pfeil
k:\lambda\mapsto a_\lambda$ be any map. Then there is a unique
$f\in\cP^W$ with $\|deg|f\le d$ and $f(\rho+\lambda)=a_\lambda$ for
all $\lambda\in\Lambda(d)$.

\Proof:
We are using induction on $d+n$. To make the
dependence on the dimension $n$ explicit we write it as an index. We
have $\Lambda_{n-1}\into\Lambda_n$ by appending a zero. Moreover let
$\cP_{n-1}^{W_{n-1}}\pfeil\cP_n^{W_n}:g\mapsto g^+$ be the homomorphism
which maps $e_i(z_{\|odd|/\|even|})\in\cP_{n-1}^{W_{n-1}}$ to
$e_i(z_{\|odd|/\|even|})\in\cP_n^{W_n}$. Then $\|deg|g^+=\|deg|g$ and
$g^+(z_1,\ldots,z_{n-1},0)=g(z_1,\ldots,z_{n-1})$.

Let $e(z):=\prod_{i:\,n-i\ \even}z_i$. Observe $\|deg|e=\noq$. Since
$e(z)$ is the one generator of $\cP_n^{W_n}$ which is not in the image
of $\cP_{n-1}^{W_{n-1}}$, every $f\in\cP_n^{W_n}$ can be uniquely
expressed as $f(z)=g^+(z)+e(z)h(z)$ with $g\in\cP_{n-1}^{W_{n-1}}$,
$\|deg|g\le\|deg|f$, $h\in\cP_n^{W_n}$, and $\|deg|h\le\|deg|f-\noq$.

Now we split $\Lambda_n(d)$ into two parts $\Lambda_n(d)^0$ and
$\Lambda_n(d)^1$ according to whether the last component $\lambda_n$
is zero or not.

For any $g\in\cP_{n-1}$ let
$g_0(z):=g(z_1+r,\ldots,z_{n-1}+r)$. Clearly, we can identify
$\Lambda_n(d)^0$ with $\Lambda_{n-1}(d)$. Then for any
$\lambda\in\Lambda_n(d)^0$ we have
$$\eqno{}
g^+(\rho_n+\lambda)=g(\lambda_1+(n-1)r,\ldots,\lambda_{n-1}+r)=
g_0(\rho_{n-1}+\lambda).
$$
Since $e(\rho+\lambda)=0$ for every $\lambda$ with $\lambda_n=0$ the
system of linear equations $f(\rho_n+\lambda)=a_\lambda,
\lambda\in\Lambda_n(d)^0$ is equivalent to the system
$g_0(\rho_{n-1}+\lambda)=a_\lambda, \lambda\in\Lambda_{n-1}(d)$. By
induction on the number of variables we conclude that is has a unique
solution.

For any $\lambda\in\Lambda_n(d)^1$ holds $e(\rho+\lambda)\ne0$ since,
by assumption, $r\not\in\QQ_{<0}$. Thus, we can define
$a_\lambda':=(a_\lambda-g^+(\rho+\lambda))/e(\rho+\lambda)$. The map
$\lambda\mapsto\tilde\lambda:=(\lambda_1-1,\ldots,\lambda_n-1)$
identifies $\Lambda_n(d)^1$ with $\Lambda_n(d-\noq)$. Thus the system
of linear equations $f(\rho+\lambda)=a_\lambda,
\lambda\in\Lambda_n(d)^1$ is equivalent to the system $\tilde
h(\rho+\tilde\lambda)=a_\lambda', \tilde\lambda\in\Lambda_n (d-\noq)$
where $\tilde h(z)=h(z_1-1,\ldots,z_n-1)$. By induction on the degree
we conclude that is has a unique solution, as well.\qed

Now, we can define interpolation polynomials as follows:

\Definition: For every $\lambda\in\Lambda$ let
$r_\lambda(z;r)$ be the unique polynomial such that
\item{$\bullet$}it is $W$\_invariant,
\item{$\bullet$}its degree is $d:=|\lambda|_\odd$,
\item{$\bullet$}for all $\mu\in\Lambda$ with $|\mu|_\odd\le d$ holds
$r_\lambda(\rho+\mu;r)=\delta_{\lambda\mu}$ (Kronecker delta).
\medskip\noindent
The normalization $r_\lambda(\rho+\lambda;r)=1$ is very natural but
there is one which is often more convenient: the ``leading''
coefficient should equal to one. To define what that means, observe
that every $W$\_orbit of a monomial contains exactly one monomial, say
$z^\nu$, such that both $\nu_\odd$ and $\nu_\even$ are
partitions. These $\nu$ are in bijection with $\Lambda$. In fact, for
every partition $\lambda\in\Lambda$ we define the composition
$[\lambda]\in\NN^n$ by
$$\eqno{E19}
[\lambda]_m:=\lambda_m-\lambda_{m+1}+\ldots+(-1)^{n-m}\lambda_n.
$$
Since $[\lambda]_m=(\lambda_m-\lambda_{m+1})+[\lambda]_{m+2}$, both
$[\lambda]_\odd$ and $[\lambda]_\even$ are in fact
partitions. Conversely, let $\nu$ be a composition such that both
$\nu_\odd$ and $\nu_\even$ are partitions. Then
$$\eqno{}
\lambda=(\nu_1+\nu_2,\nu_2+\nu_3,\nu_3+\nu_4,\nu_4+\nu_5,\ldots)
$$
is in $\Lambda$. One easily checks that these two maps are
inverse to each other. Of special interest is the first component of
$[\lambda]$ since
$$\eqno{}
[\lambda]_1=|\lambda|_\odd-|\lambda|_\even=
\lambda_1-\lambda_2+\lambda_3-+\ldots
$$
In particular,
$$\eqno{E43}
[\lambda]_1=0\hbox{ if and only if
}\lambda_1=\lambda_2,\lambda_3=\lambda_4,\ldots,\hbox{ and
}\lambda_n=0\hbox{ in case $n$ is odd.}
$$
Moreover, we have $|[\lambda]|=|\lambda|_\odd$ and therefore
$\|deg|r_\lambda(z;r)=\|deg|z^{[\lambda]}$.

\Proposition. The coefficient $C_\lambda(r)$ of $z^{[\lambda]}$ in
$r_\lambda(z;r)$ is non-zero.

\noindent With this result we can define the renormalized polynomial
$$\eqno{}
R_\lambda(z;r):={1\over
C_\lambda(r)}r_\lambda(z;r)=z^{[\lambda]}+\ldots
$$
We are proving the proposition by computing $C_\lambda(r)$
explicitly. For this we need some more notation.  A partition
$\lambda$ can be represented by its diagram, i.e., the set of all
$(i,j)\in\NN^2$ (called boxes) such that $1\le i\le n$ and $1\le
j\le\lambda_i$. The dual partition $\lambda'$ is defined by the
transposed diagram $\{(j,i)\mid (i,j)\in\lambda\}$. For every box
$s=(i,j)\in\lambda$ we define the arm\_length
$a_\lambda(s):=\lambda_i-j$ and the leg\_length
$l_\lambda(s):=\lambda_j'-i$. Then we define
$$\eqno{}
[c'_\lambda(r)]_\even:=\prod\limits_{s\in\lambda\atop
l_\lambda(s)\ \even}\kern-10pt\relax (a_\lambda(s)+1+l_\lambda(s)r).
$$
For example, we have $[c'_{(a)}(r)]_\even=a!$, $[c'_{(a\,
b)}(r)]_\even=(a-b)!b!$, and $[c'_{(1^m)}(r)]_\even=
\prod_{1\le i<m\atop i\ \even}(1+ir)$.

\Lemma LeadingTerm. For every $\lambda\in\Lambda$ holds
$C_\lambda(r)=[c'_\lambda(r)]_\even^{-1}$. In particular, we have
$$\eqno{}
R_\lambda(\rho+\lambda;r)=[c'_\lambda(r)]_\even.
$$

\Proof: We retain the notation of the proof of \cite{Existence} and
prove the lemma by a similar induction. In particular, we have an
expression $r_\lambda(z)=g^+(z)+e(z)h(z)$.

If $\lambda_n=0$ then $e(\rho+\lambda)=0$ and therefore
$g(z)=r_{\lambda'}(z_1-r,\ldots,z_{n-1}-r)$ where the prime means ``drop
the last component''. Moreover, the coefficient of $z^{[\lambda]}$ in
$r_\lambda$ equals the one in $g$ (observe $[\lambda]_n=0$). Thus, we
get by induction
$C_\lambda(r)=C_{\lambda'}(r)=[c'_{\lambda'}(r)]_\even^{-1}$. But we
also have $[c'_\lambda(r)]_\even=[c'_{\lambda'}(r)]_\even$ which
finishes this case.

If $\lambda_n\ge1$ then $g(z)=0$ and
$h(z)=e(\rho+\lambda)^{-1}r_{\tilde\lambda}(z_1-1,\ldots,z_n-1)$. One checks
$z^{[\lambda]}=e(z)z^{[\tilde\lambda]}$. Thus, by induction, the
coefficient of $z^{[\lambda]}$ is
$e(\rho+\lambda)^{-1}[c'_{\tilde\lambda}(r)]_\even^{-1}$. But
$e(\rho+\lambda)$ is the contribution of the first column of $\lambda$
to $[c'_\lambda(r)]_\even$. Thus we get
$C_\lambda(r)=[c'_\lambda(r)]_\even^{-1}$, as claimed.\qed

\noindent
The second case of the preceding proof gives the following recursion
formula which allows to reduce the computation of $R_\lambda$ to the
case $\lambda_n=0$.

\Corollary Recursion. Let $\delta:=(1,\ldots,1)$. Then for every
$\lambda\in\Lambda$ with $\lambda_n\ge1$ holds
$$\eqno{E32}
R_\lambda(z;r)=
({\textstyle\prod\limits_{n-i\ \even}}\kern-7pt
z_i\ )\cdot R_{\lambda-\delta}(z-\delta;r).
$$

\noindent
We also have the following stability result:

\Proposition Stability. For $z=(z_1,\ldots,z_n)\in k^n$ let
$z':=(z_1,\ldots,z_{n-1})\in k^{n-1}$. Then we have for any
$\lambda\in\Lambda$:
$$\eqno{}
R_\lambda(z_1,\ldots,z_{n-1},0)=
\cases{R_{\lambda'}(z_1-r,\ldots,z_{n-1}-r)&if $\lambda_n=0$;\cr
0&otherwise.\cr}
$$

\Proof: If $\lambda_n\ge1$ then $R_\lambda$ is divisible by $z_n$
(\cite{Recursion}), hence $R_\lambda|_{z_n=0}=0$. Otherwise,
$R_\lambda|_{z_n=0}$ satisfies the definition of
$R_{\lambda'}(z_1-r,\ldots,z_{n-1}-r)$.\qed

\Remark: In many circumstances it is more convenient to consider the
polynomials $\tilde R_\lambda(u;r):=R_\lambda(\rho+u;r)$. Their main
advantage is that the stability result above can now be expressed as
$$\eqno{}
\tilde R_\lambda(u_1,\ldots,u_{n-1},0)=
\tilde R_{\lambda'}(u_1,\ldots,u_{n-1})
$$
whenever $\lambda_n=0$. This means that one can form a theory of
shifted semisymmetric polynomials which is independent of the
dimension $n$. For this one defines them in infinitely many variables
as follows. Let $\cP_\infty$ be the projective limit of the polynomial
rings $k[u_1,\ldots,u_n]$ in the category of filtered algebras. An
element of $\cP_\infty$ is a possibly infinite linear combination of
monomials in $u_1,u_2,\ldots$ whose degrees are uniformly bounded. Let
$\Lambda_\infty$ be the set of all descending sequences of integers
$\lambda_1\ge\lambda_2\ge\ldots$ with $\lambda_n=0$ for
$n>\!\!\!>0$. The stability result above says that for any $\lambda$
the sequence $(\tilde
R_{(\lambda_1,\ldots,\lambda_n)}(u_1,\ldots,u_n))_{n>\!\!\!>0}$ is an
element of $\cP_\infty$. It is denoted by $\tilde R_\lambda$.

The drawback of this method that the action of the semisymmetric group
gets distorted. More precisely, $W$ acts now by $\pi\bullet
u:=\pi(u+\rho)-\rho$. For example, the simple reflection
$s_{i\,i{+}2}$ acts as $u_i\mapsto u_{i+2}-2r$ and $u_{i{+}2}\mapsto
u_i+2r$. This action extends to an action of $W_\infty$ on
$\cP_\infty$ where $W_\infty$ is the group of parity preserving
permutations of $\NN$ with finite support. It is easy to show that the
$\tilde R_\lambda$, $\lambda\in\Lambda_\infty$ form a linear basis of
$\cP_\infty^{W_\infty}$.

\medskip

Next, we present some compatibility results with shifted Jack
polynomials. First, we recall their definition from \cite{KnSa}. More
or less, one has to replace the semisymmetric group by the full
symmetric group and the odd degree by the full degree. More precisely:
for each $\lambda\in\Lambda$ we define $P_\lambda(z_1,\ldots,z_n;r)$
as the unique polynomial having the following properties:
\item{$\bullet$}$P_\lambda$ is invariant under the full symmetric
group $S_n$; \item{$\bullet$}$\|deg|P_\lambda=|\lambda|$;
\item{$\bullet$}the coefficient of $z^\lambda$ is $1$;
\item{$\bullet$}$P_\lambda(\rho+\mu)=0$ for all $\mu\in\Lambda$ with
$|\mu|\le|\lambda|$ and $\mu\ne\lambda$.

\noindent Analogously, we define $\tilde
P_\lambda(u;r):=P_\lambda(\rho+u;r)$. 

Now we show that the symmetric polynomials $P_\lambda(z;r)$ are in two
ways special cases of the semisymmetric polynomials $R_\lambda(z;r)$.

\Theorem Jack1. Let $\lambda\in\Lambda$ with
$[\lambda]_1=0$ (see \cite{E43}). Then
$$\eqno{E5}
\tilde R_\lambda(u;r)=\tilde P_{\lambda_\even}(u_\even;2r).
$$

\Proof: We show that the polynomial $P$ on the right hand side matches
the definition of $\tilde R_\lambda$. First, observe that the shifted
action of $W$ induces on the even coordinates the shifted action of
$W_\even$ with parameter $2r$. Thus, $P$ is shifted
semisymmetric. Moreover, we have
$$\eqno{}
|\lambda|_\odd=|\lambda_\odd|=|\lambda_\even|
$$
which shows that the degree of $P$ is correct. Since
$[\lambda]=(0,\lambda_2,0,\lambda_4,\ldots)$ we have
$z^{[\lambda]}=z_\even^{\lambda_\even}$. This shows that the
normalization of $P$ is correct, as well.

It remains to check the vanishing conditions. For this let
$\mu\in\Lambda$ with
$|\mu|_\odd\le|\lambda|_\odd=|\lambda_\even|$. Then
$$\eqno{}
|\mu_\even|\le|\mu_\odd|\le|\lambda_\even|.
$$
This implies $P(\mu)=0$ unless $\mu_\even=\lambda_\even$. But then
$0\le[\mu]_1=|\mu|_\odd-|\mu_\even|\le|\lambda|_\odd-|\lambda_\even|=0$
which implies $\mu_1=\mu_2$, $\mu_3=\mu_4$, etc., i.e., $\mu=\lambda$.\qed

The other connection between $R_\lambda$ and
$P_\lambda$ is:

\Theorem Jack2. For every $\mu=(\mu_1,\ldots,\mu_\noq)\in\Lambda_\noq$ holds
$$\eqno{E2}
\sum_{\lambda:\,\lambda_\odd=\mu}{\tilde R_\lambda(u;r)\over
\tilde R_\lambda(\lambda;r)}=
{\tilde P_\mu
(u_\odd;2r)\over
\tilde P_\mu(\mu;2r)}.
$$

\Proof: Let $P$ be the polynomial on the right hand side of
\cite{E2}. Again, the shifted action of $W$ induces on the odd
coordinates the shifted action of $W_\odd$ with parameter $2r$. Thus,
$P$ is shifted semisymmetric. In particular, we have an expression
$$\eqno{E3}
P=\sum_\lambda
c_\lambda\tilde R_\lambda,\quad\hbox{with } |\lambda|_\odd\le
|\mu|.
$$
Suppose there is $\lambda$ with $c_\lambda\ne0$ and
$\lambda_\odd\ne\mu$. If we choose one of minimal degree, the left
hand side of \cite{E3} evaluates at $u=\lambda$ to $c_\lambda$ while
$P(\lambda_\odd)=0$. Contradiction. Thus $c_\lambda=0$ unless
$\lambda_\odd=\mu$. In that case, the value of $c_\lambda$ is
immediately obtained by evaluating both sides of \cite{E3} at
$u=\lambda$.\qed

As a corollary we get a formula for the elementary semisymmetric
polynomials:

\Corollary Elementary. 
$$
\eqalignno{%
&\tilde R_{(1^{2m-1})}(u;r)&=
\tilde P_{(1^m)}(u_\odd;2r)-\tilde P_{(1^m)}(u_\even;2r)\cr
E4&\tilde R_{(1^{2m})}(u;r)&=\tilde P_{(1^m)}(u_\even;2r)\cr}
$$

\Proof: Formula \cite{E4} is a special case of \cite{E5}. If we put
$\lambda=(1^{2m-1})$ in \cite{E2} and use \cite{E4} we get
$$\eqno{E6}
\tilde R_{(1^{2m-1})}(u;r)=
\alpha\tilde P_{(1^m)}(u_\odd;2r)-\beta\tilde P_{(1^m)}(u_\even;2r)
$$
with two constants $\alpha$ and $\beta$. Comparison of the coefficient
of $u^{[(1^{2m-1})]}=u_1u_3\ldots u_{2m{-}1}$ implies $\alpha=1$. Next
we evaluate \cite{E6} at $u=(1^{2m})$. The left-hand side is zero by
definition. Then $(1^{2m})_\odd=(1^m)=(1^{2m})_\even$ implies
$\beta=1$.\qed

Explicit formulas for $\tilde P_{(1^m)}(u;r)$ can be found in e.g.
\cite{KnSa}~3.1. One them is:
$$\eqno{}
\tilde P_{(1^m)}(u;r)=\sum_{n\ge i_1>i_2>\ldots>i_m\ge 1}
\prod_{j=1}^m(u_{i_j}+(j-1)r)
$$
Thus, the first few elementary semisymmetric polynomials are
$$
\eqalignno{
E24&\tilde R_{(1)}(u;r)&=e_1(u_\odd)-e_1(u_\even)=(u_1+u_3+\ldots)-(u_2+u_4+\ldots),\cr
&\tilde R_{(11)}(u;r)&=e_1(u_\even)=u_2+u_4+\ldots,\cr
&\tilde R_{(111)}(u;r)&=e_2(u_\odd)-e_2(u_\even)+r\sum_{i\ \odd}(i-1)u_i
-r\sum_{i\ \even}(i-2)u_i,\cr
&\tilde R_{(1111)}(u;r)&=e_2(u_\even)+r\sum_{i\ \even}(i-2)u_i.\cr}
$$
Let $\Rq_\lambda(z;r)$ be the top homogeneous component of
$R_\lambda(z;r)$. Since the highest degree component of $\tilde
P_{(1^m)}$ is the elementary symmetric function $e_m$ we obtain:

\Corollary Elementary2. Let $\Rq_\lambda(z;r)$ be the highest degree
component of $R_\lambda(z;r)$. Then
$$
\eqalignno{
&\Rq_{(1^{2m-1})}(z;r)&=e_m(z_\odd)-e_m(z_\even)\cr
&\kern9pt\Rq_{(1^{2m})}(z;r)&=e_m(z_\even)\cr}.
$$

We conclude this section with a list of all polynomials $R_\lambda$
which are non\_elementary of degree at most $3$, i.e., with
$|\lambda|_\odd\le3$ and $\lambda_1>1$. Each $R_\lambda$ is expressed
as a polynomial in the $R_{(1^i)}$. This means, that the formulas are
valid {\it for all $n$} with the convention that $R_\lambda=0$ if the
length of $\lambda$ is greater than $n$. 

$$\eqno{}
\vcenter{\halign{$#$\hfill$\,$&$#$\hfill\cr
R_{(2)} &= R_{(1)}^2- R_{(1)}\cr
R_{(21)} &= R_{(1)}R_{(11)}-{1\over 1+2r} R_{(111)}\cr
R_{(22)} &=  R_{(11)}^2 -{2\over 1+2r} R_{(1111)}- R_{(11)}\cr
R_{(211)} &= R_{(1)}R_{(111)}- R_{(111)}\cr
R_{(2111)} &= R_{(1)}R_{(1111)}-{1\over 1+4r}R_{(11111)}\cr
R_{(221)} &= R_{(11)} R_{(111)}-{1\over 1+2r}R_{(1)}R_{(1111)}-{1\over 1+2r}R_{(11111)}- R_{(111)}\cr
R_{(2211)} &= R_{(11)}R_{(1111)}-{3\over 4r+1}R_{(111111)}-2 R_{(1111)}\cr
R_{(3)} &= R_{(1)}^3-3R_{(1)}^2+2 R_{(1)}\cr
R_{(31)} &= R_{(1)}^2R_{(11)}-{1\over 1+r}R_{(1)}R_{(111)}-R_{(1)}R_{(11)}+{1\over 1+r} R_{(111)}\cr
R_{(32)} &= R_{(1)}R_{(11)}^2-{1\over 1+r}R_{(11)}R_{(111)}- {1\over 1+r} R_{(1)}R_{(1111)}+{1\over (1+r)(1+2r)}R_{(11111)}\cr
&\phantom{=}-R_{(1)} R_{(11)}+{1\over 1+r} R_{(111)}\cr
R_{(33)} &= R_{(11)}^3-{3\over 1+r}R_{(11)}R_{(1111)}+{3\over (1+r)(1+2r)}R_{(111111)}-3R_{(11)}^2+ {6\over 1+r} R_{(1111)}\cr&\phantom{=}+2 R_{(11)}\cr}}$$

\noindent These formulas were obtained with the help of a computer. 


\beginsection MultiplicityFree. Representation theoretic interpretation

Before we study the polynomials $R_\lambda(z;r)$ further, we describe
the representation theoretic interpretation of the three special cases
which are mentioned in the introduction. For this, we recall some
basic facts about multiplicity free representations. Details appeared
for example in \cite{Mon}.

Let $G$ be a connected complex reductive group. A finite dimensional
representation $V$ of $G$ is called {\it multiplicity free\/} if every
simple $G$\_module appears in $\PP:=\CC[V]$ at most once. Equivalent
to this condition is that a Borel subgroup of $G$ has a dense orbit in
$V$. Thus, as a $G$\_module, we have a decomposition
$\PP=\oplus_{\lambda\in\Lambda_V}\PP_\lambda$ where $\PP_\lambda$ is
the simple module with {\it lowest\/} weight $-\lambda$. Then
$\Lambda_V$ is a set of dominant weights which can be shown to be a
free abelian monoid (i.e., isomorphic to $\NN^r$). Clearly, all
non\_zero polynomials in $\PP_\lambda$ have the same degree, denoted
by $|\lambda|$.

The symmetric algebra $\DD:=S^*(V)$ then decomposes accordingly as
$\DD=\oplus_{\lambda\in\Lambda_V}\DD_\lambda$ where $\DD_\lambda$ is
isomorphic to $\PP_\lambda^*$. In particular, $\lambda$ is its highest
weight. The space $\DD$ can be interpreted either as polynomial
functions on $V^*$ or as constant coefficient differential operators
on $V$. Accordingly, $\PP\otimes\DD$ can be identified with either the
algebra of polynomial functions on $V\oplus V^*$ or the algebra
$\PP\DD$ of linear differential operators on $V$ with polynomial
coefficients.

The point is now that the space of $G$\_invariants $(\PP\otimes\DD)^G$
comes with a (up to scalars) distinguished basis: we have
$$\eqno{}
(\PP\otimes\DD)^G=
\mathop{\oplus}
\limits_{\lambda,\mu\in\Lambda_V}(\PP_\lambda\otimes\DD_\mu)^G.
$$
Each summand is zero unless $\lambda=\mu$ in which case it is
one\_dimensional (Schur's Lemma). We denote a generator as $E_\lambda$
if regarded as a function on $V\oplus V^*$ (called a {\it zonal
spherical function}) and $D_\lambda$ if regarded as a differential
operator (called a {\it Capelli operator}).

The Capelli operators are easier accessible, whence we start with
them. Each differential operator $D\in(\PP\DD)^G$ acts on
$\PP_\lambda$ by a scalar denoted $c_D(\lambda)$. Recall that
$\Lambda_V$ is a set of weights and therefore sits in $\ft^*$, the
dual of the Cartan subalgebra. Let $\fa^*$ be its $\CC$\_span. Let
$W\subseteq GL(\ft^*)$ be the Weyl group and let
$\overline\rho\in\ft^*$ be the half\_sum of the positive roots. Then
the shifted action of $W$ on $\ft^*$ is defined by
$w\bullet\chi=w(\chi+\overline\rho)-\overline\rho$.

\Theorem DiffOp. {\rm(\cite{Mon} 4.4, 4.8, 4.9, 4.7)}. Let $V$ be a
multiplicity free representation.
\item{a)} Each $c_D$ is the restriction of a unique
polynomial (also denoted $c_D$) on $\fa^*$. 
\item{b)} There is a subgroup $W_V\subseteq W$ such that
$\fa^*\subseteq\ft^*$ is $W_V$\_stable with respect to the shifted
action and such that $D\mapsto c_D$ is an algebra isomorphism of
$(\PP\DD)^G$ with $\CC[\fa^*]^{W_V^\bullet}$, the space of shifted
$W_V$\_invariant polynomials on $\fa^*$.
\item{c)} The ``little Weyl group'' $W_V$ acts as a reflection group
on $\fa^*$. In particular, $(\PP\DD)^G$ and $\CC[\fa^*]^{W_V^\bullet}$
are polynomial rings.\Par

\noindent Since $(\PP\DD)^G$ has a distinguished basis we obtain a
basis $c_\lambda=c_{D_\lambda}$ of $\CC[\fa^*]^{W_V^\bullet}$. There is
a purely combinatorial characterization of the $c_\lambda$:

\Theorem Characterization. {\rm(\cite{Mon} 4.10)} The polynomial
$c_\lambda\in\CC[\fa^*]^{W_V^\bullet}$ is, up to a scalar factor,
characterized by the vanishing condition $c_\lambda(\mu)=0\hbox{ for
all }\mu\in\Lambda_V\hbox{ with } |\mu|\le|\lambda|\hbox{ and
}\mu\ne\lambda$.

One can eliminate the shifted action of $W_V$ as follows: choose a
$W_V$\_stable complement $\fa_0$ of $\fa^*$ in $\ft^*$ and let
$\overline\rho=\rho+\rho_0$ with $\rho\in\fa^*$ and
$\rho_0\in\fa_0$. The condition that $\fa^*$ is shifted $W_V$\_stable
means $w\overline\rho-\overline\rho\in\fa^*$ for all $w\in W_V$. Thus,
$\rho_0$ is $W_V$\_fixed. Therefore, we can define the shifted
$W_V$\_action as well with $\overline\rho$ replaced by
$\rho$. Actually, one can add to $\rho$ any fixed vector in $\fa^*$
without changing the shifted action. The point is now that
$p_\lambda(\chi):=c_\lambda(\chi-\rho)$ is a truly $W_V$\_invariant
polynomial on $\fa^*$.

\Corollary Characterization2. The polynomial
$p_\lambda\in\CC[\fa^*]^{W_V}$ is, up to a scalar factor,
characterized by the vanishing condition $p_\lambda(\rho+\mu)=0\hbox{ for
all }\mu\in\Lambda_V\hbox{ with } |\mu|\le|\lambda|\hbox{ and
}\mu\ne\lambda$.

Now we say a few words about the zonal spherical functions
$E_\lambda$. One of their main features is that they have many
different interpretations. First, we can consider $V\oplus V^*$ as the
cotangent bundle of $V$. Then {\it the symbol of $D_\lambda$ is
$E_\lambda$.}  This is our principal method for their study.

It is possible to define $E_\lambda$ without reference to Capelli
operators. Every differential operator $D\in\PP\DD(V)$ is also a
differential operator on $V\oplus V^*$ by acting on the first
argument. As such it is denoted $D^{(1)}$. Observe, that the
eigenspaces of $\PP\DD^G$ are then just the spaces
$\PP_\lambda\otimes\DD$. Therefore one can characterize
$E_\lambda$ as  {\it the (up to scalar) unique $G$\_invariant function $f$ on
$V\oplus V^*$ with $D^{(1)}(f)=c_D(\lambda)f$ for all\/
$D\in\PP\DD^G$.} Clearly, it suffices to let $D$ run through a set of
generators of $\PP\DD^G$.

There is also a ``Chevalley isomorphism'' for $G$\_invariant functions
on $V\oplus V^*$:

\Theorem SpherFunc. {\rm(\cite{Mon} 4.2, 4.8, 4.5)} There is $v^*\in V^*$
and a linear embedding
$\fa^*\into V$ such that the restriction map $f\mapsto f|_{\fa^*\times v^*}$
induces an isomorphism $(\PP\otimes\DD)^G\pf\sim\CC[\fa^*]^{W_V}$.
Moreover, the image of the symbol of $D\in(\PP\DD)^G$
is the highest degree component $\cq_D$ of $c_D$. In particular,
$E_\lambda$ is mapped to $\cq_\lambda$.

\noindent The subspace $\fa^*$ is constructed as follows: choose
$v^*\in V^*$ in the open $G$\_orbit. Then choose a Borel subalgebra
$\fb=\ft\oplus\fu\subseteq\|Lie|G$ such that $\fb v^*=V^*$. This is
possible, since also $V^*$ has a dense orbit for any Borel
subgroup. The surjective map $\fb\auf V^*:\xi\mapsto\xi v^*$ induces
the dual injective map $\iota:V\into\fb^*$. Via the projection
$\fb\auf\fb/\fu=\ft$ we have $\ft^*\subseteq\fb^*$. Now, one can show
that $\iota(V)\cap\ft^*=\fa^*$ which furnishes us with the desired
embedding $\fa^*\into V$.

\cite{SpherFunc} indicates another way to interpret $E_\lambda$. Let
$H^*\subseteq G$ be the isotropy subgroup of $v^*$. Then the orbit
$Gv^*$ is open in $V^*$ and isomorphic to $G/H^*$. Therefore, a
function $f$ on $V\oplus V^*$ is $G$\_invariant if and only if its
restriction to $V\times v^*$ is $H^*$\_invariant. Thus the restriction
$E_{v^*,\lambda}(v):=E_\lambda(v,v^*)$ is {\it the (up to a scalar)
unique $H^*$\_invariant function $f$ on $V$ with $D(f)=c_D(\lambda)f$
for all $D\in\PP\DD^G$.} The restriction map from $V$ to $\fa^*$
defines now an isomorphism of the algebra of $H^*$\_invariants with
$\CC[\fa^*]^{W_V}$. Thereby, the function $E_{v^*,\lambda}$ is mapped
to the highest degree component $\cq_\lambda$ of $c_\lambda$.

Another interpretation is as follows: let $K\subseteq G$ be a maximal
compact subgroup. Let $V_\RR$ be $V$ regarded as a real vector
space. It is equipped with a complex conjugation $v\mapsto\vq$. Then
we can regard $\DD$ as the algebra of polynomials in the
antiholomorphic variables $\zq_i$ and $\PP\otimes\DD$ is the algebra
of all $\CC$\_valued polynomials on $V_\RR$. Thus, the polynomial
$v\mapsto E_\lambda(v,\vq)$ is {\it the (up to a scalar)
unique $K$\_invariant function $f$ on $V_\RR$ with
$D(f)=c_D(\lambda)f$ for all $D\in\PP\DD^G$.}

Observe that also $V$ has a dense $G$\_orbit $Gv$ with isotropy group
denoted by $H$. By restriction we can interpret $E_\lambda$ also as
$G$\_invariant function on $G/H\times G/H^*$, as an $H^*$\_invariant
function on $G/H$, or as a function on $G$ which is constant on double
cosets for $H$ and $H^*$. In this last form, $E_\lambda$ can be
interpreted purely representation theoretically: Let $M$ a simple
$G$\_module which is isomorphic to $\PP_\lambda$ for some
$\lambda$. Then $M^H$ and $(M^*)^{H^*}$ are both one dimensional,
generated by vectors $m_\lambda$ and $\alpha_\lambda$,
respectively. Then $g\mapsto E_\lambda(gv,v^*)$ equals (up to a
scalar) {\it the matrix coefficient
$g\mapsto\alpha_\lambda(gm_\lambda)$.} In fact, if we identify $M$
with $\cP_\lambda$ then $m$ is just the evaluation $f\mapsto
f(v)$. Similarly, $M^*\cong\DD_\lambda$ and $\alpha$ is evaluation in
$v^*$.  Finally, $E_\lambda$ corresponds to the canonical pairing
$M\times M^*\pfeil\CC$.

\medskip

Now we are in the position to explain the representation theoretic
relevance of the polynomials $R_\lambda(z;r)$.  \medskip \noindent
{\it The case of\/ $G=GL_p(\CC)\times GL_q(\CC)$ acting on
$V:=(\CC^p\otimes\CC^q)\oplus\CC^q$.}

\noindent The following data are taken from \cite{Mon}~p.~315. Put
$n:=\|min|(2p+1,2q)$. Then $\noq=\|min|(p+1,q)$ and
$\nuq=\|min|(p,q)$. Let $\epsilon_i$ and $\epsilon_i'$ be the weights
of the defining representation of $GL_p(\CC)$ and $GL_q(\CC)$,
respectively. Moreover, let $\omega_i:=\sum_{j=1}^i\epsilon_i$,
$\omega_i':=\sum_{j=1}^i\epsilon_i'$. Then $\Lambda_V$ is the free
abelian monoid generated by $\omega_{i-1}+\omega_i'$ for
$i=1,\ldots,\noq$ and $\omega_i+\omega_i'$ for
$i=1,\ldots,\nuq$. Thus, if we put $e_{2i-1}:=\epsilon_i'$ for
$i=1,\ldots,\noq$ and $e_{2i}:=\epsilon_i$ for $i=1,\ldots,\nuq$ then
$\Lambda_V$ consists of all $\chi=\sum_{i=1}^n\lambda_ie_i$ where
$(\lambda_1,\ldots,\lambda_n)$ is a partition. The degree function is
such that $|\omega_i|=0$ and $|\omega_i'|=i$ which translates into
$|\chi|=|\lambda|_\odd$. The little Weyl group consists of all
permutations of the $\epsilon_i$ and $\epsilon_i'$ separately, i.e.,
$W_V$ is the semisymmetric group. Finally, we have
$\overline\rho=({p-1\over2},
{p-3\over2},\ldots;{q-1\over2},{q-3\over2},\ldots)$. Thus, if we
project it to the first $\noq+\nuq$ coordinates and shift it by a
suitable $W_V$\_fixed vector we arrive at $\rho=({n-1\over2},
{n-3\over2},\ldots;{n-2\over2},
{n-4\over2},\ldots)={1\over2}\sum_i(n-i)e_i$. This shows
$r={1\over2}$. In particular, $c_\chi(x)$ is a multiple of $\tilde
R_\lambda(x;{1\over2})$.

Now we describe the combinatorics in more classical term. For this, it
is convenient to write
$$\eqno{}
V=(\CC^p\oplus\CC)\otimes\CC^q=\CC^{p+1}\otimes\CC^q,
$$
i.e., $V$ is the space of $(p+1)\times q$\_matrices $X$ acted upon by
$$\eqno{}
G=GL_p(\CC)\times GL_q(\CC)\subseteq \Gq:=GL_{p+1}(\CC)\times GL_q(\CC)
$$
by $X\mapsto AXB^t$ with $A\in GL_p(\CC)\subseteq GL_{p+1}(\CC)$ and
$B\in GL_q(\CC)$.

Let $\Lambda_\infty$ be the set of infinite partitions, i.e.,
descending sequences of integers $\tau_1\ge\tau_2\ge\ldots$ with
$\tau_i=0$ for $i>\!\!\!>0$. The length $\ell(\tau)$ is the maximal
$i$ with $\tau_i\ne0$. Every $\tau\in\Lambda_\infty$ with
$\ell(\tau)\le p$ parametrizes an irreducible (polynomial)
representation $M_\tau^{(p)}$ of $GL_p(\CC)$.

Let $\noq:=\|min|(p+1,q)$. Then it is well known (see e.g. \cite{GW}
Thm.~5.2.7) that there is a decomposition of $\Gq$\_modules:
$$\eqno{}
\PP=\sum_{\tau\in\Lambda_\infty\atop\ell(\tau)\le\noq}M_\tau^{(p+1)}\otimes
M_\tau^{(q)}.
$$
Recall also the branching law of $GL_{p+1}(\CC)$ to $GL_p(\CC)$ (see
e.g. \cite{GW} Thm.~8.1.1): as a
$GL_p(\CC)$\_module we have
$$\eqno{}
M_\tau^{(p+1)}=\sum_\sigma M_\sigma^{(p)}
$$
where $\sigma$ runs through all partitions with $\ell(\sigma)\le p$
and which are ``interlaced'' with $\tau$, i.e., with
$\tau_1\ge\sigma_1\ge\tau_2\ge\sigma_2\ge\ldots$. Thus we have
$\ell(\sigma)\le\nuq:=\|min|(p,q)$. Now we can make the decomposition
of $\PP$ into simple $G$\_modules more explicit. Combine $\tau$ and
$\sigma$ to a single partition $\lambda$ by putting
$\lambda_{2i-1}:=\tau_i$ and $\lambda_{2i}:=\sigma_i$. Then, as a
$G$\_module, we have
$$\eqno{}
\PP=\sum_{\lambda\in\Lambda}
\PP_\lambda\qquad\hbox{with}\qquad \PP_\lambda=
M_{\lambda_\even}^{(p)}\otimes M_{\lambda_\odd}^{(q)}.
$$
Here, we use the fact that $\ell(\lambda)\le
n:=\noq+\nuq=\|min|(2p+1,2q)$. Therefore, one can regard $\lambda$ as
an element of $\Lambda=\Lambda_n$.

This gives also a nice interpretations of the comparison
theorems~\ncite{Jack1} and~\ncite{Jack2}. Let $V'=\CC^p\otimes\CC^q$
be the space of $p\times q$ matrices. Since $V$ projects onto $V'$ we
have $\CC[V']\subseteq\PP$. More precisely,
$$\eqno{}
\CC[V']=\sum_{\lambda\in\Lambda\atop\lambda_\even=\lambda_\odd}\PP_\lambda.
$$
Thus, every Capelli operator on $V'$ can be regarded as
a Capelli operator on $V$. This is reflected in formula~\cite{E5}.

On the other hand each $\Gq$\_invariant Capelli operator on $V$
decomposes as a sum of $G$\_invariant Capelli operators on $V$. This
is the origin of formula \cite{E2}.

We can make this fully explicit for the generators $D_{(1^a)}$. Let
$A\in V$ be a $(p+1)\times q$\_matrix. For subsets
$I\subseteq[p+1]:=\{1,\ldots,p+1\}$ and
$J\subseteq[q]:=\{1,\ldots,q\}$ of the same size $i$ let
$$\eqno{}
\|det|_I^J(A)=\|det|\big(a_{ij}\big)_{i\in I\atop j\in J}
$$
be the corresponding minor. These form a basis of
$M_{(1^i)}^{(p+1)}\otimes
M_{(1^i)}^{(q)}=\wedge^i(\CC^{p+1})^*\otimes\wedge^i(\CC^q)^*$. If $V$
is parametrized by the coordinate functions $a_{ij}$ let $\partial_A$
be the matrix with entries ${\partial\over\partial a_{ij}}$. Then the
classical $\Gq$\_invariant Capelli operators on $V$ are
$$\eqno{}
C_i:=\sum_{I\subseteq[p+1],J\subseteq[q]\atop|I|=|J|=i}
\|det|_I^J(A)\|det|_I^J(\partial_A),\qquad i=1,\ldots,\noq
$$
Now each $\wedge^i(\CC^{p+1})^*$ decomposes as a $G$\_module into two pieces:
$$\eqno{}
\wedge^i(\CC^{p+1})^*=\wedge^i(\CC^p\oplus\CC)^*
=\wedge^i(\CC^p)^*\oplus\wedge^{i-1}(\CC^p)^*
$$
Thus also $C_i$ decomposes as $C_i=D_{(1^{2i})}+D_{(1^{2i-1})}$ with
$$\eqno{E55}
D_{(1^{2i})}=\sum_{I\subseteq[p],J\subseteq[q]\atop|I|=|J|=i}
\|det|_I^J(A)\|det|_I^J(\partial_A),\quad
D_{(1^{2i-1})}=\sum_{I\subseteq[p+1],J\subseteq[q]\atop p+1\in I,|I|=|J|=i}
\|det|_I^J(A)\|det|_I^J(\partial_A)
$$
For example, for $n=3$, i.e. $p+1=q=2$, we have
$$\eqno{E63}
D_{(1)}=a_{21}{\partial\over\partial a_{21}}+a_{22}{\partial\over\partial a_{22}}
$$
$$\eqno{E61}
D_{(11)}=a_{11}{\partial\over\partial a_{11}}+a_{12}{\partial\over\partial a_{12}}
$$
$$\eqno{E62}
D_{(111)}=(a_{11}a_{22}-a_{12}a_{21})({\partial\over\partial a_{11}}
{\partial\over\partial a_{22}}-
{\partial\over\partial a_{12}}{\partial\over\partial a_{21}})
$$

Now, we explain the zonal spherical functions. We identify $V^*$ with
the space of $q\times(p+1)$\_matrices and the pairing $V\times
V^*\pfeil\CC$ is given by $(A,A^*)\mapsto\|tr|(AA^*)$. By definition,
$E_\lambda$ is a $G$\_invariant functions on $V\oplus V^*$ which is a
joint eigenvector of the differential operators $D_{(1^i)}$ (see
\cite{E55}) acting on the first factor. To make the Chevalley
isomorphism from \cite{SpherFunc} explicit, we define for $p,q\ge 1$
the following two matrices $\Xi_{p,q}=\Xi_{p,q}(z_1,\ldots,z_n)\in V$
and $\Xi_{p,q}^*\in V^*$:
$$\eqno{E54}
(\Xi_{p,q})_{ij}=\cases{
z_{2i}&if $i=j$, $i\le p$\cr
z_{2i}-z_{2i+1}&if $i<j$, $i\le p$\cr
u_{2j}&if $i=p+1$\cr
0&otherwise\cr}\quad\hbox{and}\quad
(\Xi_{p,q}^*)_{ij}=\cases{1&if $i=j$\cr
1& if $i=q<j=p+1$\cr
0&otherwise\cr}
$$
Here, we put $u_i:=z_1-z_2+z_3-+\ldots\pm z_i$ and $z_i=0$ for $i>2p+1$.

For example, for $p=3<q$, $n=7$ we have
$$\eqno{}
\Xi_{3,q}(z)=\pmatrix{z_2&z_2-z_3&z_2-z_3&z_2-z_3&\cdots\cr
         0&z_4&z_4-z_5&z_4-z_5&\cdots\cr
         0&0&z_6&z_6-z_7&\cdots\cr
         u_2&u_4&u_6&u_7&\cdots\cr},\qquad
\Xi_{3,q}^*=\pmatrix{
1&0&0&0\cr
0&1&0&0\cr
0&0&1&0\cr
0&0&0&1\cr
0&0&0&0\cr
\vdots&\vdots&\vdots&\vdots\cr}
$$
while for $p\ge q=3$, $n=6$ we have
$$\eqno{}
\Xi_{p,3}(z)=\pmatrix{
z_2&z_2-z_3&z_2-z_3\cr
0&z_4&z_4-z_5\cr
0&0&z_6&\cr
0&0&0\cr
\vdots&\vdots&\vdots\cr
u_2&u_4&u_6\cr},\qquad
\Xi_{p,3}^*=\pmatrix{
1&0&0&0&\cdots&0\cr
0&1&0&0&\cdots&0\cr
0&0&1&0&\cdots&1\cr}
$$

Let $B_p\subseteq GL_p(\CC)$ be the subgroup of upper triangular
matrices. Let $B_q\subseteq GL_q(\CC)$ be the stabilizer of the flag
$\<v_1\>,\<v_1,v_2\>,\ldots$ where $v_i:=e_i+\ldots+e_q$ and where
$e_i$ is the $i$-th canonical basis vector of $\CC^q$. Then
$B:=B_p\times B_q$ is a Borel subgroup of $G$. One can verify by a
straightforward but tedious calculation that $z\mapsto\Xi_{p,q}(z)$ is
the embedding $\fa^*\into V$ when one follows the recipe described
after \cite{SpherFunc} using $v^*=\Xi_{p,q}^*$ and the Borel subgroup
$B$.

It follows that for every zonal spherical function $E_\lambda(A,A^*)$
the restriction $E_\lambda(\Xi(z),\Xi^*)$ is proportional to
$\Rq_\lambda(z;{1\over2})$. Since, $E_\lambda$ is the symbol of
$D_\lambda$ we obtain $\Rq_\lambda(z;{1\over2})$ also from $D_\lambda$
by replacing all coordinate functions $a_{ij}$ by $\Xi_{ij}(z)$ and
all derivations ${\partial\over\partial a_{ij}}$ by $1$ if $i=j$ or
$i=p+1>j=q$ and $0$ otherwise. For example, in the case $p+1=q=2$ we
get according to \cite{E63}--\cite{E62}:
$$\eqno{}
D_{(1)}\mapsto u_3=z_1-z_2+z_3=\Rq_{(1)}(z),\qquad D_{(11)}\mapsto z_2=\Rq_{(11)}(z)
$$
$$\eqno{}
D_{(111)}\mapsto z_2u_3-(z_2-z_3)u_2=z_1z_3=\Rq_{(111)}(z).
$$

For the other interpretations of zonal spherical functions we just
mention the case when $p+1=q$ since that is the only case when the
isotropy groups $H$ and $H^*$ are reductive. In fact, in that case we
have $H=H^*=GL_{q-1}(\CC)$ embedded diagonally into $G$. Thus, the action
of $H^*$ on $V$ is just the action of $GL_{q-1}(\CC)\subseteq
GL_q(\CC)$ by conjugation on $q\times q$\_matrices. The matrix
$\Xi_{p,q}^*$ is the identity matrix $I_q$. Thus, the function
$E_\lambda(A,I_q)$ is a joint eigenfunction of the differential
operators $D_{(1^i)}$ which is invariant under conjugation by
$GL_{q-1}(\CC)$. Any conjugation invariant function is uniquely
determined by its value in $\Xi_{p,q}(z)$ and we have
$E_\lambda(\Xi_{p,q}(z),I_q)=\Rq_\lambda(z;{1\over2})$.

\medskip\noindent
{\it The case of\/ $G=GL_n(\CC)$ acting on $V=\wedge^2\CC^n\oplus\CC^n$:}

\noindent We keep the notation of the previous example. According to
the data in \cite{Mon}~p.~314, the weight monoid $\Lambda_V$ is freely
generated by $\omega_i$ for $i=1,\ldots,n$. Thus, if we set
$e_i:=\epsilon_i$ for all $i$ then $\Lambda_V$ consists of all
$\chi=\sum_{i=1}^n\lambda_ie_i$ where $(\lambda_1,\ldots,\lambda_n)$
is a partition. The degree function is given by $|\omega_i|=\lceil
{i\over2}\rceil$. Thus $|\chi|=|\lambda|_\odd$. The little Weyl group
permutes the $\epsilon_i$ with even and odd indices separately and
therefore equals the semisymmetric group. Finally,
$\overline\rho=({n-1\over2},{n-3\over2},\ldots)$. Thus we can choose
$\rho=\sum_i(n-i)e_i$ which shows $r=1$. In particular, $c_\chi(x)$ is
a multiple of $\tilde R_\lambda(x;1)$.

Again, this can be made more explicit. Observe, that $G$ is a subgroup
of $\Gq:=GL_{n+1}(\CC)$ and $V$ is the restriction of the natural
$\Gq$\_action on $\wedge^2\CC^{n+1}$ to $G$. It is known (see
e.g. \cite{GW} Thm.~5.2.11) that as a $\Gq$\_module:
$$\eqno{}
\PP=\sum_\tau M_\tau^{(n+1)}
$$
where $\tau$ runs through all partitions with $\tau_1=\tau_2$,
$\tau_3=\tau_4,\ldots$ and $\ell(\tau)\le n+1$. Now, we use again the
$GL_{n+1}-GL_n$ branching rule. For $\lambda$ to be interlaced with
$\tau$ means $\lambda_1=\tau_1$, $\lambda_3=\tau_3,\ldots$. Thus, as a
$G$\_module, we obtain
$$\eqno{}
\PP=\sum_\lambda M_\lambda^{(n)}
$$
where $\lambda$ runs through all partitions with $\ell(\lambda)\le n$.
Here $M_\lambda^{(n)}$ is sitting in $M_{\lambda^*}^{(n+1)}$ where
$\lambda^*=(\lambda_1,\lambda_1,\lambda_3,\lambda_3,\ldots)$.

An element of $V$ is represented by a skewsymmetric matrix
$A=\big(a_{ij}\big)$ of size $n+1$. For $I\subseteq[n+1]$ of {\it
even\/} size $2m$ let
$$\eqno{}
\|Pf|_I(A):=\|Pfaffian|\big(a_{ij}\big)_{i\in I\atop j\in I}.
$$
Then, the Capelli operators for $\Gq$ corresponding to simple weights
are
$$\eqno{}
C_m:=\sum_{I\subseteq[n+1]\atop|I|=2m}\|Pf|_I(A)\|Pf|_I(\partial_A).
$$
Each $\Gq$\_module $M_{(1^{2m})}^{(n+1)}$ decomposes into (at most) two
components, namely $M_{(1^{2m-1})}^{(n)}$ and
$M_{(1^{2m})}^{(n)}$. Therefore, also $C_m$ decomposes as
$C_m=D_{(1^{2m-1})}+D_{(1^{2m})}$ where
$$\eqno{E56}
D_{(1^{2m})}=\sum_{I\subseteq[n]\atop|I|=2m}\|Pf|_I(A)\|Pf|_I(\partial_A),
\quad
D_{(1^{2m-1})}=\sum_{I\subseteq[n+1]\atop n+1\in I,|I|=2m}
\|Pf|_I(A)\|Pf|_I(\partial_A).
$$
For example, for $n=3$ we get
$$\eqno{}
D_{(1)}=a_{14}{\partial\over\partial a_{14}}+
a_{24}{\partial\over\partial a_{24}}+a_{34}{\partial\over\partial a_{34}}
$$
$$\eqno{}
D_{(11)}=a_{12}{\partial\over\partial a_{12}}+
a_{31}{\partial\over\partial a_{31}}+a_{23}{\partial\over\partial a_{23}}
$$
$$\eqno{}
D_{(111)}=(a_{12}a_{34}-a_{13}a_{24}+a_{23}a_{14})\left(
{\partial\over\partial a_{12}}{\partial\over\partial a_{34}}-
{\partial\over\partial a_{13}}{\partial\over\partial a_{24}}+
{\partial\over\partial a_{23}}{\partial\over\partial a_{14}}\right)
$$

To describe the zonal spherical functions we identify $V^*$ also with
skewsymmetric matrices of size $n+1$ and pairing $V\times
V^*\pfeil\CC:(A,A^*)\mapsto{1\over2}\|tr|(AA^*)$. The function
$E_\lambda$ is a $G$\_invariant functions on $V\oplus V^*$ which is a
joint eigenvector of the differential operators $D_{(1^i)}$ defined in
\cite{E56} acting on the first argument.

To make the Chevalley isomorphism explicit we define skewsymmetric
matrices $T(z)$ and $T^*$:
$$\eqno{}
T(z):=\pmatrix{0&-\Xi_{\nuq,\noq}^t\cr\Xi_{\nuq,\noq}&0\cr},
\quad T^*:=\pmatrix{0&-\Xi_{\nuq,\noq}\cr\Xi_{\nuq,\noq}^t&0\cr}
$$
where $\Xi$ and $\Xi^*$ are defined in \cite{E54}. For example, for
$n=4$ we get (again setting $u_i:=z_1-z_2+z_3-+\ldots\pm z_i$):
{\def\!{\mskip-3mu}\nospace
$$\eqno{}
T(z)=\pmatrix{
0&0&-z_2&0&-u_2\cr
0&0&\!\!\!-z_2+z_3\!\!\!&-z_4&-u_4\cr
z_2&\!\!\!z_2-z_3\!\!\!&0&0&0\cr
0&z_4&0&0&0\cr
u_2&u_4&0&0&0\cr},\quad
T^*=\pmatrix{
0&0&\hidew{-1}&0&0\cr
0&0&0&\hidew{-1}&\hidew{-1}\cr
1&0&0&0&0\cr
0&1&0&0&0\cr
0&1&0&0&0\cr}
$$}
while for $n=5$ one has
{\def\!{\mskip-3mu}\thickmuskip0mu\thinmuskip0mu\medmuskip0mu
$$\eqno{}
\mskip0mu T(z)=\pmatrix{
0&0&0&-z_2&0&-u_2\cr
0&0&0&\!\!\!-z_2+z_3\!\!\!&-z_4&-u_4\cr
0&0&0&\!\!\!-z_2+z_3\!\!&\!\!-z_4+z_5\!\!\!&-u_5\cr
z_2&\!\!\!z_2-z_3\!\!&\!\!z_2-z_3\!\!\!&0&0&0\cr
0&z_4&\!\!\!z_4-z_5\!\!\!&0&0&0\cr
u_2&u_4&u_5&0&0&0\cr},\quad
T^*=\pmatrix{
0&0&0&\hidew{-1}&0&0\cr
0&0&0&0&\hidew{-1}&0\cr
0&0&0&0&0&\hidew{-1}\cr
1&0&0&0&0&0\cr
0&1&0&0&0&0\cr
0&0&1&0&0&0\cr}
$$}

Let $B\subseteq G$ be the Borel subgroup which is the stabilizer of
the flag $\<v_1\>,\<v_1,v_2\>,\ldots$ where $v_i=\sum_{{i\over2}\le
j\le\noq+{i\over2}}e_j$. Then one can verify that $z\mapsto T(z)$ is
the embedding of $\fa^*$ into $V$ which is induced from $v^*=T^*$ and
the Borel subgroup $B$.  Thus, we get that the restricted zonal
spherical function $E_\lambda(T(z),T^*)$ is proportional to
$\Rq_\lambda(z;1)$.

The other interpretations of spherical functions are most interesting
when $n$ is odd. Then $H=H^*=Sp_{n-1}(\CC)$ and $A\mapsto
E_\lambda(A,T^*)$ is an $Sp_{n-1}(\CC)$\_invariant functions on the
space $V$ of skewsymmetric matrices of size $n+1$ which is a joint
eigenfunction for the differential operators $D_{(1^i)}$. Any
invariant function is uniquely determined by its value at $T(z)$ and
we have $E_\lambda(T(z),T^*)=\Rq_\lambda(z;1)$.

The open $G$\_orbit in $V$ is isomorphic to
$G/H=GL_n(\CC)/Sp_{n-1}(\CC)$. Thus, the pullback of an
$Sp_{n-1}(\CC)$\_invariant function on $V$ leads to a
$Sp_{n-1}(\CC)$\_biinvariant function on $GL_n(\CC)$. Clearly, not all
of them are of this form. For this, we have to make the function
$A\mapsto\|Pf|(A)$ invertible since its zero\_set is the complement of
the open orbit. Since $E_{(1^n)}(A,T^*)=\|Pf|(A)$, this corresponds on
$\fa^*$ to make the function $\pi(z):=\Rq_{(1^n)}(z)=\prod_{i\
\odd}z_i$ invertible. We can extend the definition of $\Rq_\lambda(z)$
to every element $\lambda\in\Lambda':=
\{\lambda\in\ZZ^n\mid\lambda_1\ge\ldots\ge\lambda_n\}$ by
$\Rq_\lambda:=\pi^{-m}\Rq_{\lambda+m(1^n)}$ for $m>\!\!\!>0$ (by
\cite{Recursion}, this is independent of the choice of $m$). Then
$\Rq_\lambda(z;1)$, $\lambda\in\Lambda'$ is the radial part of an
$Sp_{n-1}(\CC)$\_biinvariant function on $GL_n(\CC)$ which is a joint
eigenfunction for all $GL_n(\CC)$\_biinvariant differential
operators. A similar result holds if $GL_n(\CC)$ is replaced by
$SL_n(\CC)$. Then $\fa^*$ should be replaced by $\{z\in\fa^*\mid
\delta(z)=1\}$ and $\lambda$ should be an element of
$\Lambda'/\ZZ(1^n)$.

\medskip\noindent
{\it The case of\/ $G=GL_p(\CC)\times GL_1(\CC)$ acting on
$(\CC^p\otimes\CC)\oplus(\CC^p)^*$.}

\noindent The action of $(A,s)\in G$ on a pair of vectors $(u,v)$ is
$(sAu,(A^t)^{-1}v)$. Here $\Lambda_V$ is generated by
$\epsilon_1+\epsilon'$, $-\epsilon_p$, and $\epsilon'$ with degrees
$1$, $1$, and $2$, respectively. Thus, if we put
$$\eqno{}
e_1=-\epsilon_p,\ e_2=\epsilon_1+\epsilon_p+\epsilon',\ e_3=-\epsilon_1.
$$
then the generating weights become $e_1$, $e_1+e_2$, and
$e_1+e_2+e_3$. In particular, the degree of $e_i$ is $1$, $0$, $1$,
respectively. The little Weyl group is generated by the permutation
which swaps $\epsilon_1$ and $\epsilon_p$, and therefore $e_1$ and
$e_3$. Thus, the Capelli operators are described by semisymmetric
polynomials in $n=3$ variables. The vector
$\overline\rho=({p-1\over2},\ldots,-{p-1\over2};0)$ equals, up to a
$W_V$\_fixed vector
$({p-1\over2},0,\ldots,0,-{p-1\over2};{p-1\over2})={p-1\over2}(\epsilon_1-\epsilon_p+\epsilon)={p-1\over2}(2e_1+e_2)$.
This shows $r={p-1\over2}$.

The concrete decomposition of $\PP$ as a $G$\_module has been worked
out in \cite{VS}, see also \cite{VK} \S11.1--11.2. Here, we give only
the fundamental Capelli operators. Denote the coordinates of $V$ by
$u_1,\ldots,u_p;v_1,\ldots,v_p$. Then

$$\eqno{}
D_{(1)}=\sum_{i=1}^pv_i{\partial\over\partial v_i},\ 
D_{(11)}=\sum_{i=1}^pu_i{\partial\over\partial u_i},\ 
D_{(111)}=(\sum_i u_iv_i)
(\sum_i{\partial\over\partial u_i}{\partial\over\partial v_i}).
$$

The zonal spherical functions have been investigated by Vilenkin--\v
Sapiro~\cite{VS} (see also \cite{VK}~11.3.2). They are eigenfunctions for
the three differential operators $D_{(1)}$, $D_{(11)}$, and
$D_{(111)}$ above. For the Chevalley isomorphism we define
$$\eqno{}
\eqalign{u(z)&=(z_2,0,\ldots,0,z_2-z_3);\cr
v(z)&=(z_1-z_2+z_3,0,\ldots,0,-z_1+z_2);\cr
u^*_0&=(1,0,\ldots,0);\cr
v^*_0&=(1,0,\ldots,0).\cr}
$$
Then $z\mapsto(u(z),v(z))$ is the embedding $\fa^*\into V$ which
corresponds to $(u^*_0,v^*_0)\in V^*$ and the Borel subgroup which
stabilizes the flag
$\<e_1+e_n\>,\<e_1+e_n,e_2\>,\ldots,\<e_1+e_n,e_2,\ldots,e_n\>$.

To be eigenfunction for $D_{(1)}$ and $D_{(11)}$ simply means to be
bihomogeneous in the $u$- and $v$\_coordinates. Since
$H^*=GL_{n-1}(\CC)$ we get the following interpretation of zonal
spherical functions: they are bihomogeneous $GL_{n-1}(\CC)$\_invariant
functions on $\CC^n\oplus(\CC^n)^*$ which are eigenfunctions for the
``Laplace operator'' $D_{(111)}$.

This interpretation has also a real form: The complexification of
$U(n-1)$ is $H^*=GL_{n-1}(\CC)$ while $V$ is the complexification of
$\CC^n$, considered as an $\RR$\_vector space. The coordinate function
$v_i$ is then simply the complex conjugate $\uq_i$ of $u_i$. Thus,
$E_\lambda(u,\uq,u^*_0,v^*_0)$ is an $U(n-1)$\_invariant
function on $\CC^n$ which is bihomogeneous in the holomorphic and the
antiholomorphic variables and which is an eigenfunction of the (now
genuine) Laplace operator $D_{(111)}=\sum_i{\partial^2\over\partial
u_i\partial\uq_i}$. In this form, the $E_\lambda$ have been studied by
Vilenkin--\v Sapiro~\cite{VS}.

\beginsection DifferenceOperators. Difference operators

For $\lambda\in k^n$ let $T_\lambda$ be the shift operator $T_\lambda
f(z):=f(z-\lambda)$. Let $\epsilon_i$ be the $i$-th canonical basis
vector of $k^n$ and $T_i:=T_{\epsilon_i}$. For reasons of clarity we
adopt the following notation: $x_i:=z_{2i-1}$, $y_i:=z_{2i}$,
$T_{x,i}:=T_{2i-1}$, and $T_{y,i}:=T_{2i}$. Then we define the
following block matrices whose entries are difference operators (where
$t$ is an indeterminate):

$$\eqno{E9}
\fX(t):=\pmatrix{
\Big[
(x_i{+}t)(x_i{+}\r)^{\noq{-}j}-x_i^{\noq{+}1{-}j}T_{x,i}
\Big]_{i=1\ldots\noq\atop j=1\ldots\noq}&
\Big[
-x_i^{\noq{-}j}T_{x,i}
\Big]_{i=1\ldots\noq\atop j=1\ldots\nuq}\cr
\Big[
(y_i{+}\r)^{\nuq{+}1{-}j}-y_i^{\nuq{+}1{-}j}T_{y,i}
\Big]_{i=1\ldots\nuq\atop j=1\ldots\noq}&
\Big[
(y_i{+}\r)^{\nuq{-}j}
\Big]_{i=1\ldots\nuq\atop j=1\ldots\nuq}\cr}
$$
$$\eqno{E10}
\fY(t):=\pmatrix{
\Big[
(x_i{+}\r)^{\noq{-}j}
\Big]_{i=1\ldots\noq\atop j=1\ldots\noq}&
\Big[
(x_i{+}\r)^{\noq{-}j}-x_i^{\noq{-}j}T_{x,i}
\Big]_{i=1\ldots\noq\atop j=1\ldots\nuq}\cr
\Big[
-y_i^{\nuq{+}1{-}j}T_{y,i}
\Big]_{i=1\ldots\nuq\atop j=1\ldots\noq}&
\Big[
(y_i{+}t)(y_i{+}\r)^{\nuq{-}j}-y_i^{\nuq{+}1{-}j}T_{y,i}
\Big]_{i=1\ldots\nuq\atop j=1\ldots\nuq}\cr}
$$
The semisymmetric Vandermonde determinant is:
$$\eqno{}
\phi(z):=\prod_{1\le i<j\le n\atop j-i\ \even}(z_i-z_j)=
\prod_{1\le i< j\le\noq}(x_i-x_j)\prod_{1\le i< j\le\nuq}(y_i-y_j).
$$
Now we define the operators
$$\eqno{E90}
X(t):=\phi(z)^{-1}\|det|\fX(t)\hbox{ and }
Y(t):=\phi(z)^{-1}\|det|\fY(t).
$$
First observe, that the entries of $\fX(t)$ and $\fY(t)$ commute if
they are in different rows. Thus, the determinants are well
defined.

\Lemma. Both $X(t)$ and $Y(t)$ act on $\cP^W$.

\Proof: Let $f\in\cP^W$. Then both $\fX(t)f$ and $\fY(t)f$ are
polynomials which are skewsymmetric with respect to both factors
$S_\noq$, $S_\nuq$ of $W$. Therefore, they are divisible by $\phi(z)$
and the quotient is $W$\_symmetric.\qed

\Lemma degree. For $f\in\cP^W$ holds $\|deg|X(t)f\le\|deg|f$ and
$\|deg|Y(t)f\le\|deg|f$.

\Proof: We use the following elementary fact: let $A=(a_{ij})$ a matrix
with entries in a filtered ring such that $\|deg|a_{ij}\le d_i'+d_j''$ for
integers $d_i'$ and $d_j''$. Then
$\|deg|\|det|A\le\sum_i(d_i'+d_i'')$.

We apply this to $\fX(t)$. Using that the operator $1-T_i$ has degree
$-1$ the entries of $\fX(t)$ have degree $d_i'+d_j''$ with
$\epsilon:=\noq-\nuq$ and
$$\eqno{E15}
\vcenter{\halign{
$#=$\hfill&
(\hfil$#$\hfil$,\ldots,$&
\hfil$#$\hfil$,\,$&
\hfil$#$\hfil$,\,$&
\hfil$#$\hfil$,\,$&
\hfil$#$\hfil$,\ldots,\,$&\hfil$#$\hfil$)$\cr
(d_1',\ldots,d_n')&0     &0           &0&-\epsilon&-\epsilon&-\epsilon\cr
(d_1'',\ldots,d_n'')  &\noq-1&           1&0 &\noq-1&\noq-2      &\epsilon\cr
}}
$$
Thus $\|deg|\|det|\fX(t)\le\sum_i(d_i'+d_i'')=\|deg|\phi(z)$.

For $\fY(t)$ one argues in the same way with
$$\eqno{E16}
\vcenter{\halign{
$#=$\hfill&
(\hfil$#$\hfil$,\ldots,$&
\hfil$#$\hfil$,\,$&
\hfil$#$\hfil$,\,$&
\hfil$#$\hfil$,\,$&
\hfil$#$\hfil$,\ldots,\,$&$#)$\cr
(d_1',\ldots,d_n')&0&0&0   &1-\epsilon&1-\epsilon&1-\epsilon\cr
(d_1'',\ldots,d_n'')  &\noq-1&1&0&\noq-2&\noq-3     &\epsilon-1\cr
}}
$$\qed

Next we derive an explicit formula for $X(t)$ and $Y(t)$.
For
$I\subseteq\{1,\ldots,n\}$ put $\epsilon_I:=\sum_{i\in I}\epsilon_i$ and
$T_If(z):=f(z-\epsilon_I)$. Thus $T_I=\prod_{i\in I}T_i$. With this
notation we define
$$\eqno{E13}
D_I:=\prod_{i\in I\atop n-i\ \rm even}\!\!\!\!z_i
\prod_{i\in I, j\not\in I\atop j-i\ \odd}\!\!(z_i-z_j-\r)
\prod_{i\in I, j\not\in I\atop j-i\ \even}\!\!(z_i-z_j)^{-1}\ T_I.
$$

Let $P_\odd$ be the set of subsets $I\subseteq\{1,\ldots,n\}$ such
that there is a $w\in W$ with $I=w\{1,\ldots,m\}$ (where $m=|I|$). Thus
$I\in P_\odd$ if and only if the number of its odd members is equal or one
more than the number of its even members. Let $P_\even\subseteq
P_\odd$ consist of those sets where these numbers are equal. This is
equivalent to $|I|$ being even. Finally, we set
$$\eqno{}
|I|_{\|o|}:=|\{i\in I\mid i\hbox{ odd}\,\}|=
\left\lceil{\textstyle|I|\over2}\right\rceil=
|\epsilon_I|_\odd.
$$

\Proposition Expansion. We have 
$$
\eqalignno{
E44&X(t)&=\sum_{I\in P_\odd}(-1)^{|I|_{\|o|}}\prod_{i\not\in I\atop
i\ \rm odd}(t+z_i)\ D_I\cr
E45&Y(t)&=\sum_{I\in P_\even}(-1)^{|I|_{\|o|}}\prod_{i\not\in I\atop
i\ \rm even}(t+z_i)\ D_I\cr}
$$

\Proof: Clearly both $X(t)$ and $Y(t)$ have expansions of the form
$\sum_Ic_IT_I$ where $I$ runs through all subsets of
$\{1,\ldots,n\}$. First we show that in $Y(t)$ only the $I\in P_\even$
contribute. For this, we subtract in \cite{E10} the columns $1,\ldots,\nuq$
from the columns $\noq+1,\ldots,\noq+\nuq=n$, respectively and obtain
$$\eqno{E11}
\|det|\fY(t)=\|det|\pmatrix{
\Big[
(x_i{+}\r)^{\noq{-}j}
\Big]_{i=1\ldots\noq\atop j=1\ldots\noq}&
\Big[
-x_i^{\noq{-}j}T_{x,i}
\Big]_{i=1\ldots\noq\atop j=1\ldots\nuq}\cr
\Big[
-y_i^{\nuq{+}1{-}j}T_{y,i}
\Big]_{i=1\ldots\nuq\atop j=1\ldots\noq}&
\Big[
(y_i{+}t)(y_i{+}\r)^{\nuq{-}j}
\Big]_{i=1\ldots\nuq\atop j=1\ldots\nuq}\cr}
$$
Next, we use the well known fact that the determinant of a block
matrix is
$$\eqno{}
\|det|\pmatrix{A&B\cr C&D\cr}=\|det|(A-BD^{-1}C)\|det|D
$$
if $D$ is invertible. In our case, the entries of $A-BD^{-1}C$ are
linear combinations (over the field of rational functions in $z$) of
$1$ and the $T_{x,i}T_{y,j}$ which proves the claim.

Since the given form of the operator $Y(t)$ is $W$\_symmetric it
suffices to check the coefficient of $T_I=T_{x,1}\ldots
T_{x,l}T_{y,1}\ldots T_{y,l}$ where $I=\{1,\ldots,2l\}$,
$l=0,\ldots,\nuq$. Every entry of the matrix \cite{E11} is of the form
$a+bT_{x,i}$ or $a+bT_{y,i}$. Thus the required coefficient is the
determinant of the matrix obtained by replacing that entry by $b$ if
$i\le l$ or $a$ if $i>l$. The ensuing matrix has the following form with
dimensions as indicated:
$$\eqno{}
\pmatrix{0_{l\times\noq}&*_{l\times\nuq}\cr*_{\noq{-}l\times\noq}&
0_{\noq{-}l\times\nuq}\cr*_{l\times\noq}&0_{l\times\nuq}\cr0_{\nuq{-}l\times\noq}&*_{\nuq{-}l\times\nuq}\cr}
$$
We interchange the first with the third block of rows. Then the
determinant gets multiplied by $(-1)^l=(-1)^{|I|_{\|o|}}$ and the
matrix acquires block diagonal form. The blocks are, up to a common
factor in each row, in Vandermonde form. Thus the formula given in the
Theorem is easily established.

The case of $X(t)$ is similar but a bit more complicated. Here we can
subtract in \cite{E9} column number $\noq+1$ through $2\noq-1$ from columns $2$
through $\noq$, respectively. Then we obtain
$$\eqno{E12}
\|det|\fX(t)=\|det|
\pmatrix{
[u_i]_{i=1\ldots\noq}&\Big[
(x_i{+}t)(x_i{+}\r)^{\noq{-}j}
\Big]_{i=1\ldots\noq\atop j=2\ldots\noq}&
\Big[
-x_i^{\noq{-}j}T_{x,i}
\Big]_{i=1\ldots\noq\atop j=1\ldots\nuq}\cr
[v_i]_{i=1\ldots\nuq}&\Big[
-y_i^{\nuq{+}1{-}j}T_{y,i}
\Big]_{i=1\ldots\nuq\atop j=2\ldots\noq}&
\Big[
(y_i{+}\r)^{\nuq{-}j}
\Big]_{i=1\ldots\nuq\atop j=1\ldots\nuq}\cr}
$$
with $u_i:=(x_i{+}t)(x_i{+}\r)^{\noq{-}1}-x_i^\noq T_{x,i}$ and
$v_i:=(y_i{+}\r)^\nuq-y_i^\nuq T_{y,i}$. Arguing as above, one notes
that all entries of $A-BD^{-1}C$ are linear combinations of $1$ and
$T_{x,i}T_{y,j}$ except those in the first column where also
$T_{x,i}$ appears. Thus, if $c_I\ne0$ then the number of odd elements
is equal or one more than the number of even elements, i.e., $I\in
P_\odd$.

To determine the correct coefficient we proceed as above. The case
$I=\{1,\ldots,2l\}$ is the same. In the case $I=\{1,\ldots,2l+1\}$ one
has to move in \cite{E12} the first column to the $\noq+1$-st place
(i.e., between the two other blocks).\qed

The main feature of the difference operators is the following cut-off
property:

\Lemma Cut-off. Expand $X(t)$ or $Y(t)$ as
$\sum_Ic_I(z)T_I$. Assume $r\ne0$. Then for any $\mu\in\Lambda$ holds:
if $\mu-\epsilon_I\not\in\Lambda$ then $c_I(\rho+\mu)=0$.

\Proof: Since $r\ne0$ (and, as always, $\rho$ dominant), the
denominator of $c_I$ does not vanish at $\rho+\mu$. If
$\mu-\epsilon_I\not\in\Lambda$ then either $\mu_n=0$ and $n\in I$ or
there is an $i<n$ with $\mu_i=\mu_{i+1}$ and $i\in I$, $i+1\not\in
I$. Now we use the precise form of $c_I(z)$ established in
\cite{Expansion}. From the definition of $D_I$, \cite{E13}, we obtain
that $c_I$ is a multiple of $z_n(z_i-z_{i+1}-\r)$, hence
$c_I(\rho+\mu)=0$.\qed

Combining all results, we obtain the main result of this paper.

\Theorem Main. Every $R_\lambda$, $\lambda\in\Lambda$, is an
eigenvector of both $X(t)$ and $Y(t)$. More precisely:
$$
\eqalignno{
E50&X(t)R_\lambda&=
{\textstyle\prod\limits_{i\ \odd}}(t+\rho_i+\lambda_i)\cdot R_\lambda\cr
E51&Y(t)R_\lambda&=
{\textstyle\prod\limits_{i\ \even}}(t+\rho_i+\lambda_i)\cdot R_\lambda\cr}
$$

\Proof: We may assume $r\ne0$. The case $r=0$ then follows by continuity.
Let $R:=X(t)R_\lambda$. \cite{degree} implies
$\|deg|R\le\|deg|R_\lambda=|\lambda|_\odd$. Let $\mu\in\Lambda$ with
$|\mu|_\odd\le|\lambda|_\odd$ and $\mu\ne\lambda$. If
$X(t)=\sum_Ic_I(z)T_I$ then $R(\rho+\mu)=
\sum_Ic_I(\rho+\mu)R_\lambda(\rho+(\mu-\epsilon_I))$.
If $\mu-\epsilon_I\in\Lambda$ then
$R_\lambda(\rho+(\mu-\epsilon_I))=0$ by definition of
$R_\lambda$. Otherwise, $c_I(\rho+\mu)=0$ by \cite{Cut-off}. Hence,
$R(\rho+\mu)=0$ which shows that $R$ is a multiple of $R_\lambda$. The
coefficient of $T_\emptyset$ in $X(t)$ is
$c(z)=\prod_{i=1}^\noq(t+z_{2i-1})$. Thus, evaluation in $z=\rho+\lambda$
gives $R=c(\rho+\lambda)R_\lambda$. The same argument works for $Y(t)$.\qed

\Corollary Commute. Let
$$\eqno{E17}
\eqalign{
X(t)&=t^\noq+X_1t^{\noq-1}+\ldots+X_\noq,\cr
Y(t)&=t^\nuq+Y_1t^{\nuq-1}+\ldots+Y_\nuq.\cr}
$$
Then $X_1,\ldots,X_\noq,Y_1,\ldots,Y_\nuq$ are pairwise commuting
difference operators.

\Example: We compute $X_1$ and $Y_1$. Any contribution to the
coefficient of $t^{\noq-1}$ in \cite{E44} comes from $I=\emptyset$,
$I=\{i\}$ with $i$ odd, and $I=\{i,j\}$ with $i$ odd, $j$ even. Since
$D_\emptyset=1$ we get
$$\eqno{E46}
X_1=\sum_{i\ \odd}z_i-\sum_{i\ \odd}D_{\{i\}}-
\sum_{i\ \odd\atop j\ \even}D_{\{i,j\}}.
$$
Similarly, in \cite{E45}, the contribution for $t^{\nuq-1}$ comes from
$I=\emptyset$ and $I=\{i,j\}$ with $i$ odd, $j$ even. Thus,
$$\eqno{E47}
Y_1=
\sum_{i\ \even}z_i-\sum_{i\ \odd\atop j\ \even}D_{\{i,j\}}.
$$

\medskip
The operators $X_1,\ldots,X_\noq,Y_1,\ldots,Y_\nuq$ defined in \cite{E17}
generate a polynomial algebra $\cR\subseteq\|End|\cP^W$. We show that
it is canonically isomorphic to $\cP^W$. More precisely:

\Proposition Coeffs. 
\item{a)}Every element $D\in\cR$ has an expansion
$$\eqno{E18}
D=\sum_{\mu\in\Psi_0} c_\mu^D(z)T_\mu,
$$
where $\Psi_0$ is the smallest $W$\_stable submonoid of $\ZZ^n$
containing $\Lambda$ and where the coefficients $c_\mu^D$ are rational
functions in $z_1,\ldots,z_n$ with poles along the hyperplanes
$z_i-z_j=a$ where $i-j\ne0$ is even and $a\in\ZZ$.
\item{b)} The coefficient $c_0^D(z)$ is in $\cP^W$ and the map
$\cR\pfeil\cP^W:D\mapsto c_0^D(z)$ is an algebra
isomorphism.
\item{c)}For every $D\in\cR$ and $\lambda\in\Lambda$
holds $D(R_\lambda)=c_0^D(\rho+\lambda)R_\lambda$.
\item{d)}If $r\not\in\QQ$ the coefficients $c_\mu^D(z)$ have the
cut-off property: Let $\lambda\in\Lambda$ with
$\lambda-\mu\not\in\Lambda$. Then $c_\mu^D(\rho+\lambda)=0$.\Par

\Proof: a) Since $\Psi_0$ contains all $W$\_translates of elements of
$\Lambda$ it contains all $\epsilon_I$ with $I\in P_\odd$. Hence, by
\cite{Expansion}, the generators $X_1,\ldots,X_\noq,Y_1,\ldots,Y_\nuq$
of $\cR$ have an expansion as claimed. This implies the result easily
for all $D\in\cR$.

b) Since $\Psi_0$ is contained in $\NN^n$ it is a pointed cone, i.e.,
$\lambda,\mu\in\Psi_0$ with $\lambda+\mu=0$ implies
$\lambda=\mu=0$. Looking at how two operators with an expansion as in
\cite{E18} multiply this implies that $D\mapsto c_0^D$ is an algebra
homomorphism. It is an isomorphism, since the generators $X_i$ and
$Y_i$ of $\cR$ are mapped to free generators of $\cP^W$.

c) The assertion needs to be checked just for the generators of $\cR$
and there it is the content of \cite{Main}.

d) The algebra $\cR$ acts on the dual space $(\cP^W)^*$ on the right
by $(\delta D)(f):=\delta(D(f))$. Let $v=(v_i)\in k^n$ such that
$v_i-v_j\not\in\ZZ$ whenever $i-j\ne0$ is even
(e.g. $v\in\rho+\Lambda$, since $r\not\in\QQ$) and let
$\delta_v:f\mapsto f(v)$ be the corresponding evaluation
function. Then $\delta_vD=\sum_\mu c_\mu^D(v)\delta_{v-\mu}$. Thus the
cut-off property is equivalent to the statement that
$\oplus_{\lambda\in\Lambda}k\delta_{\rho+\lambda}\subseteq(\cP^W)^*$
is $\cR$\_stable. It suffices to check this for generators of $\cR$
which is the content of \cite{Cut-off}.\qed

Next we study the monoid $\Psi_0$ more closely.

\Lemma Psi0. The monoid $\Psi_0$ has also the following descriptions:
\item{a)}It is generated by $\{\epsilon_i\mid i\ \odd\}\cup
\{\epsilon_i+\epsilon_j\mid i\ \odd, j\ \even\}$.
\item{b)}It consists of all $\lambda\in\NN^n$ with $[\lambda]_1\ge0$.

\Proof: a) Since $\Lambda$ is generated by all
$\epsilon_{\{1,\ldots,m\}}$, $m=1,\ldots,n$, the monoid $\Psi_0$ is
generated by all $\epsilon_I$, $I\in P_\odd$. But those can be
obtained from the given subset.

b) Let $\Psi_0'$ be the set of all $\lambda\in\NN^n$ with
$[\lambda]_1\ge0$, i.e., $|\lambda_\odd|\ge|\lambda_\even|$. We have
to show $\Psi_0=\Psi_0'$. The inclusion $\Psi_0\subseteq\Psi_0'$
follows, e.g., from a). For the converse, let $\lambda\in\Psi_0$. We
show $\lambda\in\Psi_0'$ by induction on $|\lambda|_\odd$. If
$|\lambda|_\odd=0$ then also $|\lambda_\even|=0$. Thus
$\lambda=0\in\Psi_0$. For $|\lambda|_\odd>0$ there are two cases. If
$|\lambda_\odd|>|\lambda_\even|$ then choose any odd $i$ such that
$\lambda_i>0$. Then $\lambda':=\lambda-\epsilon_i$ is in $\Psi_0'$
hence, by induction, in $\Psi_0$. Thus also
$\lambda=\lambda'+\epsilon_i\in\Psi_0$. If
$|\lambda_\odd|=|\lambda_\even|>0$ then there is $i$ odd and $j$ even
such that $\lambda_i>0$ and $\lambda_j>0$. Then
$\lambda':=\lambda-\epsilon_i-\epsilon_j$ is in $\Psi_0'$, hence in
$\Psi_0$ by induction. We conclude $\lambda\in\Psi_0$, as well.\qed

In the theory of symmetric polynomials, the containment relation
$\lambda\subseteq\mu$ for $\lambda,\mu\in\Lambda$ is defined as
$\mu-\lambda\in\NN^n$. The semisymmetric analogue is
$\mu-\lambda\in\Psi_0$ or, equivalently,
$$\eqno{}
\lambda\LE\mu\quad {\buildrel{\rm def}\over\Longleftrightarrow}\quad
\lambda\subseteq\mu\hbox{\ and\ }[\lambda]_1\le[\mu]_1.
$$

\Example: We have $(1,0,0)\subseteq(1,1,0)$ but
$(1,0,0)\not\LE(1,1,0)$. Moreover, we have
$$\eqno{}
(1,0,0),(1,1,0)\LE (2,1,0),(1,1,1).
$$
This implies in particular that $(1,0,0)$ and $(1,1,0)$ have no
supremum. Therefore, as
opposed to the classical containment relation, its semisymmetric
analogue does not form a lattice.  \medskip

Now we prove that the polynomial $R_\lambda$ vanishes at many more
points than it is supposed to by definition (Extra Vanishing Theorem).

\Theorem ExtraVanishing. For $\lambda,\mu\in\Lambda$ holds
$R_\lambda(\rho+\mu)=0$ unless $\lambda\LE\mu$.

\Proof: We may assume $r\not\in\QQ$ since the general case follows by
continuity. For fixed $\lambda$ let $\mu$ be a counterexample (i.e.,
$R_\lambda(\rho+\mu)\ne0$ and $\lambda\not\LE\mu$) which is minimal
with respect to ``$\LE$''. Since $\rho+\mu$ is not in the $W$\_orbit
of $\rho+\lambda$ there is $D\in\cR$ such that $c_0^D(\rho+\lambda)\ne
c_0^D(\rho+\mu)$. From $D(R_\lambda)=c_0^D(\rho+\lambda)R_\lambda$ we
obtain, after substituting $z=\rho+\mu$:
$$\eqno{}
(c_0^D(\rho+\lambda)-c_0^D(\rho+\mu))R_\lambda(\rho+\mu)=
\sum_{\eta\in\Psi_0}c_\eta^D(\rho+\mu)R_\lambda(\rho+\mu-\eta)
$$
If $\mu-\eta\in\Lambda$ then $R_\lambda(\rho+\mu-\eta)=0$ by
minimality of $\mu$. Otherwise $c_\eta^D(\rho+\mu)=0$ by
\cite{Coeffs} d). Contradiction.\qed

\noindent As an application we derive an explicit formula for
$R_\lambda$ when $\lambda$ is a particular kind of ``hook''.

\Corollary Hook1. Let $a,m\ge1$ be integers with $m$ odd. Then
$$\eqno{E30}
R_{(a\,1^{m-1})}=(R_{(1)}-1)(R_{(1)}-2)\ldots(R_{(1)}-a+1)R_{(1^m)}.
$$

\Proof: Denote the right-hand side by $f$. Let $\lambda=(a\,1^{m-1})$ and
$\mu\in\Lambda$ with $|\mu|_\odd\le|\lambda|_\odd=a+{m-1\over2}$ and
$f(\rho+\mu)\ne0$. Then $R_{(1^m)}(\rho+\mu)\ne0$ and
$R_{(1)}(\rho+\mu)\ne1,2,\ldots,a-1$. The Extra Vanishing
Theorem~\ncite{ExtraVanishing} implies $(1^m)\LE\mu$, hence $\mu_m\ge1$
and $[\mu]_1\ge[(1^m)]_1=1$. From $R_{(1)}(\rho+\mu)=[\mu]_1$
(\cite{Elementary}) we obtain $[\mu]_1\ge a$.  Thus
$\mu_1=[\mu]_1+(\mu_2-\mu_3)+\ldots\ge a$. We know already
$\mu_3,\mu_5,\ldots,\mu_m\ge1$. Since $|\mu|_\odd\le a+{m-1\over2}$,
equality holds throughout. This implies easily
$\mu=\lambda$. Therefore, $f$ must be a multiple of $R_\lambda$. Equality
follows from the fact that the coefficient of
$z^{[\lambda]}=z_1^az_3z_5\ldots z_m$ is $1$ in both cases.\qed

\Remarks: 1. For even $m$, the polynomials $R_{(a\,1^{m-1})}$ will be
calculated in \cite{Hook2}.

\noindent
2. For $m=1$ one obtains in particular
$R_{(a)}=R_{(1)}(R_{(1)}-1)\ldots(R_{(1)}-a+1)$ which clearly do not generate
$\cP^W$. Therefore, the polynomials $R_{(a)}$ are not a semisymmetric
analogue of the complete symmetric functions.

\beginsection TopHomogeneous. The top homogeneous components

The highest degree components $\Rq_\lambda(z;r)$ of $R_\lambda(z;r)$
are also of high representation theoretic interest (see
\cite{SpherFunc}). We show that they are eigenfunctions of {\it
differential\/} equations. More precisely, put
$$\eqno{}
\overline\fX(t):=\pmatrix{
\Big[
x_i^{\noq{-}j}\big(x_i\partial_{x_i}{+}(\noq{-}j)\r{+}t\big)
\Big]_{i=1\ldots\noq\atop j=1\ldots\noq}&
\Big[
-x_i^{\noq{-}j}
\Big]_{i=1\ldots\noq\atop j=1\ldots\nuq}\cr
\Big[
y_i^{\nuq{-}j}\big(y_i\partial_{y_i}{+}(\nuq{+}1{-}j)\r\big)
\Big]_{i=1\ldots\nuq\atop j=1\ldots\noq}&
\Big[
y_i^{\nuq{-}j}
\Big]_{i=1\ldots\nuq\atop j=1\ldots\nuq}\cr}
$$
$$\eqno{}
\overline\fY(t):=\pmatrix{
\Big[
x_i^{\noq{-}j}
\Big]_{i=1\ldots\noq\atop j=1\ldots\noq}&
\Big[
x_i^{\noq{-}j{-}1}\big(x_i\partial_{x_i}{+}(\noq{-}j)\r\big)
\Big]_{i=1\ldots\noq\atop j=1\ldots\nuq}\cr
\Big[
-y_i^{\nuq{+}1{-}j}
\Big]_{i=1\ldots\nuq\atop j=1\ldots\noq}&
\Big[
y_i^{\nuq{-}j}\big(y_i\partial_{y_i}{+}(\nuq{-}j)\r{+}t\big)
\Big]_{i=1\ldots\nuq\atop j=1\ldots\nuq}\cr}
$$
and $\overline X(t):=\phi(z)^{-1}\|det|\overline\fX(t)$, $\overline
Y(t):=\phi(z)^{-1}\|det|\overline\fY(t)$ where
$\partial_{x_i}=\partial/\partial x_i$ and
$\partial_{y_i}=\partial/\partial y_i$. These are linear differential
operators with rational coefficients.

\Theorem. Every $\Rq_\lambda$, $\lambda\in\Lambda$, is an
eigenvector of both $\overline X(t)$ and $\overline Y(t)$. More precisely:
$$
\eqalignno{
E48&\overline X(t)\Rq_\lambda&=
\prod_{i\ \odd}(t+\rho_i+\lambda_i)\,\Rq_\lambda\cr
E49&\overline Y(t)\Rq_\lambda&=
\prod_{i\ \even}(t+\rho_i+\lambda_i)\,\Rq_\lambda\cr}
$$

\Proof: Let $f\in\cP^W$ is homogeneous of degree $d$. For each
entry $a_{ij}$ of $\fX(t)$ or $\fY(t)$ we know that $a_{ij}(f)$ is a
polynomial of degree $\le d+d_i'+d_j''$ (with $d_i'$, $d_j''$ as in
\cite{E15} or \cite{E16}). Using the fact that
$$\eqno{}
\eqalign{
(1-T_{x,i})f(z)=&\partial_{x_i}f+\hbox{ lower order terms}\cr
(1-T_{y,i})f(z)=&\partial_{y_i}f+\hbox{ lower order terms}\cr}
$$
one easily calculates that the $d+d_i'+d_j''$\_degree component of
$a_{ij}(f)$ is $\aq_{ij}(f)$ where $\aq_{ij}$ is the $ij$\_entry of
$\overline\fX(t)$ or $\overline\fY(t)$, respectively. Since
$\sum_i(d_i'+d_i'')=\|deg|\phi(z)$, we get that $\Xq(t)f$ or $\Yq(t)f$ is the
$d$\_degree homogeneous component of $X(t)f$ or $Y(t)f$,
respectively. Now the assertion follows from \cite{Main}.\qed

\Remark: The operators $\Xq(t)$ and $\Yq(t)$ are the semisymmetric
analogues of the Sekiguchi-Debiard operators, \cite{Seki}, \cite{Deb},
which characterize Jack polynomials.
\medskip

If we expand $\Xq(t)$ and $\Yq(t)$ as a polynomial in $t$,
$$\eqno{E42}
\eqalign{
\Xq(t)&=t^\noq+\Xq_1t^{\noq-1}+\ldots+\Xq_\noq,\cr
\Yq(t)&=t^\nuq+\Yq_1t^{\nuq-1}+\ldots+\Yq_\nuq,\cr}
$$
we obtain as in \cite{Commute} pairwise commuting differential
operators $\Xq_1,\ldots,\Xq_\noq$, $\Yq_1,\ldots$, $\Yq_\nuq$ with
$\Rq_\lambda$ as common eigenvectors. In general, these operators seem
to be more difficult to compute explicitly than their difference
counterparts. We give a formula for the most important ones, namely
those of order one. For odd $i$ we define the following rational
function:
$$\eqno{}
u_i:=v_i\,{\prod\limits_{j\ \even}(z_i-z_j)\over\prod\limits_{j\ne
i\ \odd}(z_i-z_j)}\hbox{ where }
v_i:=\cases{z_i&for $n$ \odd\cr 1&for $n$ \even\cr}
$$

\Theorem. The following equations hold:
$$
\eqalignno{
E52&&\eta:=\Xq_1-\noq\nuq r
=\sum_iz_i{\partial\over\partial z_i}\qquad\hbox{(Euler
vector field)}\cr
E53&&\eta':=\Xq_1-\Yq_1-\nuq r=\sum_{i\ \odd}u_i{\partial\over\partial z_i}\cr}
$$
Moreover, for all $\lambda\in\Lambda$ holds
$\eta(\Rq_\lambda)=|\lambda|_\odd\Rq_\lambda$,
$\eta'(\Rq_\lambda)=[\lambda]_1\Rq_\lambda$.

\Proof: Let $E$ be the Euler vector field. By \cite{E48}, we have
$\Xq_1(\Rq_\lambda)=|\rho+\lambda|_\odd\Rq_\lambda$. From
$|\lambda|_\odd=\|deg|\Rq_\lambda$ and $|\rho|_\odd=\noq\nuq r$ it
follows that $\eta-E$ kills every $\Rq_\lambda$ and therefore every
$W$\_invariant. The (non-symmetric) polynomials are all algebraic
functions of the semisymmetric ones. Since $\eta-E$ is a derivation,
it kills all polynomials, i.e., $\eta-E=0$.

By \cite{E46} and \cite{E47}, we have
$$\eqno{}
E':=X_1-Y_1=[z]_1-\sum_{i\ \odd}D_{\{i\}}=[z]_1-\sum_{i\ \odd}u_i'T_i
$$
where
$$\eqno{}
u_i':=v_i\,{\prod\limits_{j\ \even}(z_i-z_j-r)\over\prod\limits_{j\ne
i\ \odd}(z_i-z_j)}.
$$
From \cite{E50} and \cite{E51} we obtain $E'(1)=[\rho]_1=\nuq
r$. Thus, $[z]_1-\sum_{i\ \odd}u_i'=\nuq r$ and we get
$$\eqno{}
X_1-Y_1-\nuq r=\sum_{i\ \odd}u_i'(1-T_i).
$$
This implies \cite{E53} since
$(1-T_i)(f)={\partial f\over\partial
z_i}+\hbox{lower order terms}$.\qed

\noindent The derivations $\eta$, $\eta'$ induce a bigrading on
$\cP^W$. More precisely, for integers $a,b$ let
$$\eqno{}
\cP_{a,b}^W:=\{f\in\cP^W\mid \eta(f)=af,\eta'(f)=bf\}.
$$
Then $\cP^W=\oplus_{a\ge b\ge0}\cP_{a,b}^W$. To describe $\cP_{a,b}^W$
explicitly, we have to find bihomogeneous generators of $\cP^W$. For
this, we introduce the semi\_symmetric analogue of the elementary
symmetric polynomials, namely $\e_m(z):=\Rq_{(1^m)}(z;r)$. More
explicitly we have by \cite{Elementary2}:
$$\eqno{}
\vcenter{\halign{$#$\hfill&$#$,\quad\hfill&$#$\hfill\cr
\e_{2m{-}1}(z)&=e_m(z_\odd)-e_m(z_\even)&m=1,\ldots,\noq;\cr
\e_{2m}(z)&=e_m(z_\even)&m=1,\ldots,\nuq.\cr}}
$$
Now, we consider the basis of $\cP^W$ which consists of all monomials
in the $\e_m$. More precisely, we define\footnote*{In this notation,
one has to be careful to distinguish between $\e_a$ (the index is a
number) and $\e_{(a)}$ (the index is a partition). The latter equals
$\e_1^a$.} for any $\lambda\in\Lambda$
$$\eqno{E60}
\e_\lambda:=\e_1^{\lambda_1-\lambda_2}\e_2^{\lambda_2-\lambda_3}
\ldots \e_{n-1}^{\lambda_{n-1}-\lambda_n}\e_n^{\lambda_n}.
$$
This parametrization
is chosen such that the leading term of
$\e_\lambda$ is $z^{[\lambda]}$. Then $\cP_{a,b}^W$ is spanned by all
$\e_\lambda$ with $|\lambda|_\odd=a$ and $[\lambda]_1=b$.

\Corollary Bigrading. For $\lambda\in\Lambda$ consider the expansion
$\Rq_\lambda=\sum_\mu a_{\lambda\mu}\e_\mu$.  Then only those $\e_\mu$
occur for which $|\mu|_\odd=|\lambda|_\odd$ and $[\mu]_1=[\lambda]_1$.

\Remark: This result will be generalized in \cite{Triangular2}.
\medskip
We use \cite{Bigrading} to compute $\Rq_\lambda$ for all two\_row
diagrams. We use the multinomial coefficient ${a\choose
k_1,\ldots,k_n}:={a!\over k1!\ldots k_n!}$ where $a=k_1+\ldots+k_n$.

\Theorem Tworow. For integers $a\ge b\ge0$ let $c_{ab}={a\choose
b}{-2r\choose a}$. Then
$$\eqno{E40}
\Rq_{(a\,b)}={1\over c_{ab}}\sum_\mu
{-2r\choose\mu_1}
{\mu_1\choose\mu_1{-}\mu_2,\ldots,\mu_{n{-}1}{-}\mu_n,\mu_n}\e_\mu
$$
where the sum runs through all $\mu\in\Lambda$ with $|\mu_\odd|=a$ and
$|\mu_\even|=b$.

\Proof: We use a result for the usual Jack polynomials
$\Pq_\lambda(z;r)$. Stanley (\cite{St}, see also \cite{KnSa} Prop.~3.4
for a proof in the spirit of this paper), has shown that there is a
generating series
$$\eqno{}
\sum_{a=0}^\infty v_a \Pq_{(a)}(z;r)=\prod_i(1+z_i)^{-r}
$$
where the $v_a={-r\choose a}\ne0$. Then, by the
Comparison \cite{Jack2} (with $\mu=(a,0,\ldots)$) there are constants
$w_{a,b}\ne0$ such that
$$\eqno{}
\sum_{a\ge b\ge0} w_{a,b}\Rq_{(a\,b)}(z)=\prod_{i\ \odd}(1+z_i)^{-2r}.
$$
Now, we expand the right-hand side in bihomogeneous components. For
this observe
$$\eqno{}
\prod_{i\ \odd}(1+z_i)=1+\sum_{i\ge1}e_i(z_\odd)=1+\sum_{i\ge1}\e_i(z).
$$
Thus, we get
$$\eqno{}
\eqalign{
\sum_{a\ge b\ge0} w_{a,b}\Rq_{(a\,b)}(z)&=
\sum_d{-2r\choose d}\bigg(\sum_{i\ge1}\e_i\bigg)^d=\cr
&=\sum_d\sum_{k_1+\ldots+k_n=d}{-2r\choose d}{d\choose k_1,\ldots,k_n}
\e_1^{k_1}\e_2^{k_2}\ldots\e_n^{k_n}\cr
&=\sum_{\mu\in\Lambda}
{-2r\choose \mu_1}{\mu_1\choose \mu_1{-}\mu_2,\ldots,\mu_{n{-}1}{-}\mu_n,\mu_n}\e_\mu\cr}.
$$
Now, we compare the bihomogeneous components of bidegree $(a,a-b)$ of
both sides and get formula \cite{E40} up to the scalar $c_{ab}$. But
that scalar is easily obtained by the requirement that the coefficient
of $\e_{(a\,b)}$ should be $1$.\qed

\Examples: 1. {\it The case $n=3$.} The summation in \cite{E40}
runs through all $\mu\in\Lambda_3$ with $\mu_1+\mu_3=a$ and
$\mu_2=b$. If we put $\mu_3=k$, we get $\mu=(a-k,b,k)$ with $0\le
k\le\|min|(a-b,b)$. Thus,
$$\eqno{}
\Rq_{a,b,0}=\sum_k{{a-k\choose a-b-k,b-k,k}{-2r\choose a-k}\over
{a\choose b}{-2r\choose a}}\e_{a-k,b,k}=\sum_k(-1)^k
{{a-b\choose k}{b\choose k}\over{a+2r-1\choose k}}\e_{a-k,b,k}
$$
Now, the recursion formula \cite{E32} implies
$\Rq_{a+1,b+1,c+1}=\e_3\Rq_{a,b,c}$. Thus, we obtain a formula
for $\Rq_\mu$ for arbitrary $\mu\in\Lambda_3$:
$$\eqno{E41}
\Rq_{\mu_1,\mu_2,\mu_3}=\sum_k(-1)^k
{{\mu_1-\mu_2\choose k}{\mu_2-\mu_3\choose
k}\over{\mu_1-\mu_3+2r-1\choose k}}
\e_1^{\mu_1-\mu_2-k}\e_2^{\mu_2-\mu_3-k}\e_3^{\mu_3+k}
$$
with
$$\eqno{}
\e_1=z_1-z_2+z_3,\quad\e_2=z_2,\quad\e_3=z_1z_3.
$$
This formula can be rewritten in two ways. First, as a
hypergeometric function:
$$\eqno{}
\Rq_{\mu_1,\mu_2,\mu_3}=
\e_1^{\mu_1-\mu_2}\e_2^{\mu_2-\mu_3}\e_3^{\mu_3}
\cdot{}_2F_1\Big({\mu_2{-}\mu_1,\mu_3{-}\mu_2\atop
\mu_3{-}\mu_1{-}2r{+}1}\Big|{\e_3\over\e_1\,\e_2}\Big)
$$
Secondly, we can express the sum \cite{E41} as a Jacobi
polynomial which is defined as
$$\eqno{}
P_n^{\alpha,\beta}(x):={\alpha+n\choose n}\cdot{}_2F_1\Big({-n,n+\alpha+\beta+1\atop\alpha+1}\Big|{1-x\over2}\Big).
$$
For this, we invert the order of the
summands. Since $k$ runs from $0$ to the smaller of $\mu_1-\mu_2$ and
$\mu_2-\mu_3$ we have two cases. Set
$\mu=(k_1+k_2+k_3,k_2+k_3,k_3)$. Then the first case is $k_1\le
k_2$. The substitution $k=k_1-l$ gives
$$\eqno{}
(-1)^k{{k_1\choose k}{k_2\choose
k}\over{k_1+k_2+2r-1\choose k}}=(-1)^{k_1}{{k_2\choose
k_1}\over{k_1+k_2+2r-1\choose k_1}}
{(-k_1)_l\,(k_2+2r)_l\over (k_2-k_1+1)_l\,l!}
$$
where $(a)_l=a(a+1)\ldots(a+l-1)$ is the Pochhammer symbol. Thus:
$$\eqno{}
\Rq_\mu={k_2\choose k_1}{-k_2-2r\choose k_1}^{-1}\e_2^{k_2-k_1}\e_3^{k_1+k_3}\cdot{}_2F_1
\Big({-k_1,k_2+2r\atop k_2-k_1+1}\Big|{\e_1\e_2\over\e_3}\Big).
$$
This, and a similar computation for $k_1\ge k_2$ gives
$$\eqno{}
\Rq_\mu(z)=\cases{{-k_2-2r\choose k_1}^{-1}\cdot\displaystyle
\e_2^{k_2-k_1}\e_3^{k_1+k_3}\cdot
P_{k_1}^{k_2-k_1,2r-1}(1-2{\e_1\e_2\over
\e_3})
&for $k_1\le k_2$\cr
{-k_1-2r\choose k_2}^{-1}\cdot\displaystyle
\e_1^{k_1-k_2}\e_3^{k_2+k_3}\cdot
P_{k_2}^{k_1-k_2,2r-1}(1-2{\e_1\e_2\over\e_3})
&\vrule height 20pt width 0pt for $k_1\ge k_2$\cr}
$$
These formulas are essentially due to Vilenkin--\v Sapiro~\cite{VS},
see also \cite{VK}~11.3.2.
\medskip
\noindent 2. {\it The case $n=4$.} In this case, we put $\mu_2-\mu_3=k$ and
$\mu_4=l$. Then $\mu=(a-k-l,b-l,k+l,l)$ and we get
$$\eqno{}
\eqalign{\Rq_{a,b,0,0}&={1\over c_{ab}}\sum_{k,l}{a{-}k{-}l\choose
a{-}b{-}k,b{-}k{-}2l,k,l}{-2r\choose a{-}k{-}l}\e_\mu\cr
&=\sum_{k,l}{(-a{+}b)_k\,(-b)_{k{+}2l}\over(-a{-}2r{+}1)_{k{+}l}\,k!\,l!}
\e_1^{a{-}b{-}k}\e_2^{b{-}l{-}2l}\e_3^k\e_4^l\cr}
$$
where
$$
\e_1=z_1-z_2+z_3-z_4,\quad\e_2=z_2+z_4,\quad\e_3=z_1z_3-z_2z_4,\quad\e_4=z_2z_4.
$$
This can be expressed in terms of one of Horn's hypergeometric
functions (see e.g. \cite{Bate} \S5.7.1):
$$\eqno{}
\Rq_{a,b,0,0}=\e_1^{a{-}b}\e_2^b\cdot
H_3(-b,-a{+}b,-a{-}2r{+}1;{\e_4\over\e_2^2},{\e_3\over\e_1\e_2})
$$

\beginsection Triangularity. Triangularity

In this section, we investigate vanishing properties of the
coefficients of $R_\lambda(z;r)$. For this, we consider the {\it
inhomogeneous dominance order\/}: for $\mu,\lambda\in\NN^n$ define
$$\eqno{}
\mu\le\lambda
\quad{\buildrel{\rm def}\over\Longleftrightarrow}\quad
\mu_1+\ldots+\mu_m\le\lambda_1+\ldots+\lambda_m\quad\hbox{for all
}m=1,\ldots,n.
$$
The homogeneous dominance order, commonly considered in the theory of
symmetric functions, is
$$\eqno{}
\mu\lee\lambda
\quad{\buildrel{\rm def}\over\Longleftrightarrow}\quad
\mu\le\lambda\hbox{\ and\ }|\mu|=|\lambda|.
$$
Recall, that we defined the leading term of $R_\lambda$ as
$z^{[\lambda]}$ where $[\lambda]$ is defined in \cite{E19}. The next
theorem justifies this terminology.

\Theorem Triangular1. For every
$\lambda\in\Lambda$ there are expansions
$$\eqno{}
R_\lambda(z)=\sum_{\mu\in\NN^n:\mu\le[\lambda]} a_{\lambda\mu}z^\mu
\quad\hbox{and}\quad
\Rq_\lambda(z)=\sum_{\mu\in\NN^n:\mu\lee[\lambda]}
a_{\lambda\mu}z^\mu.
$$

\Proof: For $1\le m\le n$ and $f\in\cP$ denote the total degree of $f$ in
$z_1,\ldots,z_m$ by $\|deg|_mf$. Let $\muq:=\lfloor m/2\rfloor$ and
$\moq:=m-\muq=\lceil m/2\rceil$. We show first $\|deg|_m
X(t)f\le\|deg|_mf$ and $\|deg|_m Y(t)f\le\|deg|_mf$ for all
$f\in\cP^W$.

\cite{degree} is nothing else than the case $m=n$. The general case is
the same except that the entries in the rows involving
$x_{\moq+1},\ldots, x_\noq$ and $y_{\muq+1},\ldots, y_\nuq$ have
degree $0$. Thus the degree of $\fX(t)$ can be computed by taking in
\cite{E15} or \cite{E16} the $m$ largest entries of the $d_j''$ and
the entries of $d_i'$ which correspond to
$x_1,\ldots,x_\moq,y_1,\ldots,y_\muq$. Thus
$$
\eqalignno{
&\deg_m\det\fX(t)&
=\sum_{i=1}^\moq(\noq-i)+\sum_{i=1}^\muq(\moq-i)+\muq(\nuq-\noq)=
\sum_{i=1}^\moq(\noq-i)+\sum_{i=1}^\muq(\muq-i),\cr
&\deg_m\det\fY(t)&\vtop{\hbox{$\displaystyle
=\sum_{i=1}^\moq(\noq-i)+\sum_{i=1}^\muq(\moq-i-1)+\muq(\nuq-\noq+1)=$}
\hbox{$\displaystyle
=\sum_{i=1}^\moq(\noq-i)+\sum_{i=1}^\muq(\muq-i).$}}\cr}
$$
On the other hand
$\|deg|_m\phi(z)=\sum_{i=1}^\moq(\noq-i)+\sum_{i=1}^\muq(\muq-i)$
which proves the claim.

For $\lambda\in\Lambda$ let $\cP_\lambda$ and $\cP_\lambda^\circ$ be
the intersection of $\cP^W$ with the span of all $z^\mu$ with
$\mu\le[\lambda]$ and $\mu<[\lambda]$, respectively. Then, by what we
have proved above, both $\cP_\lambda$ and $\cP_\lambda^\circ$ are
stable under $X(t)$ and $Y(t)$. The monomial symmetric polynomial
$m_{[\lambda]}(z)$ is in $\cP_\lambda$ but not in $\cP_\lambda^\circ$.
Thus, $\cP_\lambda^\circ$ is of codimension one in
$\cP_\lambda$. Because the action of $X(t),Y(t)$ is diagonalizable
there is exactly one $\nu_\lambda\in\Lambda$ such
that $R_{\nu_\lambda}$ is in $\cP_\lambda$ but not in
$\cP_\lambda^\circ$.

It remains to show $\nu_\lambda=\lambda$ for all $\lambda$. If there
exists a counterexample then choose one which is minimal with respect to
the order relation $[\nu]\le[\lambda]$. Since
$R_{\nu_\lambda}$ contains $z^{[\nu_\lambda]}$
(\cite{LeadingTerm}) we have $[\nu_\lambda]<[\lambda]$. Thus, by
minimality,
$R_{\nu_\lambda}\in\cP_{\nu_\lambda}\subseteq
\cP_\lambda^\circ$ in contradiction to the definition of $\nu_\lambda$.\qed

\Examples: 1. If $\lambda$ is of the form $(a,a,b,b,c,c,\ldots)$ then
we know from \cite{Jack1} that $R_\lambda$ is a polynomial in the even
variables $z_2,z_4,\ldots$ only. This can also be seen from
triangularity: since $[\lambda]_1=0$ we have $\mu_1=0$ for every
$z^\mu$ which occurs in $R_\lambda$. Hence $z_1$ does not
occur. By symmetry, no odd variable occurs.

2. If $\lambda$ is of the form $(a,b,b,c,c,\ldots)$ then
$[\lambda]=(a,0,b,0,\ldots)$. Hence triangularity prohibits, e.g., the
occurrence of monomials $z_1^{\mu_1} z_2^{\mu_2}\ldots$ with $\mu_1+\mu_2>a$.
\medskip

This form of triangularity seems to be optimal when the expansion of
$R_\lambda$ in monomials is considered but, since monomials are not
bihomogeneous, it does not cover the bigrading result of
\cite{Bigrading}. Therefore, we expand $R_\lambda$ in elementary
semisymmetric symmetric functions $\e_\mu$ defined in
\cite{E60}. Then, an equivalent form of \cite{Triangular1} is that for
every $\lambda\in\Lambda$ there are expansions
$$\eqno{}
R_\lambda(z)=\sum_{\mu\in\NN^n:[\mu]\le[\lambda]} a_{\lambda\mu}\e_\mu
\quad\hbox{and}\quad
\Rq_\lambda(z)=\sum_{\mu\in\NN^n:[\mu]\lee[\lambda]}
a_{\lambda\mu}\e_\mu.
$$
The point is now to define an order relation on $\Lambda$ which is
stronger than $[\mu]\le[\lambda]$.

For this, let $\Phi^+\subseteq\ZZ^n$ be the submonoid generated by all
simple roots $\epsilon_i-\epsilon_{i+2}$, $1\le i\le n-2$. Recall from
\cite{Coeffs}, that $\Psi_0\subseteq\ZZ^n$ was defined to be the
smallest $W$\_stable monoid containing $\Lambda$. We define the
semisymmetric analogue of the inhomogeneous dominance order on
$\Lambda$ as
$$\eqno{}
\mu\preceq\lambda
\quad{\buildrel{\rm def}\over\Longleftrightarrow}\quad
\lambda-\mu\in\Psi_1:=\Psi_0+\Phi^+
$$
\Lemma. The monoid $\Psi_1$ has also the following descriptions:
\item{a)}It is generated by $\{\epsilon_i-\epsilon_{i+2}\mid 1\le i\le
n-2\}\cup\{\epsilon_{n-1}+\epsilon_n,\epsilon_{2\noq-1}\}$.
\item{b)}It consists of all $\lambda\in\ZZ^n$ with $0\le\lambda_\odd$,
$0\le\lambda_\even$, and $|\lambda_\even|\le|\lambda_\odd|$.

\Proof: a) By \cite{Psi0}, the monoid $\Psi_0$ is generated by all
elements of the form $\epsilon_i+\epsilon_j$, ($i$ odd, $j$ even) and
$\epsilon_i$, ($i$ odd). Using the generators of $\Phi^+$, one can
obtain all these generators from either $\epsilon_{n-1}+\epsilon_n$ or
$\epsilon_{2\noq-1}$ alone. This shows the claim.

b) First observe that the set of generators $\Sigma$ in a) forms in
fact a linear basis of $\ZZ^n$. Now consider the set $\Sigma'$
consisting of the linear forms
$\lambda_1,\lambda_2,\lambda_1+\lambda_3,\lambda_2+\lambda_4,
\lambda_1+\lambda_3+\lambda_5,\ldots$ and
$\lambda_1-\lambda_2+\lambda_3-+\ldots$. Then the conditions in b) can
be rephrased as $\ell(\lambda)\ge0$ for all $\ell\in\Sigma'$. Observe
that $\Sigma'$ contains $|\lambda_\odd|$ which is a sum of two other
elements, thus redundant. When we remove it from $\Sigma'$ we obtain a
set $\Sigma^*$ which turns out to be the dual basis of $\Sigma$. Thus
$\Psi_1$ equals the set $\lambda\in\ZZ^n$ with $\ell(\lambda)\ge0$ for
all $\ell\in\Sigma^*$.\qed

\noindent
Since $[\lambda]_1=|\lambda_\odd|-|\lambda_\even|$, we have in particular
$$\eqno{}
\mu\preceq\lambda
\quad\Longleftrightarrow\quad
\mu_\odd\le\lambda_\odd,\ \mu_\even\le\lambda_\even\hbox{, and }
[\mu]_1\le[\lambda]_1.
$$
Next, we compare $\mu\preceq\lambda$ with $[\mu]\le[\lambda]$.

\Lemma. The monoid $\tilde\Psi_1:=\{\lambda\in\ZZ^n\mid0\le[\lambda]\}$
is generated by
$$\eqno{}
\{\epsilon_i-\epsilon_{i+2}\mid 1\le i\le
n-2\}\cup\{\epsilon_{n-1}+\epsilon_n,-\epsilon_2\}.
$$

\Proof: One easily checks that the proposed set of generators is the dual
basis to
$$\eqno{}
\{[\lambda]_1+\ldots+[\lambda]_m\mid 1\le m\le n\}.
$$\qed

\noindent The new order relation is indeed stronger than the one
considered before:

\Corollary. $\mu\preceq\lambda$ implies $[\mu]\le[\lambda]$.

\Proof: We have
$-\epsilon_{2\nuq}=-\epsilon_2+(\epsilon_2-\epsilon_4)+
\ldots+(\epsilon_{2\nuq-2}-\epsilon_{2\nuq})$. Since $2\nuq$ (resp.
$2\noq-1$) is the largest even (resp. odd) integer in $1,\ldots,n$ we
have $\epsilon_{2\noq-1}=
-\epsilon_{2\nuq}+(\epsilon_{n-1}+\epsilon_n)$. This implies
$\Psi_1\subseteq\tilde\Psi_1$ which is equivalent to the
assertion.\qed

The homogeneous version of ``$\preceq$'' is defined as
$$\eqno{}
\mu\LEE\lambda
\quad{\buildrel{\rm def}\over\Longleftrightarrow}\quad
\mu\preceq\lambda\hbox{ and }|\mu|_\odd=|\lambda|_\odd
$$
Since $|\lambda|_\odd=0$ for all $\lambda\in\Phi^+$ and
$|\lambda|_\odd>0$ for all $\lambda\in\Psi_0\setminus\{0\}$ the
definition of $\mu\LEE\lambda$ simplifies to
$\lambda-\mu\in\Phi^+$. Thus, we get
$$\eqno{}
\mu\LEE\lambda
\quad\Longleftrightarrow\quad
\mu_\odd\lee\lambda_\odd\hbox{ and }\mu_\even\lee\lambda_\even.
$$

Now we are looking at expansions of elements in $\cP^W$ in the form
$\sum_\mu a_{\lambda\mu}\e_\mu$. For technical reasons we need a
version which works for all elements in $\cP$.

\Lemma Expansions. For $f\in\cP^W$ and $\lambda\in\Lambda$ the following
statements are equivalent:
\item{a)}In the expansion
$$\eqno{}
f(u_1+u_2,u_2,u_3+u_4,u_4,\ldots)=\sum_{\mu\in\NN^n}a_\mu u^\mu
$$
(where $u_{n+1}:=0$ if $n$ is odd) only monomials $u^\mu$ with
$\mu\le[\lambda]$ and $\mu_\even\le\lambda_\even$ occur.
\item{b)}There is an expansion
$$\eqno{}
f(z)=\sum_{\mu\in\Lambda:\mu\preceq\lambda}b_\mu\e_\mu.
$$
\Par

\Proof: ``{\it b)$\Rightarrow$a)}'': Let $\deg_mf$ be the total degree of
$f$ in $u_1,\ldots,u_m$ which is the same as the degree in
$z_1,\ldots,z_m$. Then one calculates
$$\eqno{}
\deg_m\e_\lambda=\sum_{i=1}^m[\lambda]_i=
\cases{\sum_{i=1}^{m/2}\lambda_{2i-1}&if $m$ is even,\cr
\sum_{i=1}^{(m-1)/2}\lambda_{2i}+[\lambda]_1&if $m$ is odd.\cr}
$$
We conclude that $\mu\preceq\lambda$ implies
$\deg_m\e_\mu\le\deg_m\e_\lambda$. Thus, if $u^\mu$ occurs in $f$ then
$$\eqno{}
\mu_1+\ldots+\mu_m=\deg_mu^\mu\le\deg_mf\le\deg_m\e_\lambda=
[\lambda]_1+\ldots+[\lambda]_m,
$$
i.e., $\mu\le[\lambda]$.

Now let $\deg_m^u$ be the total degree of $f$ in
$u_2,u_4,\ldots,u_{2m}$. Then, due to cancellations, one has
$\deg_m^u\e_\lambda=\sum_{i=1}^m\lambda_{2i}$. The same reasoning as
above implies $\mu_\even\le\lambda_\even$ whenever $u^\mu$ occurs in $f(u)$.

``{\it a)$\Rightarrow$b)}'': Assume $\e_\mu$ occurs in the expansion
of $f$. Then, by the calculations above, we have to show
$\deg_m\e_\mu\le\deg_mf$ and $\deg_m^u\e_\mu\le\deg_m^uf$ for all $m$
(actually is suffices to consider in the first case only $m=1$ and all
even $m$).

We treat $\deg_m$ first. For this define a total order on the
monomials $u^\lambda$: first we order them by $\deg_m$ and then by the
lexicographic order on
$(\lambda_1,\lambda_3,\ldots,\lambda_2,\lambda_4,\ldots)$. One checks
that $\e_p$ has the leading monomial $u_1u_3\ldots u_p$ if $p$ is
odd and $u_2u_4\ldots u_p$ if $p$ is even. Thus, the leading monomial
of $\e_\lambda$ is $u^{[\lambda]}$. This shows in particular that if
the leading monomials of $\e_\lambda$ and $\e_\mu$ coincide then
$\lambda=\mu$. Therefore, if there were an $\e_\mu$ occurring in $f$
with $\deg_m\e_\mu>\deg_mf$ then take a maximal one. Its leading
monomial $u^\nu$ would not cancel out and would satisfy
$\deg_mu^\nu>\deg_mf$ in contradiction to {\it a)}.

For $\deg_m^u$ we argue similarly. This time the total order on the
monomials $u^\lambda$ is by lexicographic order on
$(\deg_m^uu^\lambda,\lambda_2,\lambda_4,\ldots,\lambda_1,\lambda_3,\ldots)$.
Then the leading term of $\e_p$ is $u_2u_4\ldots u_p$ if $p$ is even
and $u_2u_4\ldots u_{p-1}u_p$ if $p$ is odd. Hence the leading term of
$\e_\lambda$ is
$u_1^{\lambda_1-\lambda_2}u_2^{\lambda_2}u_3^{\lambda_3-\lambda_4}\ldots
u_n^{\lambda_n}$. These terms are again distinct for different
$\e_\lambda$'s. The rest of the argument is as above.\qed

Now we can state the better triangularity result announced earlier.

\Theorem Triangular2. For every $\lambda\in\Lambda$ there are expansions
\item{a)}$\displaystyle R_\lambda(u_1+u_2,u_2,u_3+u_4,u_4,\ldots)=
\sum_{\mu\in\NN^n}a_{\lambda\mu}u^\mu
\hbox{ where }\mu\le[\lambda],\ \mu_\even\le\lambda_\even$;
\item{b)}
$\displaystyle R_\lambda(z)=\sum_{\mu\in\Lambda:\mu\preceq\lambda}
b_{\lambda\mu}\e_\mu(z)\quad
\hbox{and}\quad
\Rq_\lambda(z)=\sum_{\mu\in\Lambda:\mu\LEE\lambda}b_{\lambda\mu}\e_\mu(z)$.

\Proof: By \cite{Expansions}, it suffices to prove {\it a)}. The degree
function $\|deg|_m$ (see last proof) is invariant under upper
triangular linear coordinate transformations.  Thus $\mu\le[\lambda]$
follows from \cite{Triangular1}.

To prove $\mu_\even\le\lambda_\even$ recall that the total degree of
$f(u_1+u_2,u_2,u_3+u_4,u_4,\ldots)$ in the coordinates $u_2,u_4,\ldots
u_{2m}$ is denoted by $\|deg|_m^uf$. We show first that the operators
$X(t)$ and $Y(t)$ preserve $\|deg|_m^u$.

The substitution $z_{2i-1}\pfeil u_{2i-1}+u_{2i}$ corresponds to
$x_i\pfeil x_i+y_i$ and $T_{y,i}\pfeil T_{y,i}T_{x,i}^{-1}$ in
$\fX(t)$ and $\fY(t)$. The same reasoning as in the proof of
\cite{Triangular1} shows that $\deg_m^u\det\fX(t)$ and
$\deg\det\fY(t)$ are bounded by $\|deg|\phi+m$. Thus we have to find a
way to decrease this estimate for the degree of $\det\fX(t)$ and
$\det\fY(t)$ by $m$.

The idea is to add a multiple of the $y_i$\_row to the
$x_i$\_row. Since the entries are in a non\_commutative ring some care
is advised. For this, we develop $\det\fX(t)$ and $\det\fY(t)$ as
$$\eqno{}
\sum\pm\det A_{i_1,j_1}^{1,\noq{+}1}\ldots
\det A_{i_m,j_m}^{m,\noq{+}m}
\det A_S^{m+1,\ldots,\noq,\noq{+}1{+}m,\ldots,n}
$$
where $\det A^U_V$ is the minor with row index in $U$ and column index
in $V$ and where the sum runs through all partitions
$\{1,\ldots,n\}=\{i_1,j_1\}\dot\cup\ldots\dot\cup\{i_m,j_m\}\dot\cup
S$. The degree of the last factor (involving $S$) is zero. Now we show
that the degree of each $2\times 2$\_minor is one less than expected
which would prove the claim.

For this we write
$$\eqno{}
\det A_{i_l,j_l}^{l,\noq+l}=
\det\pmatrix{x_{11}&x_{12}\cr x_{21}&x_{22}\cr}=(x_{11}+\alpha
x_{21})x_{22}-x_{21}(x_{12}+\alpha x_{22})-[\alpha,x_{21}]x_{22}
$$
where
$$\eqno{}
\alpha:=\cases{T_{x,l}&if $A=\fX(t)$ and $n$ is even\cr
            y_lT_{x,l}&if $A=\fX(t)$ and $n$ is odd\cr
       y_l^{-1}T_{x,l}&if $A=\fY(t)$ and $n$ is even\cr
               T_{x,l}&if $A=\fY(t)$ and $n$ is odd\cr}
$$
As mentioned above this amounts to add $\alpha$ times row \#$\noq+l$
to row \#$l$ of \cite{E11} or \cite{E12}, respectively.
Then it is easy to check that $\deg_m^u(x_{1j}+\alpha x_{2j})\le
\deg_m^u x_{1j}-1$ and
$\deg_m^u[\alpha,x_{21}]\le\deg_m^ux_{11}-1$ which proves the claim.

The rest of the proof is the same as for \cite{Triangular1}. For
$\lambda\in\Lambda$ let $\cP_\lambda$ be the space of all semisymmetric
functions $f$ in which only monomials $u^\mu$ with $\mu\le[\lambda]$
and $\mu_\even\le\lambda_\even$ occur. Let $\cP_\lambda^\circ$ the
same with additionally $\mu\ne\lambda$. Then both spaces are stable
under $X(t)$ and $Y(t)$. Moreover
$\e_\lambda\in\cP_\lambda\setminus\cP_\lambda^\circ$. We conclude as in
\cite{Triangular1}.\qed

\Examples: The improvement of strong triangularity over the weak one
is the more significant the smaller $\lambda_\even$ is. The most
extreme case is $\lambda=(a)$ where \cite{Triangular1} doesn't give
any restriction. But \cite{Triangular2} states $R_{(a)}=\sum_{i=0}^a
c_i \e_{(i)}$ which is of course also a consequence of the direct
calculation in \cite{Hook1}. A more specific example is
$\lambda=(5,2,0,\dots)$ for $n\ge10$. In that case, $R_\lambda$ has,
according to \cite{Triangular1}, $70$ independent coefficients while
\cite{Triangular2} boils that down to $27$.

\Remark: Part b) of the theorem is entirely analogous to a similar
theorem for (shifted) Jack polynomials but part a) is a bit strange
since the pretty asymmetric coordinates $u_i$ appear. A conceptual
explanation for their appearance would be very desirable.
\medskip \noindent
Now we can prove a triangularity property which is completely
intrinsic for the polynomials $R_\lambda$:

\Theorem. For every $\lambda,\mu\in\Lambda$ consider the expansion
$R_\lambda R_\mu=\sum_\tau a_{\lambda\mu}^\tau R_\tau$. Then
$a_{\lambda\mu}^\tau=0$ unless $\lambda,\mu\LE\tau\preceq\lambda+\mu$.

\Proof: First we show $\tau\preceq\lambda+\mu$ whenever
$a_{\lambda\mu}^\tau\ne0$. We have
$$\eqno{}
R_\lambda R_\mu=
\sum_{\tau_1,\tau_2}b_{\lambda\tau_1}b_{\mu\tau_2}\e_{\tau_1}\e_{\tau_2}=
\sum_\nu c_\nu\e_\nu\quad\hbox{with}\quad c_\nu=\sum_{\tau_1+\tau_2=\nu}
b_{\lambda\tau_1}b_{\mu\tau_2}.
$$
Moreover, $\tau_1\preceq\lambda$ and $\tau_2\preceq\mu$ imply
$\nu=\tau_1+\tau_2\preceq\lambda+\mu$. Now, observe that the
transformation matrix $b_{\lambda\mu}$ is upper unitriangular. Thus,
its inverse matrix has the same property, i.e., we have expansions
$\e_\lambda=\sum_{\mu\preceq\lambda}b_{\lambda\mu}'R_\mu$. Hence
$$\eqno{}
R_\lambda R_\mu=\sum_{\nu,\tau}c_\nu b_{\nu\tau}'R_\tau=\sum_\tau
a_{\lambda\mu}^\tau R_\tau
$$
with $\tau\preceq\nu\preceq\lambda+\mu$.

Now we show $\lambda\LE\tau$ whenever $a_{\lambda\mu}^\tau\ne0$. The
relation $\mu\LE\tau$ follows then by symmetry. Let $\tau_0$ be a
$\LE$\_minimal counterexample. Since $\lambda\not\LE\tau_0$, the Extra
Vanishing Theorem~\ncite{ExtraVanishing} implies
$R_\lambda(\rho+\tau_0)=0$. Hence $\sum_\tau a_{\lambda\tau}^\tau
R_\tau(\rho+\tau_0)=0$. Again by the Extra Vanishing Theorem, only
those $\tau$ with $\tau\LE\tau_0$ contribute to this sum. For those we
have $\lambda\not\LE\tau$. The minimality of $\tau_0$ implies
$\tau=\tau_0$ unless $a_{\lambda\mu}^\tau=0$. From this we derive the
contradiction $a_{\lambda\mu}^{\tau_0}R_{\tau_0}(\rho+\tau_0)=0$.\qed

\beginsection Binomial. The binomial theorem

In this section we derive a binomial type theorem for semisymmetric
functions. The proof is similar to that for symmetric functions in
\cite{Ok}.

So far, we considered values of $R_\lambda$ in the points
$z=\rho+\lambda$. Now we use that the difference operators also have a
dual vanishing property.

Recall that $\cR$ is the algebra generated by the $X_i$ and $Y_i$
where $X(t)=\sum_iX_it^i$ and $Y(t)=\sum_iY_it^i$. We introduce a
degree function on $\cR$ by letting $\|deg|z_i=0$ and
$\|deg|T_\lambda:=|\lambda|_\odd$. Thus, $\|deg|X_i=\|deg|Y_i=i$.
Observe that $(\cP^W)^*$ is a right $\|End|\cP^W$\_module, hence a
right $\cR$\_module. For $\tau\in k^n$ let $\delta_\tau\in(\cP^W)^*$
be the evaluation map $f\mapsto f(\tau)$. Then, as explained in its
proof, \cite{Coeffs}{\it d)} amounts to
$\oplus_{\lambda\in\Lambda}k\delta_{\rho+\lambda}$  being an
$\cR$\_submodule of $(\cP^W)^*$. Now, for any $\a\in k$, let
$\underline\a:=(\a,\ldots,\a)\in k^n$ and $\rho_\a:=
\rho+\underline\a=((n-i)\r+\a)_i$. Then we have

\Proposition Inter. Assume $r\ne0$. Consider the space
$M:=\oplus_{\lambda\in\Lambda}k\,\delta_{-\rho_\a-\lambda}$.
\item{a)} $M$ is an $\cR$\_submodule of $(\cP^W)^*$.  \item{b)} Define
a filtration on $M$ by putting
$\|deg|\delta_{-\rho_\a-\lambda}:=|\lambda|_\odd$. Assume $\a\not\in
-\NN-\NN\rr$. Then the map $\cR\pfeil M:D\mapsto \delta_{-\rho_\a}D$
is an isomorphism of filtered $k$-vector spaces.\Par

\Proof: Let $\lambda\in\Lambda$ and $D$ either $X_i$ or $Y_i$. Then,
the assumption $r\ne0$ makes sure that $\delta_{-\rho_\a-\lambda}D$
can be computed in the obvious way since then the denominator of $D$
does not vanish at $-\rho_\a-\lambda$.

{\it a)} It suffices to show that $\delta_{-\rho_\a-\lambda}D_I\in M$
for any $I\subseteq\{1,\ldots,n\}$. This is no problem if
$\lambda+\epsilon_I\in\Lambda$. Otherwise, there is an index $j$ with
$j\not\in I$, $i:=j+1\in I$, and $\lambda_i=\lambda_j$. But then the
factor $z_i-z_j-\r$ in $D_I$ vanishes at $z=-\rho_\a-\lambda$.

{\it b)} The map clearly preserves filtrations. Since corresponding
filtration spaces on both sides have the same dimension, it suffices
to show surjectivity. We do that by induction on the degree. By the
explicit formulas \cite{E44}, \cite{E45} we have
$$
\eqalignno{
&X_m&=(-1)^m\sum_{I\in P_\odd\atop|I|_o=m}D_I
+\sum_{I:|I|_o<m}c_I^{(1)}(z)D_I\cr
&Y_m&=(-1)^m\sum_{I\in
P_\even\atop|I|_o=m}D_I +\sum_{I:|I|_o<m}c_I^{(2)}(z)D_I.\cr}
$$
Thus, if we put $Z_{2m-1}:=(-1)^m(X_m-Y_m)$ for
$m=1,\ldots,\noq$ and $Z_{2m}:=(-1)^m Y_m$ for
$m=0,\ldots,\nuq-1$ we obtain operators with an expansion
$$\eqno{}
Z_p=\sum_{\mu\in\NN^n:\mu\preceq (1^p)}c_\mu(z)T_\mu
$$
where
$$\eqno{}
c_{(1^p)}(z)=\prod_{i\le p\atop n-i\ \even}\!\!\!\!z_i
\prod_{i\le p<j\atop j-i\ \odd}\!\!(z_i-z_j-\r)
\prod_{i\le p<j\atop j-i\ \even}\!\!(z_i-z_j)^{-1}.
$$
Let $\lambda\in\Lambda$ be non\_zero. Let $p$ be maximal with
$\lambda_p\ne0$. Put $\mu:=\lambda-(1^p)\in\Lambda$. Then we have
$$\eqno{}
\delta_{-\rho_\a-\mu}Z_r\in
c_{(1^p)}(-\rho_\a-\mu)\delta_{-\rho_\a-\lambda}+
\sum_{\nu\prec\lambda}k\delta_{-\rho_\a-\nu}.
$$
The assumptions on $r$ and $\alpha$ ensure that
$c_{(1^p)}(-\rho_\a-\mu)\ne0$. The induction hypothesis implies that
$\delta_{-\rho_\a-\lambda}$ is in the image.\qed

\Lemma. Assume $r\ne0$ and $\a\not\in-\NN-\NN\rr$. Then
$R_\lambda(-\rho_\a)\ne0$ for all $\lambda\in\Lambda$.

\Proof: Suppose $R_\lambda(-\rho_\a)=0$. For
$\mu\in\Lambda$ let $D\in\cR$ with
$\delta_{-\rho_\a-\mu}=\delta_{-\rho_\a}D$. Since $R_\lambda$ is an
eigenvector of $D$ we also have $R_\lambda(-\rho_\a-\mu)=0$. This
contradicts the fact that $-\rho_\a-\Lambda$ is Zariski dense in
$k^n$.\qed

\Remark: This lemma is only preliminary. Later, we prove the explicit
formula \cite{E20} for $R_\lambda(-\rho_\a)$.
\medskip
The binomial type theorem, announced in the beginning, is:

\Theorem Binomialformula. Assume $r\ne0$ and $\a\not\in-\NN-\NN\rr$.
Then for every $\lambda\in\Lambda$ the following formula holds:
$$\eqno{E1}
{R_\lambda(-\underline{\a}-z)\over R_\lambda(-\rho_\a)}=
\sum_{\mu\in\Lambda}(-1)^{|\mu|_\odd}
{R_\mu(\rho+\lambda)\over R_\mu(\rho+\mu)}
{R_\mu(z)\over R_\mu(-\rho_\a)}.
$$

\Proof: The polynomials $R_\lambda(-\underline{\a}-z)$ form also a
basis of $\cP^W$. Hence, every $f\in\cP^W$ has an expansion
$$\eqno{E80}
f(z)=\sum_\mu a_\mu(f)R_\mu(-\underline{\a}-z)
$$
with $a_\mu\in(\cP^W)^*$. We claim $a_\mu\in M$ with
$\|deg|a_\mu\le|\mu|_\odd$. To show this, we evaluate \cite{E80} at
$z=-\rho_\a-\mu$ and get
$$\eqno{}
\delta_{-\rho_\a-\mu}(f)=\sum_\tau a_\tau(f)R_\tau(\rho+\mu)=
a_\mu(f)+\sum_{|\tau|_\odd<|\mu|_\odd}a_\tau(f)R_\tau(\rho+\mu).
$$
Then the claim follows by induction on $|\mu|_\odd$.

It follows from \cite{Inter} that there is $D_\mu\in\cR$ with
$\|deg|D_\mu\le|\mu|_\odd$ such that $a_\mu(f)=(D_\mu
f)(-\rho_\a)$. We apply this to $f=R_\lambda$. Then
$$\eqno{E81}
a_\mu(R_\lambda)=(D_\mu
R_\lambda)(-\rho_\a)=p_\mu(\rho+\lambda)R_\lambda(-\rho_\a)
$$
where $p_\mu:=c_0^{D_\mu}\in\cP^W$ by \cite{Coeffs}{\it c)}. We have
$\|deg|p_\mu=\|deg|D_\mu\le|\mu|_\odd$. On the other side, we see
directly from \cite{E80} that
$$\eqno{E82}
a_\mu(R_\lambda)=\cases{0&if $|\lambda|_\odd\le|\mu|_\odd$ and
$\lambda\ne\mu$;\cr
(-1)^{|\mu|_\odd}&if $\lambda=\mu$.\cr}
$$
Thus \cite{E81}, \cite{E82} together and the very definition of
$R_\mu(z)$ imply
$$\eqno{}
R_\mu(-\rho_\a)p_\mu(z)={(-1)^{|\mu|_\odd}\over
R_\mu(\rho+\mu)}R_\mu(z)
$$
and therefore
$$\eqno{}
a_\mu(R_\lambda)=(-1)^{|\mu|_\odd}{R_\mu(\rho+\lambda)\over
R_\mu(\rho+\mu)}{R_\lambda(-\rho_\a)\over R_\mu(-\rho_\a)}.
$$
Now, we insert this into \cite{E80}, replace $z$ by
$-\underline{\a}-z$ and obtain \cite{E1}.\qed

\Remarks: 1. By \cite{ExtraVanishing}, only those $\mu$ with
$\mu\LE\lambda$ contribute to the sum in formula~\cite{E1}. In
particular, the sum is finite.

\noindent 2. The normalizing factor $R_\lambda(-\rho_\a)$ in the
denominator renders the formula more symmetric but causes the
restriction on $\a$. Of course, for every $\a$ there is an expansion
of $R_\lambda(-\underline{\a}-z)$ in terms of $R_\mu(z)$. It can be
easily obtained by using the explicit formula \cite{E20} to calculate
the ratio $R_\lambda(-\rho_\a)/R_\mu(-\rho_\a)$.
\medskip
There are two immediate applications of the binomial formula \cite{E1}.

\Corollary Symmetric. Assume $r\ne0$ and $\a\not\in-\NN-\NN\rr$. Then
the matrix
$$\eqno{E21}
\left({R_\lambda(-\rho_\a-\nu)\over
R_\lambda(-\rho_\a)}\right)_{\lambda,\nu\in\Lambda}
$$
is symmetric.

\Proof: Substitute $z=\rho+\nu$ in \cite{E1}. Then the right hand side
becomes clearly symmetric in $\lambda$ and $\nu$.\qed

\Corollary Involution. The matrix
$$\eqno{E22}
\left((-1)^{|\mu|_\odd}
{R_\mu(\rho+\lambda)\over
R_\mu(\rho+\mu)}\right)_{\mu,\lambda\in\Lambda}
$$
is an involution.

\Proof: By the binomial formula \cite{E1}, the matrix \cite{E22}
expresses the effect of the involution $z\mapsto-\underline{\a}-z$ on
$\cP^W$ with respect to the basis ${R_\mu(z)\over R_\mu(-\rho_\a)}$,
at least if $r\ne0$. For $r=0$ we argue by continuity.\qed

The involutory matrix \cite{E22} can be used to derive an explicit
interpolation formula (\cite{Interpol}{\it iii)} below). For this, let
$\cC(\rho+\Lambda)$ be the set of $k$\_valued functions on
$\rho+\Lambda$. For $f\in\cC(\rho+\Lambda)$ we define $\hat
f\in\cC(\rho+\Lambda)$ by
$$\eqno{}
\hat f(\rho+\lambda):=\sum_{\mu\in\Lambda}(-1)^{|\mu|_\odd}
{R_\mu(\rho+\lambda)\over R_\mu(\rho+\mu)}f(\rho+\mu)
$$
For any fixed $\lambda$ the sum is finite by the Extra Vanishing
Theorem~\ncite{ExtraVanishing}.  Let
$\cC_0(\rho+\Lambda)\subseteq\cC(\rho+\Lambda)$ be the functions with
finite support. We consider, via restriction, $\cP^W$ as a subspace of
$\cC(\rho+\Lambda)$. By \cite{Cut-off}, $\cC(\rho+\Lambda)$ is a left
$\cR$\_module, provided $r\ne0$.

\Theorem Interpol. The transformation $f\mapsto\hat f$ has the following
properties:
\item{i)}$\hat{\mkern-3mu\hat f\mkern3mu}=f$.
\item{ii)}$f\in\cP^W\Leftrightarrow\hat f\in\cC_0(\rho+\Lambda)$.
\item{iii)}For $f\in\cP^W$ holds
$f(z)=\sum_{\mu\in\Lambda}(-1)^{|\mu|_\odd}\hat f(\rho+\mu)
{R_\mu(z)\over R_\mu(\rho+\mu)}$.
\item{iv)} Assume $r\ne0$. For every $D\in\cR$ holds 
$\widehat{D(f)}=c_0^D\hat f$ and $\widehat{c_0^Df}=D(\hat f)$.

\Proof: {\it i)} follows from \cite{Involution} and the fact that the
transpose of an involutive matrix is involutive.  Let
$\chi_{\rho+\nu}\in\cC_0$ be the characteristic function of the
one-point set $\{\rho+\nu\}$. Then
$\hat\chi_{\rho+\nu}=(-1)^{|\nu|_\odd}R_\nu(\rho+\nu)^{-1}R_\nu$.
Hence, $f\mapsto\hat f$ maps a basis of $\cC_0$ to a basis of $\cP^W$
which proves {\it ii)}. Part {\it iii)} is a direct consequence of
{\it i)} and {\it ii)}.

Finally, let $D\in\cR$. The second formula in {\it iv)} follows from
the first by {\it i)}. Thus we have to prove
$$\eqno{E83}
\widehat{D(f)}(\rho+\lambda)=c_0^D(\rho+\lambda)\hat
f(\rho+\lambda)
$$
for every $D\in\cR$, $\lambda\in\Lambda$, and $f\in\cC$. If we fix $D$
and $\lambda$ then there is a finite subset $S\subset\rho+\Lambda$
such that both sides of \cite{E83} depend only on values of $f$ in
$S$. Since on $S$ every $f\in\cC$ can be interpolated by an element of
$\cP^W$ it suffices to prove \cite{E83} for $f=R_\nu$. But then we
have
$$\eqno{}
\widehat{D(R_\nu)}=[c_0^D(\rho+\nu)R_\nu]^\wedge=c_0^D(\rho+\nu)\hat
R_\nu=c_0^D\hat R_\nu.
$$
The last equality holds since $\hat R_\nu$ is a multiple of the
characteristic function $\chi_{\rho+\nu}$.
\qed

\Corollary. Assume $r\ne0$. Let $\cA\subseteq\|End|_k\cP^W$ be the
algebra generated by $\cP^W$ and $\cR$. Then there is an involutory
automorphism of $\cA$ which interchanges $\cP^W$ and $\cR$.

\Proof: The automorphism is $D\mapsto\hat D$, where $\hat
D(f):=\widehat{D(\hat f)}$. Then \cite{Interpol}{\it i)} implies that
this is an involution and part {\it iv)} implies $\hat D=c_0^D$ for
every $D\in\cR$.\qed

\beginsection EvaluationFormula. The evaluation formula

The symmetry of the matrix \cite{E21} allows to switch
the index with the argument. Using this, we obtain Pieri type formulas:

\Theorem Pieri1. Assume $r\ne0$ and $\a\not\in-\NN-\NN\rr$. Let
$D=\sum_\eta c_\eta^D(z)T_\eta\in\cR$. Then for all $\mu\in\Lambda$
holds
$$\eqno{}
c_0^D(-\underline{\a}-z){R_\mu(z)\over R_\mu(-\rho_\a)}=
\sum_{\lambda\in\Lambda}c_{\lambda-\mu}^D(-\rho_\a-\mu)
{R_\lambda(z)\over R_\lambda(-\rho_\a)}.
$$

\Proof: We substitute $z=-\rho_\a-\mu$ in the equation
$c_0^D(\rho+\nu)R_\nu(z)=D(R_\nu)(z)$ 
and apply symmetry (i.e., \cite{Symmetric}) on both sides. Thus we obtain
$$\eqno{E84}
c_0^D(\rho+\nu)
R_\nu(-\rho_\a)
{R_\mu(-\rho_\a-\nu)\over R_\mu(-\rho_\a)}=
\sum_\eta c_\eta^D(-\rho_\a-\mu)
R_\nu(-\rho_\a)
{R_{\mu+\eta}(-\rho_\a-\nu)\over R_{\mu+\eta}(-\rho_\a)}
$$
After canceling $R_\nu(-\rho_\a)$, both sides of \cite{E84} become
polynomials in $\nu$. Hence we may replace $\nu$ by $-\rho_\a-z$. Then
putting $\eta=\lambda-\mu$ yields the desired formula.\qed

\def\pprod_#1{\prod\limits_{\hbox to 0pt{\hss$\scriptstyle{#1}$\hss}}}

We are applying this to $D=X(t)$ and $D=Y(t)$. By \cite{Expansion}, the
non\_zero coefficients are
$$\eqno{}
\eqalign{
c_{\epsilon_I}^{X(t)}(z)&=(-1)^{|I|_o}\,u_I^\odd(z,t)\,v_I(z)\,w_I(z;r),\quad
I\in P_\odd\cr
c_{\epsilon_I}^{Y(t)}(z)&=(-1)^{|I|_o}\,u_I^\even(z,t)\,v_I(z)\,w_I(z;r),\quad
I\in P_\even\cr}
$$
where
$$
\eqalignno{
&u_I^{\odd/\even}(z,t)&=\pprod_{i\not\in I\atop i\ \odd/\even}(t+z_i),\qquad
v_I(z)=\pprod_{i\in I\atop n-i\ \even}z_i\cr
&w_I(z;r)&=\pprod_{i\in I,j\not\in I\atop j-i\ \odd}(z_i-z_j-r)\cdot
\pprod_{i\in I,j\not\in I\atop j-i\ \even}(z_i-z_j)^{-1}.\cr}
$$
After replacing $t$ by $\a-t$ we obtain
$$\eqno{E7}
\vcenter{\noindent$\eqalign{&
\pprod_{i\ \odd/\even}(t+z_i)\ R_\mu(z)=
\cr&
\qquad=\sum_{I\in P_{\odd/\even}}u_I^{\odd/\even}(\rho+\mu,t)
v_I(\rho_\a+\mu)w_I(\rho_\a+\mu;-r)
{(-1)^{|I|_o}R_\mu(-\rho_\a)\over R_{\mu+\epsilon_I}(-\rho_\a)}
R_{\mu+\epsilon_I}(z)\cr}$}
$$
We postpone the simplification of this formula until
section~\cite{PieriFormula}. Instead, we use \cite{E7} to derive an
explicit formula for $R_\lambda(-\rho_\a)$. To state the
result let
$$\eqno{}
[x\uparrow m]:=x(x+1)\ldots(x+m-1)
$$
be the rising factorial polynomial. For a box $s=(i,j)\in\NN^2$ of a
partition $\lambda$ recall the following notation:
$$\eqno{}
\vcenter{\halign{
$\displaystyle#$\hfill&
$\displaystyle\,:=#$\hfill&
\quad(#)\hfill&\qquad\qquad
$\displaystyle#$\hfill&
$\displaystyle\,:=#$\hfill&
\quad(#)\hfill\cr
a_\lambda(s) &\lambda_i-j &arm-length&
a'(s)&j-1         &arm-colength
\cr
l_\lambda(s) &\lambda_j'-i&leg-length&
l'(s)&i-1         &leg-colength
\cr}}
$$

\Theorem SpecialValue. Assume $r\ne0$. Then for every
$\lambda\in\Lambda$ the evaluation formula
$$\eqno{E20}
R_\lambda(-\rho_\a)=(-1)^{|\lambda|_\odd}A_\lambda(\a)B_\lambda,
$$
holds where
$$\eqno{}
A_\lambda(\a):=\prod_{1\le i\le n\atop n{-}i\ \even}
[\a+(n-i){\textstyle\r}\uparrow\lambda_i]=
\prod_{s\in\lambda\atop n{-}l'(s)\ \odd}
\big(\a+a'(s)+(n-l'(s)-1){\textstyle\r}\big)
$$
and
$$\eqno{E85}
B_\lambda:=
{\prod\limits_{1\le i<j\le n\atop j{-}i\ \odd}[(j-i+1)\r\uparrow\lambda_i-\lambda_j]
\over
\prod\limits_{1\le i<j\le n\atop j{-}i\ \even}[(j-i)\r\uparrow\lambda_i-\lambda_j]
}=
{\prod\limits_{s\in\lambda\atop n{-}l'(s)\ \even}\kern-10pt\big(a'(s)+(n-l'(s))\r\big)
\over
\prod\limits_{s\in\lambda\atop l_\lambda(s)\ \odd}
\big(a_\lambda(s)+(l_\lambda(s)+1)\r\big)}
$$

\Proof: By continuity, we may assume $\a\not\in-\NN-\NN\rr$.
We expand both sides of \cite{E7} as a polynomial in $t$ and compare
coefficients. Then the product on the left-hand side becomes
$e_d(z_{\odd/\even})$ while on the right-hand side the $u_I$\_factor
has to be replaced by $e_{d-|I|_o}(z_i|i\not\in I\
\odd/\even)$. In particular, only sets $I$ with $|I|_o\le d$ enter the
formula. It is easily checked that the set of $\epsilon_I$ with $I\in
P_\odd$ (resp. $I\in P_\even$) with $|I|_o\le d$ has a unique
maximum with respect to the order relation $[\mu]\le[\lambda]$, namely
$(1^b)=\sum_{i=1}^b\epsilon_i$ where $b=2d-1$ (resp. $b=2d$). Thus, by
\cite{Triangular1}, the monomial $z^{[\mu+(1^b)]}$ appears on the right-hand
side exactly once. Comparing its coefficient, we obtain for all
$b=1,\ldots,n$:
$$\eqno{}
R_{\mu+(1^b)}(-\rho_\a)=
(-1)^{\lceil{b\over2}\rceil}v_{(1^b)}(\rho_\a+\mu)w_{(1^b)}(\rho_\a+\mu;-r)
R_\mu(-\rho_\a)
$$
This is a recursion relation which allows to compute
$R_\mu(-\rho_\a)$ by deleting one column at a time. It follows that
$R_\mu(-\rho_\a)=(-1)^{|\mu|_\odd} A_\mu B_\mu$ with
$A_{\mu+(1^b)}=v_{(1^b)}(\rho+\underline{\a}+\mu)A_\mu$ and
$B_{\mu+(1^b)}=w_{(1^b)}(\rho_\a+\mu;-r)B_\mu$.

The first relation implies easily both formulas for $A_\mu$. The
expressions for $B_\mu$ could be derived in the same way as the
analogous formulas \cite{Mac}~VI(6.11) and ($6.11'$) for Macdonald
polynomials. Especially the second expression for $B_\lambda$ in
\cite{E85} is quite tedious to derive which is mainly
due to the parity conditions. But there is a trick to derive our
formulas directly from Macdonald's formulas (6.11) and ($6.11'$). For
any product of the form
$$\eqno{}
P=\prod_i(1-q^{a_i}t^{b_i}).
$$
with $a_i,b_i\ge0$ we define
$$\eqno{}
[P]_\even:=\prod_{\{i\mid b_i\ \even\}}(a_i+b_i\r).
$$
The map $P\mapsto[P]_\even$ is multiplicative. The point is now the
easily verified formula
$[B_{\nu/\mu}]_\even=w_{(1^b)}(\rho_\a+\mu;-r)$ with $B_{\nu/\mu}$ as
in {\it loc.~cit.}~VI(6.4). From {\it loc.~cit.}~VI(6.10) we obtain
$B_\mu=[u_0(P_\mu)]_\even$. Now the formulas for $B_\mu$ above are
nothing else than $[\cdot]_\even$ applied to {\it loc.~cit.}~VI(6.11)
and ($6.11'$).\qed

\Remarks: 1. The evaluation formula \cite{E20} are also valid for
$r=0$ provided one replaces the expressions for $B_\lambda$ in
\cite{E84} by their limits for $r\pfeil0$. The same remark holds for
all the Pieri formulas in section~\cite{PieriFormula}.
\Aussage{Conjecture}

2. The polynomials $R_\lambda(z;r)$ have, in general, coefficients in
$\QQ(r)$. Conjecturally, one obtains an integral form as follows. For
$\lambda\in\Lambda$ let
$$\eqno{}
[c_\lambda]_\even:=\prod\limits_{s\in\lambda\atop l_\lambda(s)\ \odd}
\big(a_\lambda(s)+(l_\lambda(s)+1)\r\big)
$$
and $\cR_\lambda(z;r):=[c_\lambda]_\even\, R_\lambda(z;r)$.

\Conjecture1. For all $\lambda\in\Lambda$ holds
$\cR_\lambda(z;r)\in\ZZ[\r,z]$.

\noindent The factor $[c_\lambda]_\even$ is the denominator of the
second expression for $B_\lambda$ in \cite{E84}. For $n>\!\!>0$ there
is no cancellation involving the variable $r$. Thus, at least up to a
rational factor and for big $n$, the conjectured statement appears to
be optimal. The conjecture has been tested for $n\le 6$ and
$|\lambda|_\odd\le 6$. (Shifted) Jack polynomials have also certain positivity
properties (see \cite{KnSa} and \cite{KnSa2}) but none of them seem to
generalize to semisymmetric polynomials.

\medskip
Now, we specialize the evaluation and the binomial formula to the homogeneous
polynomials $\Rq_\lambda$. The evaluation formula \cite{E20} becomes

\Theorem. For $\lambda\in\Lambda$ holds
$$\eqno{}
\Rq_\lambda(1,\ldots,1)=\cases{
B_\lambda&if $n$ is odd or $[\lambda]_1=0$;\cr
0&otherwise.\cr}
$$

\Proof: This follows from \cite{E20} by calculating
$\|lim|\limits_{\a\pfeil\infty}\a^{-|\lambda|_\odd}R_\lambda(-\rho_\a)$.
The parity of $n$ comes in because $\|deg|_\a A_\lambda(\a)$ equals
$|\lambda|_\odd$ (resp. $|\lambda|_\even$) if $n$ is odd
(resp. even).\qed

Now, the binomial formula \cite{E1} becomes:

\Theorem. Assume $r\ne0$. For all $\lambda\in\Lambda$ holds
$$\eqno{E31}
{\Rq_\lambda(\underline{1}+z)\over B_\lambda}
=\sum_\mu{R_\mu(\rho+\lambda)\over R_\mu(\rho+\mu)}
{\Rq_\mu(z)\over B_\mu}
$$
The sum is over all $\mu\in\Lambda$ with $\mu\LE\lambda$ and, in
case $n$ is even, additionally with
$[\mu]_1=[\lambda]_1$. 

\Proof: Using the evaluation formula \cite{E20}, this follows from
\cite{E1} by replacing $z$ by $\a z$, dividing by the appropriate
power $\a^N$ and taking the limit $\a\pfeil\infty$.  Here
$$\eqno{}
N=\cases{
|\lambda|_\odd-|\lambda|_\odd=0&if $n$ is odd;\cr
|\lambda|_\odd-|\lambda|_\even=[\lambda]_1&if $n$ is even.\cr}
$$\qed

\Remark: In the classical case, Lassalle~\cite{La} used the expansion
of $\Pq_\lambda(\underline{1}+z)$ to {\it define} ``generalized
binomial coefficients''. Okounkov\_Olshanski~\cite{OO2} proved later
the classical analogue of \cite{E31}. This implies in particular that
Lassalle's binomial coefficients equal ${P_\mu(\rho+\lambda)\over
P_\mu(\rho+\mu)}$. Thus, classically it is possible to define shifted
polynomials from the theory of the unshifted ones.  The theorem above
shows that this fails in the semisymmetric case if the number of
variables $n$ is even.

\beginsection PieriFormula. The Pieri formula

Now, we make the Pieri formula
completely explicit. Recall again some notation: For
$\lambda=\mu+\epsilon_I$ let $C_{\lambda/\mu}$
(resp. $R_{\lambda/\mu}$) be the set of boxes of $\mu$ which are in
the same column (resp. row) as some box of $\lambda\setminus\mu$.

\Theorem Pieri2. Assume $r\ne0$. For every $\mu\in\Lambda$ holds
$$\eqno{}
\prod_{i\ \odd/\even}\!\!\!\!\!(t+z_i)\ R_\mu(z)=
\sum_{I\in P_{\odd/\even}}\ \ \ \ 
\pprod_{i\not\in I\atop i\ \odd/\even}(t+\rho_i+\mu_i)
\ [\psi_{\lambda/\mu}']_\even\ 
R_\lambda(z)
$$
Here $\lambda=\mu+\epsilon_I$ and $I$ runs through $P_{\odd/\even}$.
Moreover,
$$\eqno{}
\eqalign{
[\psi_{\lambda/\mu}']_\even&=
{\prod\limits_{i\not\in I,j\in I\atop i<j,j-i\ \odd}
\big(\mu_i-\mu_j+(j-i-1)\r\big)\big(\lambda_i-\lambda_j+(j-i+1)\r\big)
\over
\prod\limits_{i\not\in I,j\in I\atop i<j,j-i\ \even}
\big(\mu_i-\mu_j+(j-i)\r\big)\big(\lambda_i-\lambda_j+(j-i)\r\big)}=\cr
&=\prod_{s\in C_{\lambda/\mu}-R_{\lambda/\mu}}
{[b_\lambda(s)]_\even\over[b_\mu(s)]_\even}\cr}
$$
where
$$\eqno{}
[b_\lambda(s)]_\even:=\cases{a_\lambda(s)+(l_\lambda(s)+1)\r&\vrule
width 0pt depth 5pt height 0pt for
$l_\lambda(s)$ odd,\cr
\big((a_\lambda(s)+1)+l_\lambda(s)\r\big)^{-1}&for
$l_\lambda(s)$ even.\cr}
$$

\Proof: By \cite{E7} and the evaluation formula \cite{E20} we have to
calculate
$$\eqno{}
\left(v_I(\rho_\a+\mu){A_\mu(\a)\over A_{\mu+\epsilon_I}(\a)}\right)
\left(w_I(\rho_\a+\mu;-r){B_\mu\over B_{\mu+\epsilon_I}}\right).
$$
One easily verifies that the first factor is $1$ (not too surprising,
given the fact that the result can not depend on $\a$). For the second
factor we apply the same trick as in the proof of
\cite{SpecialValue} and obtain: it is $[\cdot]_\even$ applied to the
corresponding formulas for Macdonald polynomials. The result follows
from \cite{Mac}~VI($6.7'$), (6.13), and (6.23).\qed

\Remarks: 1. A priori, it might happen that $\lambda=\mu+\epsilon_I$
is not a partition. But then $[\psi_{\lambda/\mu}']_\even=0$ and the
corresponding summand may be omitted. In fact, in that case there is
an index $i$ such that $i\not\in I$, $j:=i+1\in I$ and
$\mu_i=\mu_{i+1}$. Thus $[\psi_{\lambda/\mu}']_\even$ contains the factor
$\mu_i-\mu_j+(j-i-1)r=0$.

\noindent 2. The first expression for $[\psi_{\lambda/\mu}']_\even$
shows that it is a rational function in $\mu$. The second expression
takes the cancellation into account which occurs when the ``vertical
strip'' $\lambda\setminus\mu$ contains boxes in the same column, i.e.,
if there is an index $i$ such that $i\in I$, $i+1\in I$ and
$\mu_i=\mu_{i+1}$.  \medskip

By comparing coefficients of powers of $t$ we easily obtain
Pieri formulas involving elementary symmetric functions:

\Corollary Pieri3.
$$\eqno{E23}
e_m(z_i\mid i\ {\scriptstyle{\odd\atop\even}})\ R_\mu(z)=
\sum_I\ \ 
e_{m{-}s}(\mu_i{+}\rho_i\mid
i\not\in I,i\ {\scriptstyle{\odd\atop\even}})
\ [\psi_{\lambda/\mu}']_\even\ 
R_\lambda(z)
$$
where $\lambda=\mu+\epsilon_I$ and where $I$ runs through the elements
of $P_{\odd/\even}$ with $s:=|I|_o\le m$.

\Example: For $n=3$ we obtain:
\setbox0=\vtop{\baselineskip20pt\parindent30pt\everymath{\displaystyle}\noindent
$(z_1{+}z_3)\cdot R_{\mu_1,\mu_2,\mu_3}=
(\mu_1{+}\mu_3{+}2r)R_{\mu_1,\mu_2,\mu_3}+$

${+}R_{\mu_1{+}1,\mu_2,\mu_3}+{(\mu_2{-}\mu_3)(\mu_2{-}\mu_3{-}1{+}2r)\over(\mu_1{-}\mu_3{+}2r)(\mu_1{-}\mu_3{-}1{+}2r)}R_{\mu_1,\mu_2,\mu_3{+}1}+$

$+R_{\mu_1{+}1,\mu_2{+}1,\mu_3}+{(\mu_1{-}\mu_2)(\mu_1{-}\mu_2{-}1{+}2r)\over(\mu_1{-}\mu_3{+}2r)(\mu_1{-}\mu_3{-}1{+}2r)}R_{\mu_1,\mu_2{+}1,\mu_3{+}1}$}
$$\eqno{E25}
\box0
$$
\setbox0=\vtop{\baselineskip20pt
\parindent30pt\everymath{\displaystyle}\noindent
$
z_2\cdot R_{\mu_1,\mu_2,\mu_3}=
(\mu_2{+}r)R_{\mu_1,\mu_2,\mu_3}+
$

$
+R_{\mu_1{+}1,\mu_2{+}1,\mu_3}+{(\mu_1{-}\mu_2)(\mu_1{-}\mu_2{-}1{+}2r)\over
(\mu_1{-}\mu_3{+}2r)(\mu_1{-}\mu_3{-}1{+}2r)}
R_{\mu_1,\mu_2{+}1,\mu_3{+}1}
$}
$$\eqno{E26}
\box0
$$
\setbox0=\vtop{\baselineskip20pt\parindent30pt\everymath{\displaystyle}\noindent
$z_1z_3\cdot R_{\mu_1,\mu_2,\mu_3}=
(\mu_1{+}2r)\mu_3R_{\mu_1,\mu_2,\mu_3}+$

${+}\mu_3R_{\mu_1{+}1,\mu_2,\mu_3}+(\mu_1{+}2r){(\mu_2{-}\mu_3)(\mu_2{-}\mu_3{-}1{+}2r)\over(\mu_1{-}\mu_3{+}2r)(\mu_1{-}\mu_3{-}1{+}2r)}R_{\mu_1,\mu_2,\mu_3{+}1}+$

$+\mu_3R_{\mu_1{+}1,\mu_2{+}1,\mu_3}+(\mu_1{+}2r){(\mu_1{-}\mu_2)(\mu_1{-}\mu_2{-}1{+}2r)\over(\mu_1{-}\mu_3{+}2r)(\mu_1{-}\mu_3{-}1{+}2r)}R_{\mu_1,\mu_2{+}1,\mu_3{+}1}+
$

$+R_{\mu_1{+}1,\mu_2{+}1,\mu_3{+}1}$}
$$\eqno{}
\box0
$$

\medskip
Formula \cite{E23} can be used to give a Pieri rule for
multiplication with shifted elementary semisymmetric polynomials. For
this, we introduce the following notation. Let $I\in P_\odd$, $s=|I|$
and $f$ a semi\_symmetric polynomial in $n-s$ variables. Then we
define $f(z|I'):=f(z_{k_{s+1}},\ldots,z_{k_n})$ where $i\mapsto k_i$
is any {\it parity preserving\/} bijection from $\{s+1,\ldots,n\}$ to
$\{1,\ldots,n\}\setminus I$. For example, if $n=5$, then
$f(z|\{2,3,5\}')=f(z_4,z_1)$.

\Theorem Pieri4. For $m=0,\dots,n$ and $\mu\in\Lambda$ holds
$$\eqno{E28}
R_{(1^m)}(z)R_\mu(z)=
\sum_{s=0}^m\sum_{|I|=s}R_{(1^{m-s})}(\rho+\mu|I')
\ [\psi_{\lambda/\mu}']_\even\ R_\lambda(z).
$$
Here $\lambda=\mu+\epsilon_I$ and $I$ runs through $P_{\odd/\even}$
according to $m$ $\odd/\even$.

\Proof: From the Triangularity \cite{Triangular1} (or
\cite{Elementary} and the explicit formulas in \cite{KnSa}) follows that
$R_{(1^m)}$ is a linear combination of elementary symmetric functions in
the odd or even variables. Thus \cite{Pieri3} implies that there are
semisymmetric functions $f_p(z_1,\ldots,z_{n-p})$ such that
$$\eqno{E8}
R_{(1^m)}(z)R_\mu(z)=
\sum_I f_{|I|}(\rho+\mu|I')[\psi'_{\lambda/\mu}]_\even R_\lambda(z).
$$
Here the sum is over all $I\in P_{\odd/\even}$  if $m$ is $\odd/\even$.
Moreover, (with $p=|I|$)
$$\eqno{}
\|deg|f_p\le\|deg|R_{(1^m)}-|I|_o=
\Big\lceil{m\over2}\Big\rceil-\Big\lceil{p\over2}\Big\rceil.
$$
As one easily checks, the formula
$$\eqno{}
\Big\lceil{m\over2}\Big\rceil-\Big\lceil{p\over2}\Big\rceil=
\Big\lceil{m-p\over2}\Big\rceil
$$
holds {\it except\/} when $m$ is even and $p$ is odd, a case
which does not occur. Thus we get $\|deg|f_p\le\|deg|R_{(1^{m-p})}$ for all
$m$.

Next, we show that the vanishing conditions hold for $f_p$. For this,
let $I=\{1,\ldots,p\}$. Then
$f_{|I|}(\rho+\mu|I')=f_p(\rho_{p+1}+\mu_{p+1},\ldots,\rho_n+\mu_n)$.
Put $\lambda=\mu+\epsilon_I$. Since $[\psi'_{\lambda/\mu}]_\even\ne0$
it suffices to show that $R_\lambda$ does not occur in the expansion
of $R_{(1^m)}R_\mu$ whenever $\mu_{p+(m-p)}=\mu_m=0$. For this, we put
$z_m=\rho_m,\ldots,z_n=\rho_n$. Then the left hand side of \cite{E8}
vanishes while, by \cite{Stability}, on the right hand side those
$R_\lambda$'s which don't vanish remain linearly independent. Thus, the
coefficient in front of $R_{\mu+\epsilon_I}$ is zero which proves the
claim.

We have proved that $f_p$ is a multiple of $R_{(1^{m-p})}$. To show
equality we put $z_{m+1}=\rho_{m+1},\ldots,z_n=\rho_n$ in \cite{E8}
and then replace $z_i$ by $\rho_m+z_i$ for $i=1,\ldots,m$. Then
$R_{(1^m)}(z)$ becomes $R_{(1^m)}(z_1,\ldots,z_m)$ which is the last
elementary symmetric polynomial in the odd, respectively even,
variables. Thus \cite{E8} becomes simply a special case of
\cite{Pieri3} which implies $f_p=R_{(1^{m-p})}$.\qed

\Example: By \cite{E24} for $n=3$ holds $R_{(1)}(z)=z_1-z_2+z_3-r$. Thus
formulas \cite{E25} and \cite{E26} imply:
\setbox0=\vtop{\baselineskip20pt\parindent30pt\everymath{\displaystyle}\noindent
$R_{(1)}\cdot R_{\mu_1,\mu_2,\mu_3}=
(\mu_1{-}\mu_2{+}\mu_3)R_{\mu_1,\mu_2,\mu_3}+$

\hfill
${+}R_{\mu_1{+}1,\mu_2,\mu_3}+{(\mu_2{-}\mu_3)(\mu_2{-}\mu_3{-}1{+}2r)\over(\mu_1{-}\mu_3{+}2r)(\mu_1{-}\mu_3{-}1{+}2r)}R_{\mu_1,\mu_2,\mu_3{+}1}$%
\hskip40pt\strut

}
$$\eqno{F2}
\box0
$$
in accordance with \cite{E28}, case $m=1$. Observe also the
cancellation which occurs when one subtracts \cite{E26} from
\cite{E25}. This is reflected in the fact that in \cite{E28} for $m$
odd the a priori possible terms with $|I|=m+1$ are missing.\medskip

As a consequence of \cite{E28} we obtain a Pieri rule for the top
homogeneous parts:

\Corollary Pieri5. For every $\mu\in\Lambda$ and $m=1,\dots,n$ holds
$$\eqno{}
\e_m(z)\Rq_\mu(z)=\sum_I[\psi_{\lambda/\mu}']_\even\Rq_\lambda(z).
$$
Here $\lambda=\mu+\epsilon_I$ and $I$ runs through all subsets of
$\{1,\ldots,n\}$ consisting of $\lceil{m\over2}\rceil$ odd numbers and
$\lfloor{m\over2}\rfloor$ even numbers.

Finally, we complete the explicit computation of $R_\lambda$, started
in  \cite{Hook1}, where $\lambda$ is a hook.

\Corollary Hook2. Let $a,m\ge2$ be integers with $m$ even. Then
$$\eqno{E29}
R_{(a\,1^{m{-}1})}=(R_{(1)}-1)(R_{(1)}-2)\ldots(R_{(1)}-a+2)
\big(R_{(1)}R_{(1^m)}-{a-1\over a-1+m\r}R_{(1^{m{+}1})}\big).
$$

\Proof: From formula \cite{E28} and some short calculations we get
$$\eqno{}
R_{(1)}\cdot R_{(a{-}1\,1^{m{-}1})}=(a-2)R_{(a{-}1\,1^{m{-}1})}+
R_{(a\,1^{m{-}1})}+{mr\over(a{-}2{+}mr)(a{-}1{+}mr)} R_{(a{-}1\,1^{m})}.
$$
Thus,
$$\eqno{}
R_{(a\,1^{m-1})}=(R_{(1)}{-}a{+}2)R_{(a-1\,1^{m-1})}-
{mr\over(a{-}2{+}mr)(a{-}1{+}mr)}  R_{(a{-}1\,1^{m})}.
$$
This already implies formula \cite{E29} for $a=2$. For $a\ge2$ we are
using induction and formula \cite{E30} for $R_{(a{-}1\,1^{m})}$:
$$\eqno{}
\eqalign{
R_{(a\,1^{m-1})}=&(R_{(1)}-1)\ldots(R_{(1)}-a+2)\big(R_{(1)}R_{(1^m)}-{a-2\over
a-2+m\r}R_{(1^{m{+}1})}\big)-\cr
&-{mr\over(a-2+mr)(a-1+mr)}(R_{(1)}-1)\ldots (R_{(1)}-a+2)R_{(1^{m+1})}=\cr
=&(R_{(1)}-1)(R_{(1)}-2)\ldots(R_{(1)}-a+2)\big(R_{(1)}R_{(1^m)}-{a{-}1\over
a{-}1{+}m\r}R_{(1^{m{+}1})}\big).\cr}
$$
\qed

\hfuzz 10pt

\beginrefs

\B|Abk:Bate|Sig:Ba|Au:Bateman, H.|Tit:Higher transcendental functions,
Vol.~I|Reihe:Bateman Manuscript Project|Verlag:McGraw
Hill|Ort:New York|J:1953||

\L|Abk:BenRat|Sig:BR|Au:Benson, C., Ratcliff, G.|Tit:A
classification of multiplicity free actions|Zs:J. Algebra|Bd:181|S:152--186|J:1996||

\L|Abk:Deb|Sig:De|Au:Debiard, A.|Tit:Polyn{\^o}mes de Tch{\'e}bychev et de
Jacobi dans un espace euclidien de dimension $p$|Zs:C.R. Acad. Sc. Paris%
|Bd:296|S:529--532|J:1983||

\B|Abk:GW|Sig:GW|Au:Goodman, R.; Wallach, N.|Tit:Representations and
invariants of the classical groups|Reihe:Encycl. Math. Appl. {\bf
68}|Verlag:Cambridge U. Press|Ort:Cambridge|J:1998||

\L|Abk:Kac|Sig:Kac|Au:Kac, V.|Tit:Some remarks on nilpotent orbits%
|Zs:J. Algebra|Bd:64|S:190--213|J:1980||

\Pr|Abk:Mon|Sig:Kn1|Au:Knop, F.|Artikel:Some remarks on multiplicity
free spaces|Titel:Proc. NATO Adv. Study Inst. on Representation Theory and
Algebraic Geometry|Hgr:A.~Broer, G.~Sabidussi, eds.|Reihe:Nato ASI
Series C|Bd:514|Verlag:Kluwer|Ort:Dortrecht|S:301--317|J:1998||

\L|Abk:Cap|Sig:Kn2|Au:Knop, F.|Tit:Symmetric and non-symmetric quantum
Capelli polynomials|Zs:Comment. Math. Helv.|Bd:72|S:84-100
{\tt q-alg/9603028}|J:1997||

\L|Abk:KnSa|Sig:KS1|Au:Knop, F.; Sahi, S.|Tit:Difference equations and
symmetric polynomials defined by their zeros|%
Zs:Intern. Math. Res. Notices|Bd:10|S:473--486 {\tt q-alg/9610017}|J:1996||

\L|Abk:KnSa2|Sig:KS2|Au:Knop, F.; Sahi, S.|Tit:A recursion and a
combinatorial formula for Jack polynomials|Zs:Invent.
Math.|Bd:128|S:9--22 {\tt q-alg/9610016}|J:1997||

\L|Abk:La|Sig:La|Au:Lassalle, M.|Tit:Une formule du bin\^ome
g\'en\'eralis\'ee pour les polyn\^omes de Jack|Zs:C. R. Acad. Sci.
Paris S\'er. I Math.|Bd:310|S:253--256|J:1990||

\L|Abk:Leahy|Sig:Le|Au:Leahy, A.|Tit:A classification of
multiplicity free representations, Ph. D. thesis, Rutgers Univ.,
1995|Zs:J. Lie Theory|Bd:8|S:see also 
http://www.emis.de/journals/JLT/ index.html|J:1998||

\B|Abk:Mac|Sig:Mac|Au:Macdonald, I.|Tit:Symmetric functions and Hall
polynomials, 2nd edition|Reihe:Oxford Mathematical
Monographs|Verlag:Clarendon Press|Ort:Oxford|J:1995||

\L|Abk:Ok|Sig:Ok|Au:Okounkov, A.|Tit:Binomial formula for Macdonald
polynomials and applications|Zs:Math. Res. Lett.|Bd:4|S:533--553
{\tt q-alg/9608027}|J:1997||

\L|Abk:OO2|Sig:OO|Au:Olshanski, G.; Okounkov, A.|Tit:Shifted Jack
polynomials, binomial formula, and applications|Zs:Math. Res. Lett.
|Bd:4|S:69--78|J:1997||

\L|Abk:Seki|Sig:Se|Au:Sekiguchi, J.|Tit:Zonal spherical functions on
some symmetric spaces|Zs:Publ. RIMS, Kyoto Univ.|Bd:12|S:455--459|J:1977||

\L|Abk:St|Sig:St|Au:Stanley, R.|Tit:Some combinatorial properties of
Jack symmetric functions|Zs:Adv. Math.|Bd:77|S:76--115|J:1989||

\B|Abk:VK|Sig:VK|Au:Vilenkin, N.; Klimyk, A.|Tit:Representation of
Lie groups and special functions,
Vol.~2|Reihe:-|Verlag:Kluwer|Ort:Dordrecht|J:1993||

\L|Abk:VS|Sig:VS|Au:Vilenkin, N.; \v Sapiro, R.|Tit:Irreducible
representations of the group ${\rm SU}(n)$ of class I relative
to ${\rm SU}(n-1)$|Zs:Izv. Vys\v s. U\v cebn. Zaved.
Matematika|Bd:62|S:9--20|J:1967||

\endrefs

\bye
\endinput

\L|Abk:|Sig:|Au:|Tit:|Zs:|Bd:|S:|J:||
\B|Abk:|Sig:|Au:|Tit:|Reihe:|Verlag:|Ort:|J:||
\Pr|Abk:|Sig:|Au:|Artikel:|Titel:|Hgr:|Reihe:|Bd:|Verlag:|Ort:|S:|J:||
\bye